\documentclass[11pt, letterpaper]{article} 

\usepackage[margin=1.1in]{geometry}
\usepackage{parskip}
\usepackage{amsmath, amsfonts, amssymb, amsthm, amssymb}
\usepackage{mathrsfs}
\usepackage{bbm}
\usepackage{bm}
\usepackage{dsfont} 

\usepackage{longtable, booktabs, array}
\usepackage{graphicx, caption, subcaption, epstopdf}
\usepackage{hhline, color, xcolor}

\usepackage{algorithm, algorithmic}
\usepackage{paralist, enumitem}
\usepackage{verbatim, appendix, eurosym}

\usepackage{tikz}
\usetikzlibrary{arrows, shadows, positioning}
\tikzset{
  frame/.style={rectangle, draw, text width=6em, text centered, minimum height=4em, drop shadow, fill=white, rounded corners},
  line/.style={draw, -latex', rounded corners=3mm}
}

\usepackage[numbers, sort&compress]{natbib}
\definecolor{linkcolor}{rgb}{0.1,0,0.7}
\usepackage[unicode=true, colorlinks, citecolor=linkcolor, linkcolor=linkcolor, urlcolor=blue]{hyperref}
\usepackage{cleveref}

\newtheorem{Theorem}{Theorem}[section]
\newtheorem{Definition}[Theorem]{Definition}
\newtheorem{Proposition}[Theorem]{Proposition}
\newtheorem{lemma}[Theorem]{Lemma}

\newtheorem{Remark}[Theorem]{Remark}

\newtheorem{Assumption}[Theorem]{Assumption}

\numberwithin{equation}{section}

\newcommand{\nc}{\newcommand}
\nc{\ind}{\mathds{1}}
\def \trans{^{\scriptscriptstyle{\intercal}}}

\newcommand{\R}{\mathbb{R}}

\newcommand{\E}{\mathbb{E}}
\newcommand{\F}{\mathcal{F}}

\renewcommand{\P}{\mathbb{P}}

\newcommand{\Pc}{\mathcal{P}}

\newcommand{\dd}{\mathrm{d}}
\newcommand{\sF}{\mathscr{F}}

\allowdisplaybreaks[4]

\DeclareMathOperator{\esssup}{esssup}
\def\esssup_#1{\underset{#1}{\mathrm{ess\,sup\, }}}
\def\essinf_#1{\underset{#1}{\mathrm{ess\,inf\, }}}
\def\argmax_#1{\underset{#1}{\mathrm{arg\,max\, }}}
\def\argmin_#1{\underset{#1}{\mathrm{arg\,min\, }}}

\def\b1{\bf 1}

\def \Ac{{\cal A}}  \def \Cc{{\cal C}}
 \def \Ec{{\cal E}} \def \Fc{{\cal F}}
\def \Gc{{\cal G}}  \def \Ic{{\cal I}}
 \def \Lc{{\cal L}} 
 \def \Nc{{\cal N}} 
\def \Sc{{\cal S}}  \def \Uc{{\cal U}}
 \def \Wc{{\cal W}}

\def\bpi{\boldsymbol{\pi}}
\def\bh{{\bm h}}

\begin{document}

\title{Unified continuous-time q-learning for mean-field game and mean-field control problems}

\author{Xiaoli Wei \thanks{Email: xiaoli.wei@hit.edu.cn, Institute for Advanced Study in Mathematics, Harbin Institute of Technology, China.}
\and
Xiang Yu \thanks{Email: xiang.yu@polyu.edu.hk, Department of Applied Mathematics, The Hong Kong Polytechnic University, Kowloon, Hong Kong.}
\and
Fengyi Yuan \thanks{Email: fengyiyuan@cuhk.edu.cn, School of Science and Engineering, the Chinese University of Hong Kong (Shenzhen), Shenzhen, China}
}
\date{\vspace{-0.2in}}

\maketitle

\begin{abstract}
This paper studies the continuous-time q-learning in mean-field jump-diffusion models {in a setting where the environment simulator does not provide direct access to the population distribution}. We propose the integrated q-function in decoupled form (decoupled Iq-function) and establish its martingale characterization, which provides a unified policy evaluation rule for both mean-field game (MFG) and mean-field control (MFC) problems. Moreover, we consider the learning procedure where population distribution is updated based on the representative agent's state values. Depending on the task to solve the MFG or MFC problem, we can employ the decoupled Iq-function differently to characterize the mean-field equilibrium policy or the mean-field optimal policy respectively. Based on these theoretical findings, we devise a unified {parametric} q-learning algorithm for both MFG and MFC problems by utilizing test policies and the averaged martingale orthogonality condition. In two applications within and beyond LQ framework, we illustrate the effectiveness and efficiency of our unified parametric q-learning algorithm for both MFG and MFC learning tasks.
\end{abstract}


\section{Introduction}

Since their pioneer introduction of mean-field game (MFG) in \cite{LL2007} and \cite{Huangetal2006}, the mean-field theory has raised a lot of interests in the past two decades and has become a successful and tractable approximation for the large stochastic system. Two different types of mean-field dynamics, namely the competitive game formulation and the cooperative game formulation, lead to the MFG and mean-field control (MFC) problems with wide applications. To be more precise, a MFG problem is to find a Nash equilibrium in a large population competition game where each representative agent interacts with the population and aims to maximize his own objective function. On the other hand, MFC is concerned with a large population cooperative game where all representative agents, again interacting with the population, share the same goal to attain the social optimum, which is often interpreted as a centralized control problem by a social planner. We refer to \cite{CarD, CarD2} for further details on these two formulations.

\paragraph{{RL for discrete-time MFG and MFC}} Reinforcement learning (RL), one of the most fast growing branches in machine learning, provides a powerful way to learn solutions to MFG and MFC problems in the unknown environment via the repetitive trial-and-error procedure. Most existing RL algorithms for MFG and MFC focused on discrete-time mean-field Markov decision processes. \cite{GHXZ2019} initially proposed a discrete-time Q-learning algorithm for MFG with some convergence guarantees. Later, some extensions and different algorithms have been quickly developed, see among fictitious play algorithms in \cite{Elieetal2020, Perrinetal2020, cuikoeppl2021, Xieetal2020}, online mirror descent algorithms in \cite{Perolatetal2021}, population-dependent policies in \cite{PLPEGP2022}, and regularized MFG in \cite{AKS2023}. The solutions learnt by these RL algorithms provide approximate Nash equilibria for finite-agent games, see e.g. \cite{AKS2023b, Saldi18}. In the meantime, RL has also been employed for solving MFC problems with different features in \cite{CLT23}, \cite{GGWX23, GGWX21a, GGWX21b}, \cite{PB2021}, \cite{CTSK2021}, \cite{Mondaletal2022, Mondaletal2023}, to name a few. Recently, \cite{AFL2022, AFHR2023} studied  unified two-timescale RL algorithms to solve discrete-time MFG and MFC problems in finite and continuous state-action space, respectively. The related convergence result is also established in the subsequent work \cite{AFLZ23}. We refer to \cite{Lauriereetal2022} for a comprehensive survey on some learning MFG and MFC problems.

\paragraph{{Continuous-time RL}} To learn the optimal policy in a continuous-time framework, the conventional way is to first take time discretization and then apply some existing discrete-time RL algorithms. Recently, for single agent's control problems, \cite{WangZhou2020, Wangetal2021, JZ22a, JZ22b, jiazhou2023} laid the theoretical foundations on direct continuous-time reinforcement learning algorithms with entropy regularization. In particular, \cite{jiazhou2023} proposed the continuous-time counterpart of the conventional Q-function and developed a continuous-time q learning theory for single-agent's control problems. The incorporation of entropy regularization in the continuous-time framework, allowing to encourage the exploration, has been further generalized to various RL problems and algorithms, such as linear-quadratic control in \cite{LiLiXu} and \cite{FJ2022}, the optimal stopping problems in \cite{Dong2024}, controlled jump-diffusion models in \cite{GLZ2024} and \cite{BHYZ2026}, MFG problems in \cite{GuoXZ} and \cite{Liangetal2024} and MFC problems in \cite{FGLPS23, weiyu2025, phamwarin2025}. Different financial applications with RL methods have also been actively investigated, see among the optimal tracking portfolio in \cite{Boetal2025}, the time-consistent equilibrium mean-variance portfolio in \cite{DDJ}, and the Merton's expected utility maximization in \cite{DDJZ2023}.

As the first attempt to generalize the continuous-time q-learning in \cite{jiazhou2023} to learn MFC problems, \cite{weiyu2025} investigated the proper definition of the continuous-time integrated q-function (Iq-function) {by assuming that the environment simulator allows full access to the population distribution}. 
The integral form of the Iq-function is crucial in establishing a weak martingale characterization together with the value function based on all test policies. 
However, 
it is also observed in \cite{weiyu2025} that 
the Iq-function in the integral form cannot be used directly to learn the optimal policy. In response, \cite{weiyu2025} proposed another related q-function, called essential q-function, 
and established an integral representation between two q-functions. Thus, a continuous-time q-learning algorithm based on two q-functions (with the same parameters) is devised in \cite{weiyu2025} such that the Iq-function can be learnt by the weak martingale condition and the essential q-function can be used for the policy improvement iterations.

\noindent \paragraph{{Our contributions}} {In sharp contrast to \cite{weiyu2025}, the present paper aims to study the q-learning algorithm for both MFG and MFC in a practical situation where the simulator does not provide the population distribution to each representative agent when they exercise the learnt policy. The contributions of this paper are threefold.}

{
\begin{itemize}
\item  \textbf{Decoupled theoretical framework:} We propose a decoupled exploratory formulation of mean-field dynamics that separates the representative agent's process \eqref{equ:agent-SDE} from the population flow  \eqref{equ:population-SDE}. Within this framework, we introduce the decoupled Iq-function, defined over the product space of both individual and population state-action pairs (see Definition \ref{def:decoupled-q-function}). It integrates the representative agent's dynamics with the Fokker-Planck equation of the population distribution. This definition is fundamentally different from existing ones in the literature, including the single-agent q-function in \cite{jiazhou2023}, the Iq-function in \cite{weiyu2025}, the discrete-time IQ-function in \cite{CLT23, GGWX23, GGWX21a, GGWX21b}, and the standard discrete-time Q-function in \cite{AFL2022, AFHR2023}. We establish a martingale characterization for this decoupled Iq-function, providing a theoretical foundation for mean-field q-learning.

\item\textbf{Unity of RL theoretical foundation for MFG and MFC:} We prove that the decoupled framework provides a unified lens to study both competitive MFG and cooperative MFC. A key theoretical discovery is that both MFG and MFC share the {\it same} martingale condition for policy evaluation. The fundamental distinction between the two problems is reduced to their respective consistency conditions: the mean-field equilibrium (MEF) policy is related to the decoupled Iq-function via a fixed point of Gibbs measure (see \eqref{MFG-PI-map}), while the mean-field optimal (MFO) policy is characterized through an integral representation with an essential q-function (see \eqref{relation-q-essentialq}). This unity allows us to treat different mean-field tasks within a single learning mechanism.

\item\textbf{Unified parametric q-learning algorithm:} Based on the shared martingale condition, we develop a unified q-learning algorithm that does not require the environment simulator to provide population distribution information. Instead, we estimate population distribution flow empirically via observed representative agent's state samples and observations (see \eqref{equ:distribution-update} and \eqref{equ:distribution-update1}). By designing function approximators that satisfy problem-specific consistency constraints, our approach avoids the two-timescale schemes required in previous discrete-time methods \cite{AFHR2023}. We illustrate the algorithm's performance through applications: from the analytic mean-variance portfolio optimization (where MFE and MFO coincide) to a non-LQ crowd-aversion transport problems with jump-diffusion.
\end{itemize}}
{
\noindent \paragraph{Related studies}
Existing continuous-time RL works typically address either MFG or MFC separately. For instance, \cite{GuoXZ} investigates the effect of entropy regularizers on MFG problem, \cite{Liangetal2024} proposes an actor-critic algorithm specifically for MFG, and Frikha et al. \cite{FGLPS23} develops a policy gradient method for MFC.

The most relevant work is the unified discrete-time framework proposed by Angiuli et al. \cite{AFL2022}. While we share the common objective of designing a unified algorithm without a population distribution simulator, the methodologies are fundamentally different. The framework in \cite{AFL2022} relies on a two-timescale scheme: for MFG, the Q-function is updated at a faster rate than the population distribution, while for MFC, the roles are reversed. In contrast, our work extends the single-agent theory of Jia and Zhou \cite{jiazhou2023} to a mean-field setting via a unified decoupled formulation. By introducing a decoupled Iq-function that integrates the representative agent's dynamics with the Fokker-Planck equation of the population, we establish a continuous-time learning object analogous in spirit to the master equation. Consequently, our framework eliminates the need for alternating update rates with different timescale.
}

The rest of the paper is organized as follows. Section \ref{sec:form} introduces the exploratory formulation of the mean-field model in decoupled form. Section \ref{sec:q-func} defines the continuous-time Iq-function in decoupled form and establishes its martingale characterization. Sections \ref{sec:mfg-mfc} provides the theoretical characterizations of the MFE policy for MFG and the MFO policy for MFC, respectively. Building on these foundations, Section \ref{sec:decoupled-q-algo} presents the unified q-learning algorithm for both MFG and MFC problems using the test policies and averaged martingale orthogonality condition. Section \ref{sec:example} further studies several applications in jump diffusion models and numerically illustrates our algorithm with satisfactory performance for both MFG and MFC problems. Section \ref{sec:proofs} collects all the proofs, and Section \ref{sec:conc} concludes the paper.

\section{Formulation of Exploratory Decoupled Mean-Field Problem}\label{sec:form}

\subsection{Decoupled finite-agent formulation}
{We begin with a finite-agent formulation and show that the proposed decoupled model corresponds to its mean-field limit (see e.g., \cite{lacker17}). Consider a system of $N$ agents with mean-field interaction. In view of the homogeneity, we fix an agent $i \in \{1,  \ldots, N\}$. The agent $i$ uses a feedback policy $\hat{\bm \pi}$, while all other agents use a common feedback policy ${\bm \pi}$. The state process of agent $i$ evolves as
\begin{equation*}
\dd X_s^{t, x^i, \hat{\bm \pi}} = b(s, X_s^{t, i, \hat{\bm \pi}}, \mu^{-i, {\bm \pi}}_s, a_s^{i, \hat{\bm \pi}}) \dd s + \sigma(s, X_s^{t, i, \hat{\bm \pi}}, \mu^{-i, {\bm \pi}}_s, a_s^{i, \hat{\bm \pi}}) \dd W_s^i + \Gamma(s, X_s^{t, i, {\bm \pi}},\mu^{-i, {\bm \pi}}_s, a_s^{i, \hat{\bm \pi}}) \dd{\widetilde N_s^i},
\end{equation*}
where $X_t^{t, x^i, \hat{\bm \pi}} = x^i$ is the agent $i$'s initial state, $W^i$, $1 \leq i \leq N$ are independent Brownian motions, $\widetilde N^i$, $1 \leq i \leq N$ are independent Poisson random measures,
and $\mu^{-i, {\bm \pi}}_s = \frac{1}{N-1}\sum_{j \neq i} \delta_{X_s^{j, {\bm \pi}}}$ denotes the empirical distribution of the other agents under policy $\hat{\bm \pi}$.  We write $\mu^{-i}:=\mu_t^{-i,{\bm\pi}}$ for its value at the initial time $t$. The objective function of agent $i$ is
\begin{equation*}
J_d(t, x^i, \mu^{-i}; \hat{\bm \pi}, {\bm \pi}) = \E\left[\int_t^T e^{-\beta(s-t)} r(s,X_s^{t, x^i, \hat{\bm \pi}} , \mu_s^{-i, \bm \pi}, a_s^{i, \bm {\widehat \pi}}) \dd s + e^{-\beta (T -t)}g(X_T^{t, x^i, \hat{\bm \pi}}, \mu_T^{-i, \bm \pi})\right].
\end{equation*}
This objective function can be used to define Nash equilibria or social optima, but we will not pursue these finite-$N$ notions here. Instead, we focus on their mean-field counterparts introduced below.
\begin{Remark}[Information structure]
In this finite-agent formulation, we allow feedback policies that depend on all information available to the agents, including their own state process and aggregate information about the population $\mu_s^{-i, {\bm \pi}}$.
\end{Remark}

 As $N \to \infty$, the empirical measure $\{\mu_s^{-i, {\bm \pi}}\}_{s \geq t}$ is expected to converge to a deterministic flow of probability measures describing the population distribution, see e.g., \cite{lacker17}. This motivates the mean-field formulation introduced in the next subsection.
}

\subsection{Decoupled mean-field formulation}
{ In the limiting model, let us introduce a decoupled formulation by separating the population dynamics from that of the representative agent, allowing them to adopt different policies. This formulation provides a unified framework that accommodates both MFG and MFC problems. The corresponding formulations will be introduced in Sections \ref{sec:mfg} and \ref{sec:mfc}, respectively.
}

Let $(\Omega, \Fc, \P)$ be a complete probability space that supports a $n$-dimensional Brownian motion $W=(W_s)_{s \geq 0}$ and an independent compound Poisson process $N$ with a finite time-dependent intensity $\eta(s)$.  We denote by $\F^{W, N} = (\Fc_s^{W, N})_{s \geq 0}$ the $\P$-completion filtration generated by $W$ and $N$. Assume that there exists a sub-$\sigma$-algebra $\Gc$ of $\Fc$ such that $\Gc$ and $\Fc_\infty^{W, N}$ are independent and $\Gc$ is ``rich enough". We denote by $\F = (\Fc_t)_{t \geq 0}$ the filtration defined as $\Fc_t = \Gc \vee \Fc_t^{B, N}$.

To explore the environment in continuous-time reinforcement learning, we inject a collection of pairwise independent exploration noise $\{U_t\}_{0 \leq t \leq T}$, independent of $W$ and $N$, to continuously sample actions. But this causes some measurability issues as $(t, \omega)$ $\mapsto$ $U_t(\omega)$ is not joint measurable in the usual probability space $(\Omega, \Fc, \P)$, see Remark 2.1 in \cite{STZ24}.  To deal with these measure-theoretical difficulties, we use the Fubini extension in \cite{Sun2006} such that there exists a probability space $(\hat\Omega, \hat\Fc, \hat\P)$ that is rich enough to support a {\it continuum} of essentially pairwise independent random variables $\{U_t\}_{0 \leq t \leq T}$. The probability space $(\hat\Omega, \hat\Fc, \hat\P)$ represents the randomization of actions for the purpose of exploration during the RL procedure. We then expand $(\Omega, \Fc, \F, \P)$  to $(\Omega^e, \Fc^e, \F^e, \P^e)$, where $\Omega^e = \Omega \times \hat\Omega$, $\P^e = \P \otimes \hat\P$, $\F^e = (\Fc^e_t)_{0 \leq t\leq T}$ with $\Fc_t^e = \Fc_t \vee \sigma(U_t, 0 \leq s \leq t)$. We denote by $\E^e$ (resp. $\hat\E$) the expectation under $\P^e$ (resp. $\hat\P$). For any random variable $X$ defined on $(\Omega^e, \Fc^e, \P^e)$, we denote by $\P^e_{X}$ and $\E^e[X] = \E[\hat\E[X]]$  the probability distribution and the expectation of the random variable $X$.
{
\begin{Remark}
Continuous-time RL involves sampling actions from the policy. One way is to use Fubini extension to construct a continuum of pairwise independent random variables, while an alternative approach is to work with discretely sampled processes for the approximation, see e.g. \cite{jiazhou2023}. In order to avoid additional technical details, we adopt the Fubini extension framework in this paper. We note that this it is also widely used in the graphon literature; see, e.g., \cite{Aurelletal2022, Caretal2022}.
\end{Remark}
}
We consider a pair of SDEs, namely the controlled dynamics of the population and that of the representative agent. All agents in the whole population take a feedback relaxed policy ${\bm \pi}$.  The representative agent is allowed to either conform to or deviate from the policy of the population. {Let $\Pi$ stand for the set of admissible policies satisfying the following definition.

\begin{Definition} \label{def:pi}A policy ${\bm \pi}$ is called admissible if
\begin{enumerate}[label=\upshape(\roman*)]
\item ${\bm \pi}(\cdot|t, x, \mu) \in \Pc(\Ac)$, $\mbox{supp}  {\bm \pi}(\cdot|t, x, \mu) = \Ac$ for every $(t, x, \mu) \in [0, T] \times \R^d \times \Pc_2(\R^d)$, and ${\bm \pi}$ is jointly measurable with respect to $(t, x, \mu) \in [0, T] \times \R^d \times \Pc_2(\R^d)$.
\item There exits a constant $C>0$ independent of $(t, a)$ such that for any $\bm \pi \in \Pi$ and any $x, x' \in \R^d$ and $\mu, \mu' \in \Pc_2(\R^d)$
\begin{equation*}
\int_{\Ac} |{\bm \pi}(a|t, x, \mu) - {\bm \pi}(a|t, x', \mu')|\dd a \leq C\left(|x - x'| + \Wc_2(\mu, \mu')\right),
\end{equation*}
  where $\Wc_2$ denotes the 2-Wasserstein metric on $\mathcal{P}_2(\R^d)$.
\item  There exists a constant $C> 0$ such that the entropy of $\bm \pi \in \Pi$ have polynomial growth in $x$ and $\mu$
\begin{equation*}
\big|\mathcal{E}_{\bm \pi}(t, x, \mu)\big| \leq C\left( 1 + |x|^2 +M_2(\mu)^2\right).\\
\end{equation*}
where $M_2(\mu) = (\int_{\R^d} |x|^2 \mu(\dd x))^{1/2}$, and the Shannon entropy $\mathcal{E}_{\bm \pi}$ is defined by: $\mathcal{E}_{\bm \pi}(t, x, \mu) = - \int_{\Ac} \log {\bm \pi}(a|t, x, \mu) {\bm \pi}(a|t, x, \mu)\dd a$.
\end{enumerate}
\end{Definition}

\begin{Remark}[Information structure]
In the mean-field setting, the representative agent's policy at time $t$ is allowed to depend on her own state process as well as the current population distribution (i.e. a deterministic flow). Under this full-information structure, the decoupled formulation problem is fully observed on the extended state space $(t, x, \mu) \in [0, T] \times \R^d \times \Pc_2(\R^d)$.
\end{Remark}

}

{We impose the following assumptions on model coefficients throughout the paper.
\begin{Assumption} \label{ass:b-sigma}\begin{enumerate}[label=\upshape(\roman*)]
\item $b$, $\sigma$, $\Gamma$ and $r$ are jointly continuous in $(t, x, \mu, a) \in [0, T] \times \R^d \times \Pc_2(\R^d) \times \Ac$, and $g$ is jointly continuous in $(x, \mu) \in \R^d \times \Pc_2(\R^d)$.
\item There exists a constant $C >0$ such that for $f \in \{b, \sigma, \Gamma\}$, and all $t, t' \in [0, T]$, $x, x' \in \R^d$, $\mu, \mu' \in \Pc_2(\R^d)$, $a \in \Ac$, it holds that
    \begin{equation*}
    \begin{aligned}
    &|f(t, x, \mu, a) -f(t', x', \mu', a)| \leq C\left(|t- t'| +|x - x'| + \Wc_2(\mu, \mu')\right),\\
    &|r(t, x, \mu, a) - r(t', x', \mu, a)| \leq C\left(1 + |x| + |x'| + \|\mu\|_2 + \|\mu'\|_2\big) \big(|t- t'| + |x - x'| + \Wc_2(\mu, \mu')\right).
    \end{aligned}
    \end{equation*}

\item There exists a constant $C > 0$ such that for each $(t, x, \mu, a) \in [0, T] \times \R^d \times \Pc_2(\R^d) \times \Ac$, it holds that
\begin{equation*}
|r(t, x, \mu, a)| + |g(x, \mu)| \leq C\left(1 + |x|^2 + M_2(\mu)^2 + |a|^2\right).
\end{equation*}

\end{enumerate}
\end{Assumption}
}

For a fixed policy ${\bm \pi} \in \Pi$, the dynamics of the population is a controlled McKean-Vlasov SDE
\begin{equation}\label{equ:population-SDE}
\dd X_s^{{\bm \pi}} = b(s, X_s^{{\bm \pi}}, \P^e_{X_s^{\bm \pi}},  a_s^{{\bm \pi}}) \dd s + \sigma(s, X_{s}^{{\bm \pi}}, \P^e_{X_s^{\bm \pi}}, a_s^{{\bm \pi}}) \dd W_s + \Gamma(s, X_{s-}^{{\bm \pi}}, \P^e_{X_{s-}^{\bm \pi}}, a_{s-}^{\bm \pi}) d{\widetilde N_s}, 
\end{equation}
with $X_{t-}^{{\bm \pi}} = \xi \sim \mu$, $a_s^{\bm \pi} \sim {\bm \pi}(\cdot|s, X_s^{{\bm \pi}}, \P^e_{X_s^{\bm \pi}})$,
where $d{\widetilde N}_t = dN_t - \eta(t)dt$ is a compensated Poisson random measure that is independent of $W$. {Under Assumption \ref{ass:b-sigma} and Definition \ref{def:pi}, the SDE \eqref{equ:population-SDE} has a strong solution denoted by $X^{t, \xi, \bm \pi} = \{X_s^{t, \xi, \bm \pi}, t \leq s \leq T\}$}. We will denote by ${\bm \mu}^{\bm \pi} = (\mu_s^{\bm \pi})_{t \leq s \leq T}$ the corresponding flow of the controlled state process in \eqref{equ:population-SDE}, i.e. $\mu^{\bm \pi}_t = \mu$ and ${\mu}^{\bm \pi}_s = \P^e_{X_s^{t, \xi, \bm \pi}}, s \geq t$. 

Given the population distribution flow ${\bm \mu}^{\bm \pi}$ in \eqref{equ:population-SDE}, the dynamics of the representative agent who takes the action randomly sampled from a policy $\widehat{\bm \pi} \in \Pi$ is then given by
\begin{equation}\label{equ:agent-SDE}
\dd X_s^{t, x, \xi, \widehat{\bm \pi}} = b(s, X_s^{t, x, \xi, \widehat{\bm \pi}}, \mu_s^{\bm \pi}, a_s^{\widehat {\bm \pi}})\dd s + \sigma(s, X_s^{t, x, \mu, \widehat{\bm \pi}}, \mu_s^{\bm \pi}, a_s^{\widehat{\bm \pi}}) \dd W_s + \Gamma(s, X_{s-}^{t, x, \mu, \widehat{\bm\pi}}, \mu_s^{\bm \pi}, a_{s-}^{\widehat{\bm \pi}}) \dd{\widetilde N_s}, 
\end{equation}
with $X_{t-}^{t, x, \xi, \widehat{\bm \pi}} = x \in \R^d$. It is shown that {under Assumption \ref{ass:b-sigma}}, we have $\E^e[\sup_{t \leq s \leq T} |X_s^{t, x, \xi, \widehat{\bm \pi}} - X_s^{t, x, \xi', \widehat{\bm \pi}}|^2] \leq C \Wc_2(\P^e_{\xi}, \P^e_{\xi'})$. In other words, $X_s^{t, x, \xi, \widehat{\bm \pi}}$ depends on the random variable $\xi$ only via its law $\mu$, see \cite{HL16, Li18}. Hence we denote by $X_s^{t, x, \mu, \widehat{\bm \pi}}$ the solution to \eqref{equ:agent-SDE}. 

{
\begin{Remark}
It is observed that the population dynamics \eqref{equ:population-SDE} is a McKean-Vlasov SDE. The evolution of its law induces a nonlinear measure-flow equations or nonlinear controlled Fokker-Planck equations; see e.g. \cite{Bar2023, Dau2023, BFY2013, HamJet2025}. It is worth mentioning that the representative agent's dynamics \eqref{equ:agent-SDE} is a standard controlled SDE, where the measure flow is taken from \eqref{equ:population-SDE} and treated as exogenous rather than being self-generated. Therefore, \eqref{equ:agent-SDE} is referred to as a decoupled SDE in the sense that the mean-field flow is frozen, and this setup is called a decoupled formulation.
\end{Remark}
}

To facilitate the future analysis, it is convenient to consider the exploratory average formulation of \eqref{equ:population-SDE} and \eqref{equ:agent-SDE}:
 \begin{align}
 &\left\{
\begin{aligned}\label{equ:population-average-SDE}
\dd X_s^{t, \xi, {\bm \pi}} = & b_{{\bm \pi}}(s, X_s^{t, \xi, {\bm \pi}}, \mu_s^{\bm \pi})\dd s + \sigma_{{\bm \pi}}(s, X_s^{t, \xi, {\bm \pi}}, \mu_s^{\bm \pi}) \dd W_s \\
&+\int_\Ac\Gamma(s, X_{s-}^{t,\xi, {\bm \pi}},  \mu_{s}^{\bm \pi},a) \widetilde{\Nc}_{\bm \pi}(\dd s,\dd a),\\
X_{t-}^{t,\xi, {\bm \pi}}\sim& \mu,
\end{aligned}\right.
\\
 &\left\{
\begin{aligned}\label{equ:agent-average-SDE}
\dd X_s^{t, x, \mu, \widehat{\bm \pi}} =&  b_{\widehat {\bm \pi}}(s, X_s^{t, x,\mu, \widehat{\bm \pi}},  \mu_s^{\bm \pi})\dd s + \sigma_{\widehat {\bm \pi}}(s, X_s^{t, x, \mu, \widehat{\bm \pi}},  \mu_s^{\bm \pi}) \dd W_s\\
&+\int_\Ac \Gamma(s, X_{s-}^{t, x, \mu, \widehat{\bm \pi}},  \mu_s^{\bm \pi}, a)  \widetilde{\Nc }_{\widehat{\bm \pi}}(\dd s,\dd a),\\
X_{t-}^{t,x,\mu, \widehat{\bm \pi}}=&x,
\end{aligned}
\right.
\end{align}
where for every policy ${\bm h} \in \Pi$, we denote by $\widetilde{\Nc}_{\bm h}(\dd s,\dd a)$ a compensated random measure with compensator ${\bm h}(a|s, x, \mu)\eta(s)\dd a\dd s$, and $b_{{\bm h}}$ and $\sigma_{{\bm h}}$ are defined by
\begin{equation*}
\begin{aligned}
 b_{{\bm h}}(s, x, \mu) := \int_{\Ac} b(s, x, \mu, a){\bm h}(a|s,x,\mu)\dd a, \;\; \sigma_{{\bm h}}(s, x, \mu) := \left(\int_{\Ac} \sigma\sigma\trans(s, x,\mu, a) {\bm h}(a|s, x, \mu) \dd a\right)^{1/2}.
\end{aligned}
\end{equation*}

 By the uniqueness of the controlled martingale problem, the exploratory formulation \eqref{equ:population-SDE}-\eqref{equ:agent-SDE} and exploratory average formulation \eqref{equ:population-average-SDE}-\eqref{equ:agent-average-SDE} correspond to the same controlled martingale problem and thus admit the same solution in law.

 \begin{Remark}\label{rmk:martingale-problem}
 \eqref{equ:population-average-SDE} is called the exploratory average formulation of \eqref{equ:population-SDE}. The martingale formulation of \eqref{equ:population-average-SDE} is described as:
 \begin{equation*}
 f(X_s^{\bm \pi}) - f(\xi) - \int_t^s \int_{\Ac} \mathcal{L}^{u, a, \P^e_{X_u^{\bm \pi}}}[f(\cdot)](X_u^{\bm \pi})\pi(a|u, X_u^{\bm \pi}, \P^e_{X_u^{\bm \pi}})\dd a \dd u, \; s \geq t,
 \end{equation*}
 is a $\P^e$-martingale, where for any bounded and twice continuously differentiable function $f: \R^d \to \R$
\begin{align}\label{generatorL}
\mathcal{L}^{t, a, \mu}[f(\cdot)](x) &= b(t, x, \mu, a)\trans \partial_x f(x) + \frac{1}{2}{\rm Tr}\big(\sigma\sigma\trans(t, x, \mu, a) \partial_{xx} f(x)\big) \\
&\;\;\;+ \big[f(x + \Gamma(t, x, a, \mu)) - f(x) - \Gamma(t, x, a, \mu)\trans \partial_x f(x)\big]\eta(t). \nonumber
\end{align}
It is similar to the relaxed control formulations in \cite{lacker17, BCP20, DPT22}. But their motivations are different: relaxed control formulations in \cite{lacker17,BCP20, DPT22} are mainly used for the existence of optimal controls and the limit theory, whereas the exploratory average formulation here is motivated by the exploration in RL. {We allow for jump-diffusion dynamics in order to incorporate discontinuous shocks in the state process, primarily for modeling generality and for the numerical illustrations in Section~\ref{sec:example}.}
 Although we consider the compound Poisson process here, all arguments and conclusions can be generalized to the Poisson random measure.
 \end{Remark}

{ We consider the entropy-regularized objective functional of the representative agent. This functional represents the expected cumulative reward to be maximized by the representative agent interacting with the population, where the agent adopts the policy $\widehat{\bm \pi}$ and the population evolves according to a fixed policy ${\bm \pi}$. At this stage, the population distribution is treated as exogenous. This functional depends on the agent's private state as well as the initial distribution of the  population, which motivates the term decoupled formulation.} 
\begin{align}\label{decoupled-J}
J_d(t, x, \mu; \widehat{\bm \pi}, {\bm \pi}) =&  \E\biggl[\int_t^T e^{-\beta(s-t)} \left(r(s, X_s^{t, x, \mu, \bm {\widehat\pi}}, \mu_s^{\bm \pi}, a_s^{\bm {\widehat \pi}}) + \gamma \Ec_{\widehat{\bm \pi}}(s, X_s^{t, x, \mu, \bm {\widehat\pi}}, \mu_s^{\bm \pi}) \right)\dd s\\
&\;\;\;\;\;\; + e^{-\beta (T -t)}g( X_T^{t, x, \mu, \bm {\widehat\pi}}, \mu_T^{\bm \pi})\Big| X_t^{t, x, \mu, \bm {\widehat\pi}} = x, \mu_t^{\bm \pi} = \mu\biggl].\nonumber
\end{align}
{ Based on this functional \eqref{decoupled-J}, equilibrium notions will be introduced later in \eqref{MFG-objective} for MFG and \eqref{MFC-objective} for MFC. We incorporate an entropy regularization associated with the representative agent's policy into the objective functional. Such regularization is widely adopted in the learning literature. It ensures that the optimal policy remains non-degenerate distribution over the action space, thereby encouraging sufficient exploration, see e.g., \cite{WangZhou2020, Wangetal2021, JZ22b, GuoXZ, STZ24}.}
In the sequel, we will 
use $X^{P, {\bm \pi}}$ and $X^{R, \widehat{\bm\pi}}$ to represent $X^{t, \xi, {\bm \pi}}$ and $X^{t, x, \mu, \widehat{\bm \pi}}$ for simplicity of notations. Here the superscripts $P$ and $R$ represent the population and the representative agent, respectively. We call $J_d(t, x, \mu; \widehat{\bm \pi}, {\bm \pi})$ the decoupled value function associated with the pair of policies $(\widehat{\bm \pi}, {\bm \pi})$.

\subsection{{Regularity of Decoupled Value Function}}

For the later learning procedure, it is necessary for us to first discuss its regularity under some reasonable assumptions. We introduce the space $\Cc^{2, 2}(\R^d \times \Pc_2(\R^d))$ the set of functions $J: \R^d \times \Pc_2(\R^d) \to \R$ such that $\partial_x J(x, \mu), \partial_{xx} J(x, \mu)$, $\partial_\mu J(x, \mu)(v)$ and $\partial_v \partial_\mu J(x, \mu)(v)$ exist and are jointly continuous.  Furthermore, we need the following two notations: for ${\bm h} \in \Pi$, we also denote $\Gamma^{\eta}_\bh(t,x,\mu) = \eta(t)\int_\Ac \Gamma(t,x,\mu,a)\bh(a|t,x,\mu) \dd a$.
\begin{Assumption}\label{assump} Given ${\bm h}  \in \{\widehat{\bm \pi}, {\bm\pi}\}$, we have
\begin{itemize}
\item[(i)] For $ f\in \{b_{\bm h}, \sigma_{\bm h}, \Gamma^\eta_\bh\}$, the derivatives $\partial_xf(t,x,\mu)$, $\partial_{x}^2 f(t,x,\mu)$, $\partial_{\mu}f(t,x,\mu)(v)$ and $\partial_{v}\partial_{\mu} f(t,x,\mu)(v)$ exist for any $(t,x,v,\mu)\in [0,T]\times  (\R^d)^2 \times \mathcal{P}_2(\mathbb{R}^d)$, are bounded and locally Lipschitz continuous with respect to $x,\mu, v$ uniformly in $t\in[0,T]$. For all $(t,x,\mu)\in[0,T]\times\mathbb{R}^d\times\mathcal{P}_2(\mathbb{R}^d)$, we have $|f(t,x,\mu)|\leq C(1+|x|+M_2(\mu))$.

\item[(ii)] For any $t\in [0,T]$, $r_{\bm h }(t,\cdot)$, $\mathcal{E}_{\bm h}(t,\cdot)$ and $g(\cdot)\in\mathcal{C}^{2,2}(\mathbb{R}^d\times\mathcal{P}_2(\mathbb{R}^d))$.

\item[(iii)] There exists some constant $C<+\infty$, such that for any $(t,x,v,\mu)\in [0,T]\times {\mathbb{R}^d} \times (\R^d)^2 \times \mathcal{P}_2(\mathbb{R}^d)$, we have
\begin{equation*}
\begin{aligned}
&|r_{\bm h}(t,x,\mu)|+|\mathcal{E}_{{\bm h}}(t,x,\mu)|+|g(x,\mu)|\leq C\left(1+|x|^2+M_2(\mu)^2\right),\\
&|\partial_xr_{\bm h}(t,x,\mu)|+|\partial_x\mathcal{E}_{{\bm h}}(t,x,\mu)|+|\partial_xg(x,\mu)|\leq C\left(1+|x|+M_2(\mu)\right),\\
&|\partial_{\mu}r_{\bm h}(t,x,\mu)(v)|+|\partial_{\mu}\mathcal{E}_{{\bm h}}(t,x,\mu)(v)|+|\partial_\mu g(x,\mu)(v)|\leq C\left(1+|x|+|v|+M_2(\mu)\right),\\
&|\partial_{v}\partial_{\mu}r_{\bm h}(t,x,\mu)(v)|+|\partial_x^2r_{\bm h}(t,x,\mu)|+  |\partial_v\partial_{\mu}\mathcal{E}_{{\bm h}}(t,x,\mu)(v)|\\
&\quad\quad+ |\partial_x^2 \mathcal{E}_{{\bm h}}(t,x,\mu)|+ |\partial_v\partial_{\mu}g(x,\mu)(v)|+|\partial_x^2 g(x,\mu)|\leq C\left(1+M_2(\mu)\right).
\end{aligned}
\end{equation*}

\end{itemize}
\end{Assumption}

Further, to impose assumptions on the jump term $\Gamma$, we need some additional notations. Consider a measurable mapping $\Psi:[0,1)\times \Pc_2(\R^d)\to \Pc_2(\R^d)$ such that for any $\mu\in \Pc_2(\R^d)$, $\Psi(\cdot,\mu)_\#({\rm Uniform([0,1])})=\mu$ (c.f. Lemma 5.29 of \cite{CarD}). We fix a reference measure $\mu_0\in \Pc_2(\R)$ with positive density and denote
\[
\Gamma_{\bm h}(t,x,\mu,z):= \Psi \Big(\mu_0((-\infty,z]),\Gamma(t,x,\mu,\cdot)_\#{\bm h}(\cdot|t,x,\mu) \Big).
\]

We note that, if $Z\sim \mu_0$, $\mu_0((-\infty,z])\sim {\rm Uniform([0,1])}$. Consequently, $\Gamma_{\bm h}(t,x,\mu,\cdot)_\# \mu_0=\Gamma(t,x,\mu,\cdot)_\#{\bm h}(\cdot|t,x,\mu)$.

\begin{Assumption}\label{Gassump}
For each $\bh\in \Pi$ and $z\in \R$, the derivatives $\partial_x\Gamma_{\bh}(\cdot,\cdot,\cdot,z)$, $\partial_{x}^2 \Gamma_{\bh}(\cdot,\cdot,\cdot,z)$, $\partial_{\mu}\Gamma_{\bh}(\cdot,\cdot,\cdot,z)(\cdot)$ and $\partial_{v}\partial_{\mu} \Gamma_{\bh}(\cdot,\cdot,\cdot,z)(\cdot)$ exist. Moreover, they are bounded and Lipschitz continuous, with bounds and Lipschitz constants of the form $C|z|$.
\end{Assumption}

Let us introduce the space of $\Cc^{1, 2, 2}([0, T] \times \R^d \times \Pc_2(\R^d))$ as the space of all functions $J: [0, T] \times \R^d \times \Pc_2(\R^d)$ such that $\frac{\partial J}{\partial t}(t, x, \mu)$, $\partial_x J(t, x, \mu), \partial_{xx} J(t, x, \mu)$, $\partial_\mu J(t, x, \mu)(v)$ and $\partial_v \partial J(t, x, \mu)(v)$ are jointly continuous. The following Proposition characterizes the decoupled value function in terms of a parabolic PDE over the enlarged state space $\R^d \times \Pc_2(\R^d)$.

\begin{Proposition}\label{prop:regularity}
Under Assumptions \ref{assump} and \ref{Gassump}, $J_d(\cdot; \widehat{\bm \pi}, {\bm \pi})$ defined in \eqref{decoupled-J} is of $C^{1,2,2}([0,T]\times\mathbb{R}^d\times\mathcal{P}_2(\mathbb{R}^d))$ and satisfies the dynamic programming equation
\begin{align}
&\frac{\partial J_d}{\partial t}(t, x, \mu; \widehat{\bm \pi}, {\bm \pi}) - \beta J_d(t, x, \mu; \widehat {\bm \pi}, {\bm \pi}) + \int_{\Ac} \Big(\mathcal{L}^{t, a, \mu}[J_d(t, \cdot, \mu; \widehat {\bm \pi}, {\bm \pi})](x)\nonumber\\
&+ r(t, x, \mu, a) - \gamma \log \widehat {\bm \pi}(a|t, x, \mu)\Big) \widehat {\bm \pi}(a|t, x, \mu) \dd a \label{DP-equation-J}\\
&+\E^e_{\xi \sim \mu}\Big[\int_{\Ac} \mathcal{L}^{t, a, \mu}[\frac{\delta J_d}{\delta \mu}(t, x, \mu; \widehat {\bm \pi}, {\bm \pi})(\cdot)](\xi){\bm \pi}(a|t, \xi, \mu)\dd a\Big] = 0, \nonumber
\end{align}
with $J_d(T, x, \mu; \widehat{\bm\pi}, {\bm\pi}) = g(x, \mu)$, where $\mathcal{L}^{t, a, \mu}[f](x)$ acting on the $x$ variable of the function $f(x)$ is defined in \eqref{generatorL},
and $\frac{\delta J_d}{\delta \mu}(t, x, \mu; \widehat {\bm \pi}, {\bm \pi})(\cdot)$ stands for the linear functional derivative of $J_d$ with respect to $\mu$.
\end{Proposition}
The proof of Proposition \ref{prop:regularity} is reported in Section \ref{app:proof-prop-regularity}. Define the Hamiltonian operator
\begin{equation*}
\begin{aligned}
H(t, x, \mu, a, p, q, M) & :=  b(t, x, \mu, a) \trans p + \frac{1}{2}{\rm Tr}\big(\sigma\sigma\trans(t, x, \mu, a)q\big) + r(t, x, \mu, a)\\
& + \big[M(x + \Gamma(t, x, a, \mu)) - M(x) - \Gamma(t, x, a, \mu)\trans p\big]\eta(t). \end{aligned}
\end{equation*}
It then follows that
\begin{equation*}
\int_{\Ac}\Big(\mathcal{L}^{t, a, \mu}[J_d(t, \cdot, \mu; \widehat {\bm \pi}, {\bm \pi})](x) + r(t, x, \mu, a)\Big) \widehat {\bm \pi}(a|t, x, \mu) \dd a = H(t, x, \mu, a, \partial_xJ_d, \partial_{xx}J_d, J_d),
\end{equation*}
which comes from the representative agent. The nonlocal term $\E^e_{\xi \sim \mu}[\cdot]$ in \eqref{DP-equation-J} stems from the population, which vanishes when there is no mean-field interaction in the system.

\begin{Remark}
\begin{enumerate}
	\item Note that our state dynamics \eqref{equ:population-average-SDE}-\eqref{equ:agent-average-SDE} are not standard mean-field type SDE with jumps as in \cite{HL16} and \cite{Li18}. Indeed, our jump measure $\tilde N_\bh$ has a compensator which depends on state variables $(t,x,\mu)$. In the mean-field game setting, they are even different between population and representative agent dynamic due to the different choices of policies $\bpi$ and $\hat \bpi$. However, by introducing $\Psi$ and $\Gamma_\bh$ in this section, we transform this nonstandard mean-field system to the one with a consistent jump measure, yet preserve the distributions of state processes. This is achieved by moving all dependences on $(t,x,\mu)$ and $\bh$ to the new jump coefficient $\Gamma_\bh$ and utilizing the relation $\Gamma_\bh(t,x,\mu,\cdot)_\#\mu_0 = \Gamma(t,x,\mu,\cdot)_\#\bh(\cdot|t,x,\mu)$. For details, see Appendix \ref{app:proof-prop-regularity}.
	\item Although \cite{HL16} and \cite{Li18} require some bounds and Lipschitz constants related to $\Gamma$ in the form of $C (1\wedge |z|)$, their results can be easily extended to our case because we have $\mu_0\in \Pc_2(\R)$, a stronger assumption than $\int_\R 1\wedge |z|^2 \mu_0(d z)$ as imposed in \cite{HL16} and \cite{Li18}.
\end{enumerate}
\end{Remark}

\begin{Remark}
As shown in \cite{FGLPS23}, Assumption \ref{assump} is satisfied if $b$ and $\sigma$ are linear functions of $x$ and $\bar{\mu}:=\int_{\R^d} z\mu(d z)$, $r$ and $g$ are quadratic functions of $x$ and $\bar{\mu}:=\int_{\R^d} z\mu(d z)$, and ${\bm \pi}(a|t,x,\mu)$ is Gaussian with its mean and variance as smooth functions of $(t,x,\mu)$. If we further assume $d=1$, previous conditions on model parameters also imply Assumption \ref{Gassump}. Indeed, suppose $\Gamma(t,x,\mu,a)=\Gamma_0(t,x,\mu)a$, and ${\bm \pi}(\cdot|t,x,\mu)=\Nc(m(t,x,\mu),\Sc^2(t,x,\mu))$. Then,
\begin{equation*}
\Gamma(t,x,\mu,\cdot)_\# {\bm \pi}(\cdot|t,x,\mu)=\Nc(\hat{m}(t,x,\mu),\hat{\Sc}^2(t,x,\mu)),
\end{equation*}
where $\hat{m}:=\Gamma_0m$, $\hat{\Sc}=\Gamma_0\Sc$. As $d=1$, the mapping $\Psi$ can be taken as $\Psi(e,\mu)=Q_\mu(e)$, where $Q_\mu$ is the quantile function of $\mu$. In particular, if we take $\mu_0=\Nc(0,1)$, we get
\begin{equation*}
\begin{aligned}
\Gamma_{\bm \pi}(t,x,\mu,z)&=\Psi\Big(\Phi(z),\Nc(\hat{m}(t,x,\mu),\hat{\Sc}^2(t,x,\mu))\Big)
=\hat{m}(t,x,\mu)+\hat{\Sc}(t,x,\mu)z.
\end{aligned}
\end{equation*}
Clearly, Assumption \ref{Gassump} is satisfied if we impose enough assumptions on $m$, $\Sc$ and $\Gamma_0$.
\end{Remark}

\section{Continuous Time Decoupled Iq-Function and Martingale Characterization}\label{sec:q-func}
In this section,  {we start with the definition of a decoupled IQ-function parametrized by a time step $\Delta t >0$ and then motivates the notion of decoupled Iq-function that is independent of $\Delta t$}.
It is revealed that a unified definition of the Iq-function in decoupled form, which will be called \textit{decoupled Iq-function} from here onwards, can be utilized for the learning tasks in both MFG and MFC problems. {We further provide martingale characterization of the decoupled Iq-function.}

\subsection{{Decoupled IQ-function}}
{Recall that in the classical single-agent setting in \cite{jiazhou2023}, the Q-function is interpreted as the value function under a local perturbation of the policy over a short time interval. In contrast to the classical single-agent setting, the representative agent in a mean-field system interacts with the population. It is not sufficient to perturb only the policy of the representative agent as in \cite{jiazhou2023}. One must additionally specify how the population evolves during the perturbation period.}

To this end, we consider a pair of perturbation of policies. Let $\widehat{\bm \pi}^R \in \Pi$ denote the perturbed policy of the representative agent as in \cite{jiazhou2023}: she takes a constant action $a$ on $[t, t+ \Delta t)$ and follows the policy $\widehat {\bm \pi}$ on $[t+\Delta t, T)$; let ${\bm \pi}^P \in \Pi$ be the perturbed policy of the population as follows: the population takes the policy ${\bm h} \in \Pi$ on $[t, t + \Delta t)$ and the policy ${\bm \pi}$ on $[t+\Delta t, T)$.

The state process $X^{P, {\bm \pi}^P}$ on $[t, T)$ of the population is the solution to \eqref{equ:population-average-SDE} under the perturbed policy ${\bm \pi}^P$.
Given the population flow ${\bm \mu}^{{\bm\pi}^P} = (\mu_s^{{\bm\pi}^P})_{t \leq s \leq T}$, the state process $X^{R, \widehat{\bm \pi}^R}$ of the representative agent is described by \eqref{equ:agent-SDE} under the perturbed policy $\widehat{\bm \pi}^R$.
{We then introduce the IQ-function in decoupled form (decoupled IQ-function for short), which  represents the expected total discounted reward of the representative agent under the above perturbed policies $(\widehat{\bm \pi}^R, {\bm \pi}^P)$. For $(t, x, \mu, a, {\bm h})$ $\in$ $[0, T] \times \R^d \times \Pc_2(\R^d) \times \Ac \times \Pi $, define}
\begin{equation*}
\begin{aligned}
& Q_{d, \Delta t} (t, x, \mu, a, {\bm h}; \widehat {\bm \pi}, {\bm \pi}) \nonumber\\
:=& \E^e\biggl[\int_t^{t + \Delta t} e^{-\beta(s-t)}r(s, X_s^{R, \widehat{\bm \pi}^R}, \mu_s^{{\bm \pi}^P}, a)\dd s + \int_{t + \Delta t}^T \int_{\Ac} e^{-\beta(s-t)} \Big(r(s, X_s^{R, \widehat{\bm \pi}^R}, \mu_s^{{\bm \pi}^P}, a^{\widehat{\bm \pi}}_s)\\
 &- \gamma \log \widehat{\bm \pi}(a^{\widehat{\bm \pi}}_s|s, X_s^{R, \widehat{\bm \pi}^R}, \mu_s^{{\bm \pi}^P})\Big) \dd s + e^{-\beta(T-t)} g(X_T^{R, \widehat{\bm \pi}^R}, \mu_T^{{\bm \pi}^P})\biggl].
 \end{aligned}
\end{equation*}
{Note that, unlike the classical Q-function, which depends only the state-action pair $(x, a)$ of the representative agent, the decoupled IQ-function additionally depends on the lifted state-action pair $(\mu, {\bm h})$ of the population.}

By the standard flow property of $\mu_s^{{\bm \pi}^P}$ and $X_s^{R, \widehat{\bm \pi}^R}$, we can rewrite $Q_{d, \Delta t}$  as
\begin{equation*}
\begin{aligned}
Q_{d, \Delta t} (t, x, \mu, a, {\bm h}; \widehat {\bm \pi}, {\bm \pi})  =& \E^e\biggl[\int_t^{t + \Delta t} \int_{\Ac}e^{-\beta(s-t)}r(s, X_s^{R, \widehat{\bm \pi}^R}, \mu_s^{{\bm \pi}^P}, a)\dd s\biggl]\\
& + e^{-\beta \Delta t} J_d(t + \Delta t, X_{t + \Delta t}^{R, \widehat{\bm \pi}^R}, \mu_{t + \Delta t}^{{\bm \pi}^P}; \widehat{\bm \pi}, {\bm \pi}).
\end{aligned}
\end{equation*}
Applying It\^o's formula to $J_d(s, X_{s}^{R, \widehat{\bm \pi}^R}, \mu_{s}^{{\bm \pi}^P}; \widehat{\bm \pi}, {\bm \pi})$ between $t$ and $t + \Delta t$, we get the first-order approximation of $Q_{\Delta t}$
\begin{equation*}
\begin{aligned}
 &Q_{d, \Delta t} (t, x, \mu, a, {\bm h}; \widehat{\bm \pi}, {\bm \pi})\\
=& J_d(t, x, \mu; \widehat {\bm \pi}, {\bm \pi}) + \Delta t \biggl( \frac{\partial J_d}{\partial t}(t, x, \mu; \widehat{\bm \pi}, {\bm \pi}) - \beta J_d(t, x, \mu; \widehat {\bm \pi}, {\bm \pi}) + \mathcal{L}^{t, a, \mu}[J_d(t, \cdot, \mu; \widehat{\bm \pi}, {\bm \pi})](x) \nonumber\\
& + r(t, x, \mu, a)  + \E^e_{\xi \sim \mu}\Big[\int_{\Ac}\mathcal{L}^{t, a, \mu}[\frac{\delta J_d}{\delta \mu}(t, x, \mu; \widehat{\bm \pi}, {\bm \pi})(\cdot)](\xi) {\bm h}(a|t, \xi, \mu)\dd a\Big]
\biggl) + o(\Delta t).
\end{aligned}
\end{equation*}


\subsection{{Decoupled Iq-function}} As $\Delta t \to 0$, the decoupled IQ-function coincides with decoupled value function and hence can not be used to compare the action $a$, we focus on the first-order approximation, which leads to the following definition of the decoupled Iq-function.
\begin{Definition}\label{def:decoupled-q-function}
Fix policies $(\widehat{\bm \pi}, {\bm \pi}) \in \Pi^2$. For any $(t, x, \mu, a, {\bm h}) \in [0, T] \times \R^d \times \Pc_2(\R^d) \times \Ac \times \Pi$, the continuous time decoupled Iq-function is defined by
\begin{equation*}
\begin{aligned}
q_d(t, x, \mu, a, {\bm h}; \widehat{\bm \pi}, {\bm \pi}) := &  \frac{\partial J_d}{\partial t}(t, x, \mu; \widehat {\bm \pi}, {\bm \pi}) - \beta J_d(t, x, \mu; \widehat {\bm \pi}, {\bm \pi}) \\
&  + \mathcal{L}^{t, a, \mu}[J_d(t, \cdot, \mu; \widehat{\bm\pi}, {\bm \pi})](x) + r(t, x, \mu, a) \\
& + \E^e_{\xi \sim \mu}\Big[\int_{\Ac}\mathcal{L}^{t, a, \mu}[\frac{\delta J_d}{\delta \mu}(t, x, \mu; \widehat{\bm \pi}, {\bm \pi})(\cdot)](\xi) {\bm h}(a|t, \xi, \mu)\dd a\Big].
\end{aligned}
\end{equation*}

\end{Definition}
With the introduction of the decoupled Iq-function, the equation \eqref{DP-equation-J} can be rewritten as
\begin{equation}\label{DPE-of-J-in-terms-q}
\int_{\Ac}\big(q_d(t, x, \mu, a, {\bm \pi}; \widehat{\bm \pi}, {\bm \pi}) - \gamma \log\widehat{\bm \pi}(a|t, x, \mu)\big)\widehat{\bm \pi}(a|t, x, \mu)\dd a=0.
\end{equation}

The following result is a martingale characterization of the decoupled value function and the decoupled Iq-function associated with any pair of policies $(\widehat{\bm \pi}, {\bm \pi})$. 
\begin{Theorem}[Characterization of the decoupled Iq-function]\label{thm:martingale-decoupledJq}Given a pair of policies $(\widehat{\bm\pi}, {\bm\pi}) \in \Pi^2$. Let a continuous function $\widehat J_d: [0, T] \times \R^d \times \Pc_2(\R^d) \to \R$ and a continuous function $\widehat q_d: [0, T] \times  \R^d \times \Pc_2(\R^d) \times \Ac \times \Pi \to \R$ be given. Then $\widehat J_d$ and $\widehat q_d$ are the decoupled value function and decoupled Iq-function associated with $(\widehat{\bm \pi}, {\bm \pi}) \in \Pi^2$ if and only if
\begin{enumerate}[label={{\upshape(\roman*)}}]
\item (Consistency condition) $\widehat J_d$ and $\widehat q_d$ satisfy
\begin{equation}\label{thm:martingale-characterization-consistency}
\widehat J_d(T, x, \mu) = g(x, \mu),\;  \int_{\Ac} \Big[\widehat q_d(t, x, \mu, a, {\bm \pi}) - \gamma \log \widehat {\bm \pi}(a|t, x, \mu)\Big] \widehat {\bm \pi}(a|t, x, \mu)\dd a = 0.
\end{equation}
\item (Martingale condition) For any $(t, x, \mu, a, {\bm h}) \in [0, T] \times \R^d \times \Pc_2(\R^d) \times \Ac \times \Pi$, the following process
\begin{equation}\label{equ:martingale_characterization}
M_s^{t, x, \mu, {\bm h}}: = e^{-\beta s} \widehat J_d(s, X_s^{R, {\bm h}}, \mu_s^{{\bm h}}) + \int_{t}^s e^{-\beta t'}\Big[ r(t', X_{t'}^{R, {\bm h}}, \mu_{t'}^{{\bm h}}, {a}_{t'}^{{\bm h}})- \widehat q_d (t', X_{t'}^{R, {\bm h}}, \mu_{t'}^{{\bm h}}, {a}^{{\bm h}}_{t'}, {\bm h})\Big]\dd t'
\end{equation}
is a $\F^e$-adapted martingale, where $X_s^{R, {\bm h}}$ and $\mu_s^{{\bm h}}$, $t \leq s \leq T$ are solutions to \eqref{equ:agent-SDE} and \eqref{equ:population-average-SDE} under the policy ${\bm h}$, respectively.
\end{enumerate}
\end{Theorem}

Theorem \ref{thm:martingale-decoupledJq} provides a policy evaluation rule for the decoupled value function and the decoupled Iq-function. The two equations in the consistency condition \eqref{thm:martingale-characterization-consistency} correspond to the terminal condition and the dynamic programming equation of the decoupled value function. Besides the consistency condition, we need to utilize {\it all} test polices to generate samples in the martingale condition to fully characterize the decoupled value function and the decoupled Iq-function.  The proof of Theorem \ref{thm:martingale-decoupledJq} is  deferred to Section \ref{app:proof-martingale-decoupledJq}.



\section{Mean-Field Games and Mean-Field Controls}\label{sec:mfg-mfc}

\subsection{Mean-field games}\label{sec:mfg}
Fix $(t, x, \mu) \in [0, T] \times \R^d \times \Pc_2(\R^d)$.
 We say that ${\bm \pi}^{MFG, *} \in \Pi$ is a mean-field equilibrium (MFE) policy if
\begin{equation}\label{MFG-objective}
J_d(t, x, \mu; {\bm \pi}^{MFG, *}, {\bm \pi}^{MFG, *}) = \sup_{\widehat{\bm\pi} \in \Pi} J_d(t, x, \mu; \widehat{\bm\pi}, {\bm \pi}^{MFG, *}).
\end{equation}


When the policy of the population ${\bm \pi}$ is fixed (in other words, ${\bm \pi}$ is  viewed as an exogenous parameter), it is a classical optimal control problem. Similar to \cite{jiazhou2023}, the optimal policy $\widehat{\bm \pi}^*$ of the representative agent given ${\bm \pi}$ is a Gibbs measure in terms of the Hamiltonian. Fix $\vartheta = (t, x, \mu)$, we consider the operator $\Ic^{MFG}: \Pi \times \Pi \to \Pi$ for each pair $(\widehat {\bm \pi}, {\bm\pi}) \in \Pi^2$ that
\begin{equation}\label{MFG-PI-map}
\Ic^{MFG}(\widehat {\bm \pi}, {\bm \pi})(a|\vartheta) = \frac{\exp\left\{\frac{1}{\gamma}H\big(\vartheta, a, \partial_x J_d(\vartheta; \widehat{\bm \pi}, {\bm \pi}), \partial_{xx}^2 J_d(\vartheta; \widehat{\bm \pi}, {\bm \pi}),J_d(\cdot; \widehat{\bm \pi}, {\bm \pi})\big)\right\} }{\int_{\Ac}\exp\left\{\frac{1}{\gamma}H\big(\vartheta, a, \partial_x J_d(\vartheta; \widehat{\bm \pi}, {\bm \pi}), \partial_{xx}^2 J_d(\vartheta; \widehat{\bm \pi}, {\bm \pi}), J_d(\cdot; \widehat{\bm \pi}, {\bm \pi})\big) \right\}\dd a}.
\end{equation}

The following result is a direct consequence of policy improvement for classical optimal control problem, see \cite{jiazhou2023}, and the proof is omitted.

\begin{Proposition}[Policy improvement for MFG] \label{prop:PI-MFG}For each pair $(\widehat {\bm \pi}, {\bm\pi}) \in \Pi^2$, define $\widehat {\bm \pi}' = \mathcal{I}^{MFG}(\widehat {\bm \pi}, {\bm \pi})$, with $\mathcal{I}^{MFG}$ given in \eqref{MFG-PI-map}.
Then $J_d(t, x, \mu; \widehat {\bm \pi}', {\bm \pi}) \geq J_d(t, x, \mu; \widehat {\bm \pi}, {\bm \pi})$.
If $\widehat{\bm \pi}^* \in \Pi$ satisfies $\Ic^{MFG}(\widehat{\bm \pi}^*, {\bm \pi}) = \widehat{\bm \pi}^*$, then $J_d(t, x, \mu; {\widehat {\bm \pi}}^*, {\bm \pi}) = \sup_{\widehat {\bm \pi} \in \Pi} J_d(t, x, \mu; \widehat {\bm \pi}, {\bm \pi})$ and $\widehat{\bm \pi}^*$ is called the best response to ${\bm \pi}$.
\end{Proposition}
By \eqref{MFG-objective}, if there exists ${\bm \pi}^{MFG, *}$ such that $\Ic^{MFG}({\bm \pi}^{MFG, *}, {\bm \pi}^{MFG, *})= {\bm \pi}^{MFG, *},$ then it is the MFE policy. {If the MFE policy is unique, the corresponding $J_d^{MFG, *}(t, x, \mu): = J_d(t, x, \mu; {\bm \pi}^{MFG,*}, {\bm \pi}^{MFG, *})$ is referred to as the game value.} Substituting ${\bm \pi}^{MFG, *}$ into \eqref{DP-equation-J}, we derive the master equation for the game value:
\begin{align}
&\frac{\partial J_d^{MFG, *}}{\partial t}(\vartheta) - \beta J_d^{MFG, *}(\vartheta) +\int_{\R^d \times \Ac} \mathcal{L}^{t, a, \mu}\left[\frac{\delta J_d^{MFG, *}}{\delta \mu}(\vartheta)(\cdot)\right](y){\bm \pi}^{MFG, *}(a|t, y, \mu)\dd a \mu(\dd y)
\nonumber\\
 &+ \gamma \log \left[\int_{\Ac} \exp\Big\{\frac{1}{\gamma} H(\vartheta, a, \partial_x J_d^{MFG, *}(\vartheta), \partial_{xx} J_d^{MFG, *}(\vartheta), J_d^{MFG, *}(\vartheta))\Big\} \dd a\right] = 0.\label{master-equation-J}
\end{align}

In the context of RL in which the Hamiltonian is unknown, we elaborate how to utilize the decoupled Iq-function and its martingale characteriztation to learn the MFE policy.

Combining Theorem \ref{thm:martingale-decoupledJq} with the characterization of the MFE policy, we immediately obtain the verification argument of the decoupled value function and the decoupled Iq-function for MFG. The proof of Proposition~\ref{cor:MFG-martingale} is deferred to Section~\ref{subsec:martingale-MFG}.

{\begin{Proposition}\label{cor:MFG-martingale}
Let a continuous function $\widehat J^*_d: [0, T] \times \R^d \times \Pc_2(\R^d) \to \R$ with $\widehat J^*_d(T, x, \mu) = g(x, \mu)$ , a continuous function $\widehat q^*_d: [0, T] \times  \R^d \times \Pc_2(\R^d) \times \Ac \times \Pi \to \R$ and a policy $\widehat{\bm \pi}^*$ be given.
\begin{enumerate}[label={{\upshape(\roman*)}}]
\item (Consistency condition) $\widehat J^*_d$, $\widehat q^*_d$ and $\widehat {\bm \pi}^*$ satisfy \begin{equation}\label{equ:mfg-consistency}
 \int_{\Ac} \exp\left\{\frac{1}{\gamma}\widehat q^*_d(t, x, \mu, a, \widehat {\bm \pi}^*)\right\}\dd a= 1, \;\;\; \widehat {\bm \pi}^*(a|t, x, \mu) =  \frac{\exp \left\{\frac{1}{\gamma} \widehat q^*_d(t, x, \mu, a, {\bm h})\right\}}{\int_{\Ac} \exp \left\{\frac{1}{\gamma} \widehat q^*_d(t, x, \mu, a, {\bm h})\right\}\dd a}.
\end{equation}
\item for any $(t, x, \mu, a, {\bm h}) \in [0, T] \times \R^d \times \Pc_2(\R^d) \times \Ac \times \Pi$, $\widehat J_d^*$ and $\widehat q_d^*$ satisfy the martingale condition \eqref{equ:martingale_characterization}.
\end{enumerate}
Then $\widehat {\bm \pi}^*$ is a MFE policy and $\widehat J^*_d$, $\widehat q^*_d$ are respectively $J_d^{MFG, *}$ and $q_d^{MFG, *}$ associated with $\widehat {\bm \pi}^*$.
\end{Proposition}
}

\subsection{Mean-Field Control}\label{sec:mfc}

The MFC problem is also called the McKean-Vlasov control problem. The objective of MFC is to maximize the collective reward for the population
\begin{equation}\label{MFC-objective}
\sup_{{\bm \pi} \in \Pi}J^{MFC}(t, \mu; {\bm \pi}) = \sup_{{\bm \pi} \in \Pi}\int_{\R^d} J_d(t, x, \mu;{\bm \pi}, {\bm\pi})\mu(\dd x).
\end{equation}
 We call ${\bm \pi}^{MFC, *}$ a mean-field optimal (MFO) policy of  \eqref{MFC-objective} if
 $J^{MFC}(t, \mu; {\bm \pi}^{MFC, *}) = \sup_{{\bm \pi} \in \Pi} J^{MFC}(t, \mu; {\bm \pi})$.


Denote $J(t, x, \mu; {\bm \pi}): = \frac{\delta J^{MFC}}{\delta \mu}(t, \mu; {\bm \pi})(x)$ as the linear functional derivative of $J^{MFC}(t, \mu; {\bm \pi})$ for the simplicity of notation. Then
$J(t, x, \mu; {\bm \pi}) = J_d(t, x, \mu; {\bm \pi}, {\bm \pi}) + \int_{\R^d} \partial_\mu J_d(t, v, \mu; {\bm \pi}, {\bm \pi})(x) \mu(\dd v)$.
 Recall that $\vartheta= (t, x, \mu)$.
 As in \cite{weiyu2025}, let us consider the map $\mathcal{I}^{MFC}: \Pi \to \Pi$
\begin{equation}\label{equ:policy_improvemet_map}
 \mathcal{I}^{MFC}({\bm \pi})(a|\vartheta) = \frac{\exp\Big\{\frac{1}{\gamma}H\big(\vartheta, a,  \partial_x J(\vartheta; {\bm \pi}), \partial_{xx}^2  J(\vartheta; {\bm \pi}), J(\cdot; {\bm \pi})\big)\Big\}}{\int_{\Ac}\exp\Big\{\frac{1}{\gamma}H\big(\vartheta, a,  \partial_x J(\vartheta; {\bm \pi}), \partial_{xx}^2  J(\vartheta; {\bm \pi}), J(\cdot; {\bm \pi})\big)\Big\}\Big\} \dd a}.
\end{equation}

As shown in \cite{weiyu2025},  a new policy ${\bm \pi}'$ generated by ${\bm \pi}' = \mathcal{I}^{MFC}({\bm \pi})$ improves the previous policy ${\bm \pi}$.

\begin{Proposition}[Policy improvement for MFC]\label{prop:PI-MFC}For any given ${\bm \pi} \in \Pi$, let us define ${\bm \pi}' = \mathcal{I}^{MFC}({\bm \pi})$, with $\mathcal{I}^{MFC}$ given in \eqref{equ:policy_improvemet_map}.
Then it holds that $J^{MFC}(t, \mu; {\bm \pi}') \geq J^{MFC}(t, \mu; {\bm \pi})$.
Moreover, if the map $\mathcal{I}^{MFC}$ in \eqref{equ:policy_improvemet_map} has a fixed point ${\bm \pi}^{MFC, *} \in \Pi$, then ${\bm \pi}^{MFC, *}$ is a MFO policy of \eqref{MFC-objective}.
\end{Proposition}


\medskip
In the cooperative MFC framework, agents aim for a social optimum. This leads to the integrated q-function without the entropy term in \cite{weiyu2025}, which is defined by
\begin{equation}
q^{MFC}(t, \mu, {\bm h}; {\bm \pi}) = \E^e\left[q_d(t, x, \mu, a, {\bm h}; \widehat{\bm \pi}, {\bm \pi})|_{a = a^{\bm h}, \widehat{\bm \pi} = {\bm \pi}, x = \xi}\right].
\end{equation}
For the purpose of policy iteration, we need to introduce the essential q-function, which is closely related with the Hamiltonian of MFC and can be used for the policy improvement.
\begin{Definition}\label{realtion:qc-barq-qd} 
If there exists a function $q_{e}^{MFC}: [0, T] \times \R^d \times \Pc_2(\R^d) \times \Ac \to \R$ independent of ${\bm h}$ such that
\begin{equation}\label{relation-q-essentialq}
\E_{\mu, {\bm h}}^e\left[q_d(t, \xi, \mu, a^{\bm h}, {\bm h}; {\bm \pi}, {\bm \pi})\right] = \E_{\mu, {\bm h}}^e\left[q_{e}^{MFC}(t, \xi, \mu, a^{\bm h}; {\bm \pi})\right],
\end{equation}
we shall call $q_{e}^{MFC}$ the essential q-function. Here  we denote by $\E_{\mu, {\bm h}}^e [f(\xi, a^{\bm h})]$ the integral of the function $f$ with respect to ${\bm h}$ and $\mu$, that is, $\E_{\mu, {\bm h}}^e[f(\xi, a^{\bm h})] = \int_{\R^d \times \Ac} f(x, a){\bm h}(a|t, x, \mu) \dd a\mu(\dd x)$.
\end{Definition}
Such an essential q-function $q_{e}^{MFC}$ always exists. For example, we can define $q_{e}^{MFC}$ by
\begin{equation*}
\begin{aligned}
 q_{e}^{MFC}(t, x, \mu, a; {\bm \pi}) =& \frac{\partial J_d}{\partial t}(t, x, \mu; {\bm \pi}, {\bm \pi}) - \beta J_d(t, x, \mu; {\bm \pi}, {\bm \pi})\\
  &+ H\left(t, x, \mu, a, \partial_x J(t, x, \mu; {\bm \pi}), \partial_{xx}^2 J(t, x, \mu; {\bm \pi}), J(\cdot;{\bm \pi})(x)\right),
\end{aligned}
\end{equation*}
where recall that $J(t, x, \mu; {\bm \pi})$ is the linear functional derivative of $J^{MFC}(t, \mu; {\bm \pi})$.
The decoupled value function and the decoupled Iq-function associated with the MFO policy ${\bm \pi}^{MFC, *}$ are respectively denoted by $J^{MFC, *}_d$ and $q^{MFC, *}_d$. The proof of Proposition~\ref{cor:MFC-martingale} is deferred to Section~\ref{subsec:martingale-MFC}.

\begin{Proposition}\label{cor:MFC-martingale}
 Let $\widehat J^*_d: [0, T] \times \R^d \times \Pc_2(\R^d) \to \R$ with $\widehat J^*_d(T, x, \mu) = g(x, \mu)$, and $\widehat q^*_d: [0, T] \times  \R^d \times \Pc_2(\R^d) \times \Ac \times \Pi \to \R$ be given continuous functions. Then $\widehat J^*_d$, $\widehat q^*_d$ and $\widehat {\bm \pi}^*$ are respectively $J^{MFC, *}_d$, $q^{MFC, *}_d$ and a MFO policy ${\bm \pi}^{MFC, *}$ if and only if
 \begin{enumerate}[label={{\upshape(\roman*)}}]
\item (Consistency condition)
 \begin{equation}\label{equ:mfc-consistency-q}
  \int_{\Ac} \big(\widehat{q}^*_d(t, x, \mu, a,\widehat {\bm \pi}^*) - \gamma \log \widehat {\bm \pi}^*(a|t, x, \mu)\big)\widehat {\bm \pi}^*(a|t, x, \mu) \dd a =0,
  \end{equation}
  \begin{equation}
 \E_{\mu, {\bm h}}[\widehat q^*_d(t, \xi, \mu, a^{\bm h}, {\bm h})]  = \E_{\mu, {\bm h}}[\widehat q^*_{e}(t, \xi, \mu, a^{\bm h})],\;\;\; \widehat {\bm \pi}^*(a|t, x, \mu) =  \frac{\exp \left\{\frac{1}{\gamma} \widehat q^*_e(t, x, \mu, a)\right\} }{\int_{\Ac} \exp\left\{\frac{1}{\gamma}\widehat q^*_e(t, x, \mu, a)\right\}\dd a }, \label{equ:mfc-consistency-pi}
 \end{equation}
for some $\widehat q^*_e: [0, T] \times \R^d \times \Pc_2(\R^d) \times \Ac \to \R$.
\item for any $(t, x, \mu, {\bm h}) \in [0, T] \times \R^d \times \Pc_2(\R^d) \times \Pi$, $\widehat J_d^*$ and $\widehat q_d^*$ satisfy the martingale condition \eqref{equ:martingale_characterization}.
\end{enumerate}
\end{Proposition}

\begin{Remark}[Comparison between MFG and MFC]Both Propositions \ref{cor:MFG-martingale} and \ref{cor:MFC-martingale} characterize the decoupled value function and decoupled Iq-function associated with the MFE policy or the MFO policy in terms of the consistency condition and the martingale condition \eqref{equ:martingale_characterization} using all test policies. Notably, MFG and MFC share the same martingale condition, which will be the foundation for developing a unified learning algorithm for MFG and MFC. Although the decoupled Iq-function is not used directly for the MFO policy of MFC problem, its connection with essential  q-function via the integral representation \eqref{relation-q-essentialq} is vital for the q-learning of MFC.The subtle difference in the consistency condition leads to different constraints on the paramaterized approximators of the decoupled value function and the decoupled Iq-function which will be discussed in the next section.
\end{Remark}

\section{Unified q-Learning Algorithm for MFG and MFC}\label{sec:decoupled-q-algo}
In this section,  we devise the unified q-learning algorithm {for MFG and MFC problems} according to Propositions \ref{cor:MFG-martingale} and \ref{cor:MFC-martingale}.

\noindent \paragraph{{Function approximation under consistency condition}} First, we choose the {parametrized} function approximators $J^\theta_d$ and $q^\psi_d$ such that the consistency condition holds. In both MFG and MFC formulations, the parameterized function approximator $J^\theta_d$ of $J_d^{MFG, *}$ (resp. $J^{MFC, *}_d$) satisfies the terminal condition $J^\theta_d(T, x, \mu) = g(x, \mu)$. The key is the function approximator of $q^{MFG, *}_d$ (resp. $q^{MFC, *}_d$).  By Definition \ref{def:decoupled-q-function}, let the parameterized approximator of $q^{MFG, *}_d$ (resp. $q^{MFC, *}_d$) take a separable form:
\begin{equation}\label{param-q}
q^\psi_d(t, x, \mu, a, {\bm h}) = q_1^\psi(t, x, \mu, a) + \int_{\R^d \times \Ac} q_2^\psi(t, x, \mu, a, y) {\bm h}(a|t, y, \mu)\dd a\mu(\dd y)
\end{equation}
for some parametrized functions $q_1^\psi$ and $q_2^\psi$. 

In MFG, the consistency condition \eqref{equ:mfg-consistency} shows that $q_d^\psi$ and ${\bm \pi}^\psi$ should satisfy
\begin{align}\label{mfg:param-pi}
{\bm \pi}^\psi(a|t, x, \mu) &= \frac{\exp\big\{\frac{1}{\gamma}  q_1^\psi(t, x, \mu, a)\big\}}{\int_{\Ac}\exp\big\{\frac{1}{\gamma}  q_1^\psi(t, x, \mu, a)\big\}\dd a},\\
\gamma \log \int_\Ac \exp\big(\frac{1}{\gamma} q_1^\psi(t, x, \mu, a)\big)\dd a &+ \int_{\R^d \times \Ac} q_2^\psi(t, x, \mu, a, y) {\bm \pi}^\psi(a|t, y, \mu)\dd a\mu(\dd y)=0.\label{mfg:param-q}
\end{align}

In MFC, first by the integral representation between $q^{MFC, *}_d$ and $q_e^{MFC, *}$ in \eqref{relation-q-essentialq}, we choose the function approximator $q_e^\psi$ of $q_e^{MFC, *}$ in the form
\begin{equation}\label{mfc:param-q-e}
q_e^\psi(t, x, \mu, a) = q_1^\psi(t, x, \mu, a) + \int_{\R^d} q_2^\psi(t, y, \mu, a, x) \mu(\dd y).
\end{equation}
From the consistency condition \eqref{equ:mfc-consistency-q}-\eqref{equ:mfc-consistency-pi}, we have
\begin{align}\label{mfc:param-pi}
{\bm \pi}^\psi(a|t, x, \mu) &= \frac{\exp\big\{\frac{1}{\gamma}  q_e^\psi(t, x, \mu, a)\big\}}{\int_{\Ac}\exp\big\{\frac{1}{\gamma}  q_e^\psi(t, x, \mu, a)\big\}\dd a},\\
\int_{\Ac} \big(q^\psi_d(t, x, \mu, a, {\bm \pi}^\psi) &- \gamma \log {\bm \pi}^\psi(a|t,x, \mu)\big) {\bm \pi}^\psi(a|t, x, \mu)\dd a =0.\label{mfc:param-q}
\end{align}

\noindent \paragraph{{Unified averaged martingale orthogonality condition}} With careful parameterizations \eqref{param-q}-\eqref{mfc:param-q} for MFG and MFC, it suffices to maintain the martingale condition in \eqref{equ:martingale_characterization} under all test policies.

However, it is impossible to exhaust all test policies ${\bm h} \in \Pi$ in the practical implementation to ensure that \eqref{equ:martingale_characterization} holds. In response, we restrict to a parametrized  test policy space $\Pi^\Psi: = \{{\bm h}^{\tilde\psi}\}_{\tilde\psi \in \Psi(\psi)}$, where $\Psi(\psi) \subset \R^{L_\psi}$. We take $\Psi(\psi)$ as $\Psi$ without loss of generality. Denote
\begin{align}\label{equ:martingale_characterization_para}
M_t^{{\tilde\psi}, \theta, \psi}: = e^{-\beta t} J^\theta_d(t,X_t^{R, {\bm h}^{\tilde\psi}}, \mu_t^{{\bm h}^{\tilde\psi}}) + &  \int_{0}^t e^{-\beta t'}\Big[ r(t', X_{t'}^{R, {\bm h}^{\tilde\psi}}, \mu_{t'}^{{\bm h}^{\tilde\psi}} , {a}^{{\bm h}^{\tilde\psi}}_{t'})\\
&\;  - q^\psi_d (t', X_{t'}^{R, {\bm h}^{\tilde\psi}}, \mu_{t'}^{{\bm h}^{\tilde\psi}} , {a}^{{\bm h}^{\tilde\psi}}_{t'}, {\bm h}^{\tilde\psi})\Big]dt'. \nonumber
\end{align}
Now, we state a parameterized version of Propositions \ref{cor:MFG-martingale} and \ref{cor:MFC-martingale}. 
\begin{Proposition} \label{prop:para-martingale-condition} Let a continuous function $J^\theta_d: [0, T] \times \R^d \times \Pc_2(\R^d) \to \R$ with $J^\theta_d(T, x, \mu) = g(x, \mu)$, a continuous function $q^\psi_d: [0, T] \times  \R^d \times \Pc_2(\R^d) \times \Ac \times \Pi^\Psi \to \R$ and a policy ${\bm \pi}^\psi$ be given. Suppose that $q^\psi_d$ takes the form  \eqref{param-q} and satisfies \eqref{mfg:param-pi}-\eqref{mfg:param-q} (resp. \eqref{mfc:param-q-e}-\eqref{mfc:param-q}). Furthermore, there exists $(\theta^*, \psi^*) \in \Theta \times \Psi$ such that for any $(t, x, \mu, a, \tilde\psi) \in [0, T] \times \R^d \times \Pc_2(\R^d) \times \Ac \times \Psi$, the process $M_t^{{\tilde\psi}, \theta^*, \psi^*}$ in \eqref{equ:martingale_characterization_para} is a $\F^e$-adapted martingale. Then
 ${\bm \pi}^{\psi^*} = {\bm \pi}^{MFG, *}$ (resp. ${\bm \pi}^{MFC, *}$), $J_d^{\theta^*} = J^{MFG, *}_d$ (resp. $J_d^{MFC, *}$), and $q^{\psi^*}_d = q^{MFG, *}_d$ (resp. $q_d^{MFC, *}$) restricted on $\Pi^{\Psi}$.
\end{Proposition}

It is noted that applying the martingale condition to test policies in $\Pi^\Psi$ rather than $\Pi$ leads us to conclude that $q^\psi$ is equal to decoupled Iq-function only within the restricted set $\Pi^\psi$. Nevertheless, this suffices to determine the MFE policy and the MFO policy, as well as the associated decoupled value functions.

{
\begin{Remark}
Proposition~\ref{prop:para-martingale-condition} describes the ideal case where the parametrized function classes for the decoupled value function, the decoupled Iq-function, and the test policies are rich enough so that the consistency conditions and the martingale conditions can be
satisfied by the true solution; in this case the learned parameters recover the optimal ones and there is no approximation error.

In general, when the optimal objects do not lie in the chosen parametrized classes, a quantitative analysis of the
resulting function-approximation error would require additional assumptions on the richness of the parametrized
classes and on the regularity of the optimal solutions, and is beyond the scope of this paper. We leave a rigorous
approximation-error analysis for parametrized continuous-time mean-field RL
for future research. In particular, some recent studies on RL with linear and nonlinear function approximations (e.g. \cite{LongHan2023}) might be helpful to the study of our continuous-time mean-field setting.

\end{Remark}
}

According to Proposition \ref{prop:para-martingale-condition}, we need to consider the following averaged martingale orthogonality conditions (AMOC for short) under all parameterized test function $\xi = (\xi^{\tilde\psi})_{\tilde\psi \in \Psi}$, indexed by $\tilde\psi$, that
\begin{align}\label{MOC}
0 = \int_{\Psi} \left[\int_0^T \E^e \left[\xi_t^{\tilde\psi} d M_t^{{\tilde\psi}, \theta, \psi}\right] \right] \Uc(\dd \tilde\psi), \tag{AMOC}
\end{align}
where $\mathcal{U}$ denotes the uniform distribution on the parameter space $\Psi$. The following lemma shows that  \eqref{MOC} is a necessary and sufficient condition for $M_t^{{\tilde\psi}, \theta, \psi}$ being a martingale for any $\tilde\psi \in \Psi$. The proof of Lemma \ref{lemma:AMOC} is similar to Proposition 2 in \cite{JZ22b}, and hence it is omitted here.

\begin{lemma} \label{lemma:AMOC} The process \eqref{equ:martingale_characterization_para} is a $\F^e$-adapted martingale for every ${\bm h}^{\tilde\psi}$, $\tilde\psi \in \Psi$, if and only if \eqref{MOC} holds.
\end{lemma}

While all above theoretical analyses of q-learning algorithms are conducted in a continuous setting, the actual implementations of algorithms are discrete, using a finite number of test policies and a fixed time mesh size. Precisely, we discretize $[0, T]$ with fixed mesh size $\Delta t$, where $\{t_k = k \Delta t\}_{0 \leq k \leq K}$ is a time discretization of $[0, T]$, and choose $M$ test policies ${\bm h}^{\tilde\psi_1}$, $\ldots$, ${\bm h}^{\tilde\psi_M}$, where $\tilde\psi_m$, $1 \leq m \leq M$ are uniformly drawn from $\Psi$. The next theorem answers the question whether q-learning algorithm with these discretizations converges to the solution of \eqref{MOC} as $M \to +\infty$ and $\Delta t \to 0$. For any continuous-time process $\{X_t\}_{0 \leq t \leq T}$, we denote by $\{\bar X_{t_i}\}_{0 \leq i \leq K-1}$ its discrete-time version in the sequel.

We impose the following assumptions on the test functions and parameterized functions.
\begin{Assumption}
\label{assu:MOC-converegence}
\begin{itemize}
\item [(i)] A test function $\xi = \{\xi^{\tilde\psi}\}_{\tilde\psi \in \Psi}$ is measurable with respect to $\tilde\psi$ and satisfies
\begin{equation*}
\int_{\Psi} \E^e[\int_0^T\|\xi_t^{\tilde\psi}\|^2dt]^{1/2} \Uc(d\tilde\psi) < +\infty,  \; \; \; \sup_{\tilde\psi \in \Psi}\E^e[|\xi^{\tilde\psi}_{t'} -  \xi^{\tilde\psi}_t|^2] \leq C_{\tilde\psi, \theta, \psi}|t - t'|^\alpha
\end{equation*}
for any $t, t' \in [0, T]$, where $C_{\tilde\psi, \theta, \psi}$ is a continuous function of $\tilde\psi$, $\theta$ and $\psi$, and $\alpha >0$ is a given constant.
\item [(ii)] $J_d^\theta$ and $q_d^\psi$ are sufficiently smooth functions such that all derivatives required exist in classical sense. Moreover, for all $(\theta, \psi) \in \Theta \times \Psi$, $\check q_d^\theta$ and $q_d^\psi$ are square integrable in the sense that
\begin{equation*}
\begin{aligned}
&\int_{\R^d \times \Ac} \big(|\check q_d^{\theta}(t, x, \mu, a, {\bm h}^{\tilde\psi})|^2 + |q_d^\psi(t, x, \mu, a, {\bm h}^{\tilde\psi})|^2\big){\bm h}^{\tilde\psi}(a|t, x, \mu)\dd a\mu(\dd x)\\
 \leq &  C_{\tilde\psi, \theta, \psi}(1 + |x|^2 + M_2(\mu)^2),
\end{aligned}
    \end{equation*}
     where $\check q_d^\theta$ is defined similarly as $q_d$ in Definition \ref{def:decoupled-q-function} but with $J_d$ replaced by $J_d^\theta$.
\end{itemize}
\end{Assumption}

\begin{Theorem} \label{thm:approximation-M-Deltat} Denote by $(\theta^*_{\Delta t, M}, \psi^*_{\Delta t, M})$ the solution to discrete averaged martingale orthogonality condition (DAMOC for short) with $M$ test policies and mesh size $\Delta t$:
\begin{align}\label{DAMOC}
\frac{1}{M} \sum_{m=1}^M \sum_{i=0}^{K-1} \xi_{t_i}^m \big(\bar M_{t_{i+1}}^{{\bm h}^{\tilde\psi_m}, \theta, \psi} - \bar M_{t_i}^{{\bm h}^{\tilde\psi_m}, \theta, \psi}\big) = 0. \tag{DAMOC}
\end{align}
Then under Assumption \ref{assu:MOC-converegence}, any convergent subsequence of $(\theta^*_{\Delta t, M}, \psi^*_{\Delta t, M})$ converges to the solution to \eqref{MOC},
\begin{equation*}
(\theta^*, \psi^*) : = \lim_{\Delta t \to 0, M \to \infty} (\theta^*_{\Delta t, M}, \psi^*_{\Delta t, M}).
\end{equation*}
\end{Theorem}
The proof of Theorem \ref{thm:approximation-M-Deltat} is given in Section \ref{app:proof-thm-approximation-M-Deltat}.

\paragraph{{Population estimation and algorithm implementation}} 
{A practical challenge in the learning procedure is that the population distribution is not directly provided by the environment simulator. Instead, the representative agent must estimate the distribution flow based on the local observations of state trajectories.} We provide two ways to estimate the population distribution. Let $X_{t_k}^{j, m}$, $1 \leq j \leq N, 1 \leq m \leq M$ denote the observation of the state of the representative player following the test policy ${\bm h}^{\tilde\psi_m}$ at time $t_k$ for episode $j$. First, similar to \cite{FGLPS23}, \cite{AFL2022}, we can update the state distribution for each episode $j =1, \ldots, N$ as
\begin{equation}\label{equ:distribution-update}
\mu_{t_k}^{j, m} = (1 - \rho^{j}) \mu_{t_k}^{j-1, m} + \rho^j \delta_{X_{t_k}^{j, m}},
\end{equation}
where $(\rho^j)_j$ is a sequence of learning rates. Alternatively, inspired by the way of updating the mean in the example of mean-variance portfolio optimization in \cite{jiazhou2023}, we can update the state distribution per $l$ epoches:
\begin{equation}\label{equ:distribution-update1}
\mu_{t_k}^{j, m} = (1 - \rho) \mu_{t_k}^{j-l-1, m} + \frac{\rho}{l} \sum_{i=j-l+1}^{j} \delta_{X_{t_k}^{i, m}}.
\end{equation}

Therefore, the representative agent in the state $X_{t_k}^{j, m}$  estimates the population state distribution according to \eqref{equ:distribution-update} or \eqref{equ:distribution-update1}. By interacting with the environment (e.g. by environment simulator), he takes an action $a_{t_k}^{j, m} \sim {\bm h}^{\tilde\psi_m}(\cdot|t_k, X_{t_k}^{j, m}, \mu_{t_k}^{j, m})$, observes the reward $r_{t_k}^{j, m}$ and moves to the new state $X_{t_{k+1}}^{j, m}$. 
We will omit the superscript $j$ in the following if there is no confusion.

Once samples $(X_{t_k}^m, a_{t_k}^m, \mu_{t_k}^m, r_{t_k}^m, X_{t_{k+1}}^m)_{0 \leq t \leq K -1}$  generated by the test policy ${\bm h}^{\tilde\psi_m}$ are available,  we choose in \eqref{DAMOC}  the conventional test functions below
\begin{equation}\label{test_orthogonal}
\begin{aligned}
\xi_{t_k}^{\tilde\psi_m} =\frac{\partial J^\theta_d}{\partial \theta}\left(t, X_{t_k}^m, \mu_{t_k}^{m}\right),
\mathrm{and} \quad &\zeta_{t_k}^{\tilde\psi_m} =\frac{\partial q^\psi_d}{\partial \psi}\left(t,X_{t_k}^m, \mu_{t_k}^{m}, a^{m}_{t_k},\bh^{\tilde\psi_m}\right).
\end{aligned}
\end{equation}

Then, with learning rates $\alpha_\theta$ and $\alpha_\psi$ for $\theta$ and $\psi$ respectively, we obtain the following update rules of parameters:
\begin{equation*}
\theta \leftarrow \theta + \alpha_\theta \frac{1}{M} \sum_{m=1}^M \Delta^m \theta, \; \psi \leftarrow \psi + \alpha_{\psi} \frac{1}{M}\sum_{m=1}^M\Delta^m \psi,
\end{equation*}
where $G_{t_k}^m$, $\Delta^m \theta$ and $\Delta^m \psi$ are defined by
\begin{align}\label{method2-expression-Gd-tk-T}
G_{t_k}^m = &J^{\theta}_d(t_{k+1},X^m_{t_{k+1}},\mu^m_{t_{k+1}})-J^{\theta}_d(t_{k},X^m_{t_{k}},\mu^m_{t_k}) \\
  &+\big(r^m_{t_k}-\beta J^{\theta}_d(t_{k},X^m_{t_{k}},\mu^m_{t_k}) - q^\psi_d(t_k,X^m_{t_k},\mu^m_{t_k},a^m_{t_k},h^{\tilde{\psi}_m})\big) \Delta t,   \nonumber\\
\Delta^m \theta =& \sum_{k=0}^{K-1} e^{-\beta t_k}G_{t_k}^m \frac{\partial {J}^\theta_d}{\partial \theta}(t_k, X_{t_k}^m, \mu_{t_k}^m),\\
\Delta^m \psi   = & \sum_{k=0}^{K-1} e^{-\beta t_k}G_{t_k}^m \frac{\partial q^\psi_d}{\partial\psi}\left(t_k,X^m_{t_k},\mu^m_{t_k},a^m_{t_k},h^{\tilde{\psi}_m}\right).
\label{method2-expression-Deltadm-tildepsi}
\end{align}
We provide the pseudo-code algorithm in Algorithm \ref{algo:decoupled-q} and illustrate the procedure of unified q-learning algorithm in Figure \ref{figure:procedure-q} as below.

\begin{algorithm}[t!]
\caption{Offline Unified {Parametric} q-Learning Algorithm}
\textbf{Inputs}: horizon $T$, time step $\Delta t$, number of mesh grids $K$, { number of test policies ${\bm h}^\psi$}, initial learning rates $\alpha_{\theta},\alpha_{\psi}$,  temperature parameter $\gamma$, and functional forms of parameterized  value function $J^{\theta}_d$ with $J_d(T, x, \mu) = g(x, \mu)$, $q^{\psi}_d$ satisfying \eqref{mfg:param-pi}-\eqref{mfg:param-q} in MFG or \eqref{mfc:param-q-e}-\eqref{mfc:param-q} in MFC.

\textbf{Required program}: environment simulator $(x', r) = \textit{Environment}_{\Delta t}(t, x, a)$ that takes current time--state pair $(t, x)$ and the action $a$ as inputs and generates new state $x'$ at time $t+\Delta t$ and the reward $r$ at time $t$ as outputs.

\textbf{Initialization}: state distribution $\mu_{t_k}^m$ on $\R^d$ for $k=0, 1, \ldots, K$, $m=1, \ldots, M$, parameters $\theta$ and $\psi$.

\textbf{Learning procedure}:
\begin{algorithmic}[1]
\FOR{episode $j=1$ to $N$}
\STATE{Initialize $k = 0$.
\FOR{$m=1$ to $M$}
\STATE{Initialize $X_0^m \sim \mu_0^m {:= \mu}$.}
\STATE{Draw $\tilde\psi^m$ from $\psi \cdot \mathcal{U}([p(j), q(j)])$ and set the test policy ${\bm h}^{\tilde\psi^m}$.}
\WHILE{$k < K$} \STATE{
Update state distribution $\mu_{t_k}^m = (1 - \rho^j) \mu_{t_k}^m + \rho^j \delta_{X_{t_k}^m}^m$ or according to \eqref{equ:distribution-update1}.
Generate $a_{t_k}^m \sim {{\bm h}}^{\tilde\psi^m}(\cdot|t_k, X_{t_k}^m, \mu_{t_k}^m)$.
Apply $a_{t_k}^m$ to  environment simulator $(x', r) = Environment_{\Delta t}(t_k, X_{t_k}^m, a_{t_k}^m)$, and observe  new state $x$ and the reward $r$ as output. Store $X_{t_{k+1}}^m \leftarrow x$ and $r_{t_k}^m \leftarrow r$.
Update $k \leftarrow k + 1$.
}
\ENDWHILE	
\ENDFOR

For every $k=0,1,\cdots,K-1$, compute $G_{t_k}^m$, $\Delta^m \theta$ and $\Delta^m \psi$ according to \eqref{method2-expression-Gd-tk-T}-\eqref{method2-expression-Deltadm-tildepsi}.
Update $\theta$ and $\psi$ by $\theta \leftarrow \theta + \alpha_\theta \frac{1}{M} \sum_{m=1}^M \Delta^m \theta, \psi \leftarrow \psi + \alpha_{\psi} \frac{1}{M} \sum_{m=1}^M \Delta^m \psi$.
}

\ENDFOR
\end{algorithmic}
\label{algo:decoupled-q}
\end{algorithm}

\begin{figure}[t!]
\centering
\begin{tikzpicture}[font=\tiny\sffamily\bfseries, node distance=2cm and 2cm, auto]
    \node (A) [draw, rectangle,  fill=red!10, text width=2cm, align=center] { approximators: $J_d^\theta$  and $q_d^\psi$};
    \node (B) [draw, rectangle, fill=blue!10, right= 3.5cm of A, text width=2cm, align=center, xshift=1cm, yshift=1cm] {approximated policy: ${\bm \pi}^\psi$};
    \node (C) [draw, rectangle, fill=blue!10, right=of A, text width=2cm, align=center, xshift=1cm, yshift=-1cm] { approximated essential q: $q_e^\psi$};
    \node (D) [draw, rectangle, fill=blue!10, text width=2cm, align=center, right=1.2cm of C] {approximated policy: ${\bm \pi}^{\psi}$};
    \node (E) [draw, rectangle, fill=yellow!10, text width=3.2cm, align=center, right=of D, yshift=1cm] {Samples: $\big(X_{t_k}^m, a_{t_k}^m, \mu_{t_k}^m, r_{t_k}^m, X_{t_{k+1}}^m\big)_k$ };

    \draw[->, >=latex] (A.east) to node[midway, above, sloped] {MFG: \eqref{mfg:param-pi}-\eqref{mfg:param-q}} (B.west);
    \draw[->, >=latex] (A.east) to  node[midway, below, sloped] {MFC: \eqref{mfc:param-q-e}}(C.west);
    \draw[->, >=latex] (B.east) to node[midway, above, sloped] {test policy ${\bm h}^{\tilde\psi^m}$}(E.west);
    \draw[->, >=latex] (C) -- node[midway, above]{\eqref{mfc:param-pi}-\eqref{mfc:param-q}}(D);
    \draw[->, >=latex] (D.east) to node[midway, below, sloped]{test policy ${\bm h}^{\tilde\psi^m}$} (E.west);

    \path (E) edge [->, >=latex, bend left=40] node[midway, below]{\eqref{DAMOC}}(A);

\end{tikzpicture}
\caption{Flowchart of the unified q-learning algorithm for MFG and MFC}\label{figure:procedure-q}
\end{figure}

\section{Numerical Examples}\label{sec:example}

\subsection{MFG and MFC under the mean-variance criterion}\label{MV-exm}
In this example, we apply the RL algorithm for the well-known portfolio optimization under the mean-variance criterion. To this end, let us consider the wealth process that satisfies the jump-diffusion SDE that
\begin{equation}\label{equ:MV_exploratory_average_SDE}
\dd X_s = a_s\big(b\dd s +\sigma \dd W_s + \Gamma \dd \widetilde N_s\big),\quad  a_s \sim {\bm \pi}(\cdot|s, X_s, \mu_s^{\bm \pi}),\; s \geq t,
\end{equation}
where $a_s$ is the wealth amount invested in the risky asset at time $s$, $b$ is the excess return, $\sigma$ is the volatility, $\Gamma$ is the size of jump risk, and $d{\widetilde N}_s = \dd N_t - \eta \dd t$ is a compensated Poisson random measure with constant parameter $\eta$. The coefficients of the dynamics \eqref{equ:MV_exploratory_average_SDE} do not depend on the distribution of the solution. Several continuous time RL algorithms have been developed in the past few years to learn the mean-variance portfolio policy when the model parameters are unknown. In particular, to overcome the time-inconsistency issue rooted in the dynamic mean-variance formulation, \cite{WangZhou2020} studied the RL algorithm for the minimum variance problem in a diffusion model subjecting to expectation constraint and transformed the learning problem into an unconstrained one by using the Lagrange multiplier. Later, \cite{GLZ2024} investigated the same pre-committed strategy for the minimum variance problem in the jump-diffusion model. \cite{DDJ} developed the RL algorithm for learning the time-consistent equilibrium portfolio under the mean-variance criterion. \cite{weiyu2025} lifted the state space and formulated the time-consistent learning MFC problem under the mean-variance criterion using the q-learning algorithm for the social planner in the diffusion setting.

In the present example, we again lift the state space but in the decoupled form so that the optimal control problem for the representative agent is time consistent. Our goal is to apply the decoupled Iq-function and unified q-learning algorithm to learn the portfolio policy in both MFG and MFC problems. Note that the forms of dynamics for the representative agent and the rest of the population are the same except the initial value.
The decoupled value function is given by
\begin{equation*}
\begin{aligned}
J_d(t, x, \mu; \widehat{\bm \pi}, {\bm \pi}) =& \E^e\biggl[X_T^{t, x, \mu, \widehat{\bm \pi}}-\lambda(X_T^{t, x, \mu, \widehat{\bm \pi}} - \E^e[X_T^{t, \xi, \bm \pi}])^2\\
&- \gamma  \int_t^T \int_{\Ac}\log \widehat{\bm \pi}(a|s, X_s^{t, x, \mu, \widehat{\bm \pi}},  \mu_s^{\bm \pi}) \widehat {\bm \pi}(a|s, X_s^{t, x, \mu, \widehat{\bm \pi}},  \mu_s^{\bm \pi})\dd a\dd s \biggl].
\end{aligned}
\end{equation*}

We aim to solve the mean-variance portfolio optimization problem in the context of both MFG and MFC. The detailed derivations in both formulations are given in the {online companion}. Interestingly, due to the specific LQ structure, it turns out that the decoupled value functions (and therefore the decoupled Iq functions) for both MFG and MFC problems are the same, i.e., the MFG problem under the mean-variance criterion is a potential MFG. The MFE policy for MFG and the MFO policy for MFC coincide with the same expression. For this reason, we can drop the superscript ``MFG" and ``MFC", and obtain the explicit expressions of $J_d^*$ and $q_d^*$ in both MFG and MFC problems as follows:
\begin{equation*}
\begin{aligned}
J_d^*(t, x, \mu) &= -\lambda e^{\frac{b^2}{\sigma^2 + \eta \Gamma^2}(t-T)}(x - \bar\mu)^2 + x +\frac{\gamma b^2}{4(\sigma^2 + \eta \Gamma^2)} (t- T)^2 \\
& \;\;\; - (t -T) \frac{\gamma}{2} \log \frac{\pi\gamma}{(\sigma^2 + \eta \Gamma^2)\lambda} + \frac{1}{4\lambda} e^{(-\frac{b^2}{\sigma^2 + \eta \Gamma^2}(t-T)} - \frac{1}{4\lambda},
\end{aligned}
\end{equation*}
and
\begin{equation*}
\begin{aligned}
q_d^*(t, x, \mu, a, {\bm h})
= & -\lambda(\sigma^2 + \eta \Gamma^2) e^{\frac{b^2}{\sigma^2 + \eta \Gamma^2}(t-T)}\Big(a + \frac{b}{\sigma^2 + \eta \Gamma^2}(x - \bar\mu) \\
&-\frac{b}{2\lambda(\sigma^2 + \eta \Gamma^2)} e^{-\frac{b^2}{\sigma^2 + \eta \Gamma^2}(t-T)}\Big)^2 - \frac{b^2}{\sigma^2 + \eta \Gamma^2} (x - \bar\mu) \\
&- \frac{\gamma}{2} \log \frac{\pi \gamma}{(\sigma^2 + \eta \Gamma^2) \lambda} + \frac{\gamma b^2}{2(\sigma^2 + \eta \Gamma^2)}(t -T) \\
&+2\lambda b e^{\frac{b^2}{(\sigma^2 + \eta \Gamma^2)}(t-T)} (\bar\mu- x)\E_{\mu, {\bm h}}\big[a^{\bm h}\big].
\end{aligned}
\end{equation*}
As a consequence, both the MFE policy and the MFO policy are in the same expression of
\begin{equation}\label{mean-variance-pi*}
{\bm \pi}^{*}(\cdot|t, x, \mu) = \Nc\Big(-\frac{b}{\sigma^2 + \eta \Gamma^2} \big(x - \bar\mu - \frac{1}{2\lambda}e^{-\frac{b^2}{\sigma^2 + \eta \Gamma^2}(t-T)}\big), \frac{\gamma}{2\lambda (\sigma^2 + \eta \Gamma^2)}e^{-\frac{b^2}{\sigma^2 + \eta \Gamma^2}(t-T)}\Big).
\end{equation}

\begin{Remark}\label{rmk:MV}
From the explicit MFE policy (the MFO policy) in \eqref{mean-variance-pi*} for the MFG (MFC) under the mean-variance criterion, it is straightforward to observe that the larger jump size $\Gamma$ or the larger jump intensity $\eta$ leads to the smaller variance of the randomized Gaussian policy $\pi^*$ in \eqref{mean-variance-pi*}. { As shown in \eqref{mean-variance-pi*}, the variance of the optimal policy is inversely related to the jump size $\Gamma$. This suggests that the jump component introduces additional randomness into the state dynamics, enhancing exploration in the state space. Consequently, a smaller policy variance suffices to achieve adequate exploration.} Therefore, in this particular LQ model with a larger jump risk, we may need to choose a relatively larger temperature parameter $\gamma$ to encourage exploration at the action level for RL.  {Furthermore, as the temperature parameter $\gamma$ $\to$ $0$, the variance of both the MFE and the MFO policies vanishes. Intuitively, the level of exploration diminishes and policies become increasingly deterministic. In the limiting case $\gamma =0$, they reduce to the solutions of the classical MFG and MFC problems.}
\end{Remark}

When the model parameters $b$, $\sigma$, $\Gamma$, and $\eta$ are unknown, we can accordingly consider the exact parameterized functions $J_d^\theta$ and $q_d^\psi$ that
\begin{align}
&\begin{aligned}
J_d^\theta(t, x, \mu) =& -\lambda e^{\theta_1(t-T)}(x-\bar\mu)^2 +  x + \frac{\gamma \theta_1}{4}(t-T)^2 + \theta_2 (t-T)\\
& + \frac{1}{4\lambda} e^{-\theta_1(t-T)} - \frac{1}{4\lambda},
\end{aligned} \label{parametrization_J_d} \\
&\begin{aligned}
q_d^{\psi}(t, x, \mu, a, {\bm h})  =& -\frac{\exp(\psi_1 + \psi_2(t- T))}{2} \big(a + \psi_3\big(x - \bar\mu) +\psi_4 e^{-\psi_2(t-T)}\big)^2 \\
&- \frac{\gamma}{2}\log (2\pi \gamma) + \frac{\gamma}{2} {\psi_1} + \frac{\psi_2 \gamma}{2}(t - T) - \psi_2(x-\bar\mu)\\
&+ \psi_5 e^{\psi_2(t - T)} (\bar\mu- x)\E_{\mu, {\bm h}}^e\big[a^{\bm h}\big],
\end{aligned} \label{parametrization_q_d}
\end{align}
where $\theta = (\theta_1, \theta_2)\trans \in \R^2$, and $\psi = (\psi_1, \ldots, \psi_5)\trans \in \R^5$.
It follows that
\[
{\bm \pi}^\psi(\cdot|t, x, \mu) = \Nc( -\psi_3\big(x - \bar\mu) - \psi_4\exp(-\psi_2(t-T))\big), \gamma\exp(-\psi_1-\psi_2(t-T))).
\]
We also parametrize the test policy ${\bm h}$ in the same form of ${\bm \pi}^\psi$ but with different parameters that
\begin{equation}\label{equ:MV-parametrizedh}
{\bm h}^{\tilde \psi}(\cdot|t, x, \mu) = \Nc(-\tilde \psi_3\big(x - \bar\mu) -\tilde\psi_4\exp(-\tilde \psi_2(t-T)), \gamma \exp(-\tilde\psi_1 -\tilde\psi_2(t-T))).
\end{equation}
It remains to require the approximators to satisfy the consistency condition  \eqref{mfg:param-q} for MFG and the consistency condition \eqref{mfc:param-q} for MFC, respectively. As a result, we need to impose $\psi_5 = -\frac{\psi_2}{\psi_4}$.

The true values of the parameters are given by
\begin{equation*}
\begin{aligned}
\theta_1^* &= \frac{b^2}{\sigma^2+\eta \Gamma^2}, \;\theta_2^* =  -\frac{\gamma}{2} \log \frac{\pi\gamma}{(\sigma^2+\eta\Gamma^2)\lambda}, \\
\psi_1^* &= \log (2\lambda (\sigma^2+\eta\Gamma^2)), \;\psi_2^* = \frac{b^2}{\sigma^2+\eta \Gamma^2}, \;\psi_3^* = \frac{b}{\sigma^2+\eta \Gamma^2},\; \psi_4^* = -\frac{b}{2(\sigma^2+\eta\Gamma^2)\lambda}.
\end{aligned}
\end{equation*}


{\bf The Simulator Environment$_{\Delta t}$.} Denote $\bar\mu_{\text{true}, t_k}$ and $\bar\mu_{\text{estimate}, t_k}$ as the true mean and the estimated mean by the representative agent of the distribution of $X_{\text{true}, t_k}$ under the policy ${\bm h}^{\tilde\psi}$, respectively. Then in the learning phase, for each time step $t_k$, the representative agent takes action $a_{t_k}$ $\sim {\bm h}^{\tilde\psi}(\cdot|t_k, X_{t_k}, \bar\mu_{\text{estimate}, t_k})$, and the state at $t_{k+1}$ is generated by
\begin{equation}\label{MV_representative_simulator}
X_{t_{k+1}} = X_{t_k} + a_{t_k}b \triangle t + a_{t_k}\sigma \sqrt{\Delta t}\triangle W + a_{t_k}(\Gamma\triangle N_{\eta,\Delta t}- \Gamma\eta \Delta t),
\end{equation}
where $\Delta t$ denotes the time step, $\triangle W\sim \Nc(0,1)$, and $\triangle N_{\eta,\Delta t} \sim {\rm Poisson}(\eta \Delta t)$.

{As a numerical safeguard during training, each sampled test-policy action is replaced by $\operatorname{clip}(a,-9,9)$ before the same realized action is used in both the simulator and the stochastic-gradient update. The analytical Gaussian population-action mean appearing in $q_d^\psi$ is kept unchanged, and no action clipping is applied in the independent held-out evaluation.}

Next, we apply Algorithm \ref{algo:decoupled-q} to the mean-variance portfolio optimization problem, with parametrizations depicted in \eqref{parametrization_J_d} and \eqref{parametrization_q_d}. For the simulator, we choose $T=1$, $b=0.25$, $\sigma=0.5$, $\Gamma=0.5$ and $\eta=1$. Moreover, we set the known parameters as $\gamma=0.5$, $\lambda=2$ and $\beta=0$. We also choose the time discretization $K=500$ (so that $\Delta t = 0.002$), the number of episodes $N=20000$, and the number of test policies $M=20$. The lower bounds and upper bounds of uniform distribution for the test policies (see Line 5 of Algorithm \ref{algo:decoupled-q}) are $p_1(j)=p_2(j)=p_3(j)=p_4(j)=0$, $q_1(j)=q_2(j)=q_3(j)=q_4(j)=\frac{2}{j^{0.05}}$. The learning rates are given by
\begin{equation*}
\begin{aligned}
&\alpha_\theta(j)=\left\{
\begin{aligned}
&\left( \frac{0.003}{j^{0.41}},\frac{0.02}{j^{0.31}}     \right),&1\leq j<15000,\\
&\left( \frac{0.001}{j^{0.41}},\frac{0.02}{3\times j^{0.31}} \right), &j\geq 15000,
\end{aligned}\right.\\
&\alpha_\psi(j) = \left\{
\begin{aligned}
&\left( \frac{0.005}{j^{0.31}},\frac{0.01}{j^{0.31}},\frac{0.01}{j^{0.31}},\frac{0.01}{j^{0.31}}           \right), &1\leq j<15000,\\
&\left( \frac{0.0025}{j^{0.31}},\frac{0.005}{j^{0.31}},\frac{0.005}{j^{0.31}},\frac{0.005}{j^{0.31}}           \right), &j\geq 15000.
\end{aligned}
\right.
\end{aligned}
\end{equation*}

\medskip

\begin{figure}[t!]
\centering
\begin{minipage}{0.98\textwidth}
\centering
\includegraphics[width=\linewidth]{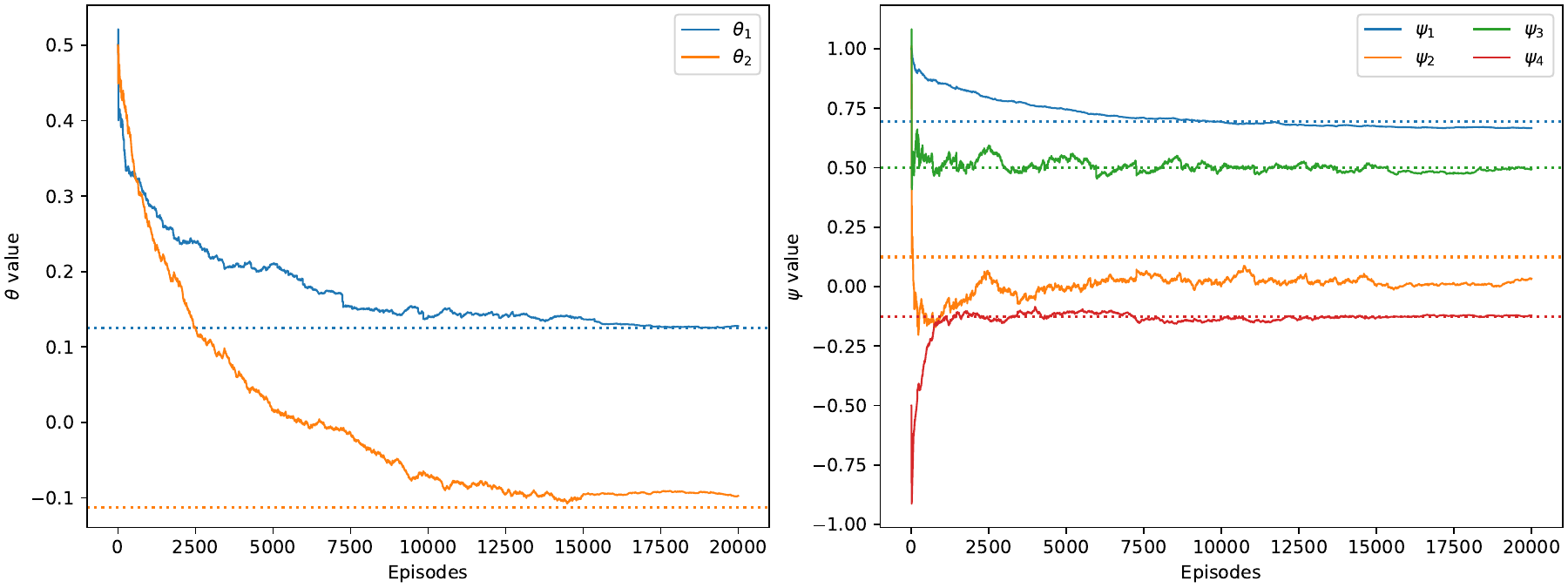}
\caption{Parameter trajectories. Dotted lines indicate the analytical parameter values.}
\label{MVfigure:parameters}
\end{minipage}
\vspace{0.8em}

\begin{minipage}{0.6\textwidth}
\centering
\includegraphics[width=\linewidth]{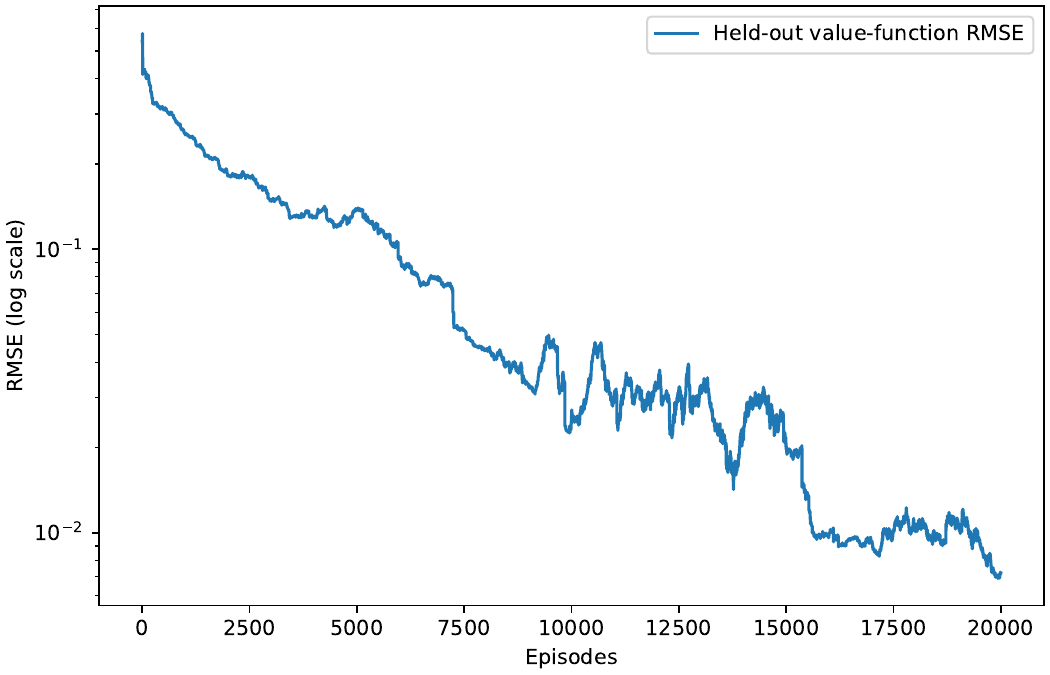}
\caption{{Held-out value-function RMSE over learning episodes. The vertical axis is logarithmic, and the same independently generated held-out sample is used at every episode.}}
\label{MVfigure:value}
\end{minipage}
\end{figure}

The initial state distribution is $\mu_0=\Nc(0,1)$, and the estimated population means are updated with $\rho^j=1/j$. We initialize the parameters at $\theta_0=(0.5,0.5)$ and $\psi_0=(1,1,1,-0.5)$. {These are generic values that are not centered at the analytical parameters. In particular, $\psi_{4,0}$ is initialized away from zero because the consistency condition requires $\psi_5=-\psi_2/\psi_4$. The coordinate-wise learning rates accommodate the different empirical scales of the parameter updates, and their coefficients are reduced from episode $15000$ onward to mitigate late-stage fluctuations.} The resulting parameter trajectories are presented in Figure \ref{MVfigure:parameters}, and the learnt and true parameter values are summarized in Table \ref{parameters_val}. At the final episode, the normalized errors $\|\hat\theta-\theta^*\|_2/\|\theta^*\|_2$ and $\|\hat\psi-\psi^*\|_2/\|\psi^*\|_2$ are $0.0935$ and $0.1092$, respectively.

\begin{table}[htbp]
\centering
\caption{Learnt versus true parameter values} \label{parameters_val}
\begin{tabular}{@{}r@{\quad}r@{\quad}r@{\quad}r@{\quad}r@{\quad}r@{\quad}r@{}}
\toprule
             & $\theta_1$& $\theta_2$ & $\psi_1$ & $\psi_2$ & $\psi_3$ & $\psi_4$ \\ \midrule
True Value   & 0.1250   & -0.1129    & 0.6931   & 0.1250   & 0.5000   & -0.1250   \\
Learnt Value & 0.1277   & -0.0974    & 0.6669   & 0.0336   & 0.4950   & -0.1214 \\
\end{tabular}
\end{table}

\begingroup
To evaluate the value-function approximation independently of the training trajectories, we generate a held-out population of $L=10000$ agents under the mean control induced by the optimal Gaussian policy. The agents have independent Brownian and Poisson increments, and the held-out sample is generated independently of the training trajectories. Denote the resulting states and their empirical mean by $X_{t_k}^{*,i}$ and $\bar\mu_{t_k}^*:=L^{-1}\sum_{i=1}^L X_{t_k}^{*,i}$, respectively. The same held-out inputs are used to evaluate every parameter iterate. For episode $j=0,\ldots,N$, we report
\begin{equation*}
\operatorname{RMSE}_J(j):=\left(\frac{1}{L(K+1)}\sum_{i=1}^L\sum_{k=0}^K\left|J_d^{\theta_j}(t_k,X_{t_k}^{*,i},\bar\mu_{t_k}^*)-J_d^{\theta^*}(t_k,X_{t_k}^{*,i},\bar\mu_{t_k}^*)\right|^2\right)^{1/2}.
\end{equation*}
Figure \ref{MVfigure:value} plots this held-out RMSE on a logarithmic vertical scale. In the reported run, the RMSE is $0.5433$ at initialization and $0.007223$ after $20000$ episodes.

Keeping all other settings fixed, we repeat the full $N=20000$-episode experiment with $M\in\{5,10,20,40\}$. Each value of $M$ is evaluated in a single training run under the same initialization and random-number-stream protocol, and all runs use common independently generated held-out evaluation inputs. Write
$E_\theta:=\|\hat\theta-\theta^*\|_2/\|\theta^*\|_2$ and $E_\psi:=\|\hat\psi-\psi^*\|_2/\|\psi^*\|_2$.

\begin{table}[htbp]

\centering
\caption{Sensitivity to the number of test policies $M$ in one fixed training realization.\label{tab:MV-M-sensitivity}}
{\begin{tabular}{@{}r@{\qquad}r@{\qquad}r@{\qquad}r@{}}
\toprule
$M$ & $E_\theta$ & $E_\psi$ & $\operatorname{RMSE}_J(N)$ \\ \hline
$5$  & $0.3672$ & $0.2489$ & $0.065394$ \\
$10$ & $0.0590$ & $0.1627$ & $0.012267$ \\
$20$ & $0.0935$ & $0.1092$ & $0.007223$ \\
$40$ & $0.0472$ & $0.1241$ & $0.013517$ \\
\end{tabular}}{}
\end{table}

Table \ref{tab:MV-M-sensitivity} shows that choosing $M$ too small can lead to large errors: all three error measures are substantially larger at $M=5$. Once $M$ is increased to moderate values such as $10$ or $20$, the errors decrease markedly and the algorithm already yields comparatively stable results in the reported experiment. Since each episode simulates and processes $M$ test policies, increasing $M$ further raises the computational cost. The choice $M=20$ used in the main experiment therefore provides a practical balance between numerical accuracy and computation. Across the four runs, clipping affects at most $0.000958\%$ of the training actions and is inactive after episode $112$.
\endgroup

\clearpage

{

\subsection{Crowd-aversion transport with jump-diffusion dynamics}

We next consider a non-LQ crowd aversion transport problem adapted from the crowd-aversion transport model in \cite{BHYZ2026}. In that work, the example is used to illustrate a continuous-time policy gradient method for a diffusion model. Here we keep the same target-seeking and crowd-avoidance mechanism, but make the modifications needed for the present q-learning framework: we work on a finite horizon, add compound Poisson shocks to the state dynamics, and learn the policy through the martingale characterization of the decoupled Iq-function rather than through a policy-gradient update. The MFC experiment is presented as the main illustration in this section; we also report the corresponding MFG run as an auxiliary comparison under the same numerical setting.

The state is two-dimensional. Under a randomized policy ${\bm \pi}$, the representative state follows the controlled jump-diffusion
\begin{equation}\label{sec72:state-continuous}
\dd X_t=a_tdt+\sigma \dd W_t+\dd J_t,\quad a_t\sim{\bm \pi}(\cdot|t,X_t,\mu_t),
\end{equation}
where $W$ is a two-dimensional Brownian motion, $J_t=\sum_{\ell=1}^{N_t}Z_\ell$ is a compound Poisson process with intensity $\lambda_J$ and marks $Z_\ell\sim\Nc(0,\sigma_J^2I_2)$, and $\mu_t=\Lc(X_t)$. In the discretization of the simulator, the time grid is $t_k=k\Delta t$ and \eqref{sec72:state-continuous} is implemented as
 \begin{equation*}
X_{t_{k+1}}=X_{t_k}+a_{t_k}\Delta t+\sigma\sqrt{\Delta t}\,\xi_k+\sum_{\ell=1}^{\Delta N_k}Z_{k,\ell},
\end{equation*}
where $\xi_k\sim\Nc(0,I_2)$ and $\Delta N_k\sim{\rm Poisson}(\lambda_J\Delta t)$. The population distribution is approximated by the empirical measure
\begin{equation*}
\widehat\mu_{t_k}^{N}=\frac1N\sum_{i=1}^N\delta_{X_{t_k}^i},
\end{equation*}
where the particles follow the same learned policy as the representative agent.

Following \cite{BHYZ2026}, the cost includes control effort, attraction to a target, crowd aversion, and a mild confinement term. In the present finite-horizon experiment, crowd aversion is imposed as a running cost. For a policy ${\bm \pi}$, we write
\begin{align*}
\mathcal{C}({\bm \pi})=
\E\bigg[\int_0^T e^{-\beta t}\left\{\frac{c_a}{2}|a_t|^2+w_CK_\ell(X_t,\mu_t)\right\}dt+e^{-\beta T}\left\{
\frac{w_T}{2}|X_T-x_{\rm tar}|^2
+\frac{w_R}{2}|X_T|^2\right\}\bigg],
\end{align*}
where
\begin{equation*}
K_\ell(x,\mu)=\int_{\R^2}\exp\left(-\frac{|y-x|^2}{2\ell^2}\right)\mu(\dd y).
\end{equation*}
Equivalently, the reward in \eqref{decoupled-J} is the negative of the running control and crowd costs, and the terminal reward is
\begin{equation}\label{sec72:terminal-reward}
g(x,\mu)=-\frac{w_T}{2}|x-x_{\rm tar}|^2-\frac{w_R}{2}|x|^2.
\end{equation}

We apply Algorithm \ref{algo:decoupled-q} without imposing an exact closed-form structure on the value function or the decoupled Iq-function. The value approximation enforces the terminal condition in \eqref{sec72:terminal-reward}, while a finite-rank parameterization of the decoupled and essential q-functions enforces the MFC integral and Gibbs-policy consistency conditions. The test policies are Gaussian mean perturbations of the current q-induced policy, and the empirical population distribution is updated across episodes according to \eqref{equ:distribution-update}. Detailed parameterizations and implementation formulas for both MFC and MFG problems are provided in the Online Companion.

\begin{figure}[!t]
\centering
\begin{subfigure}{0.48\textwidth}
    \centering
    \includegraphics[width=\linewidth]{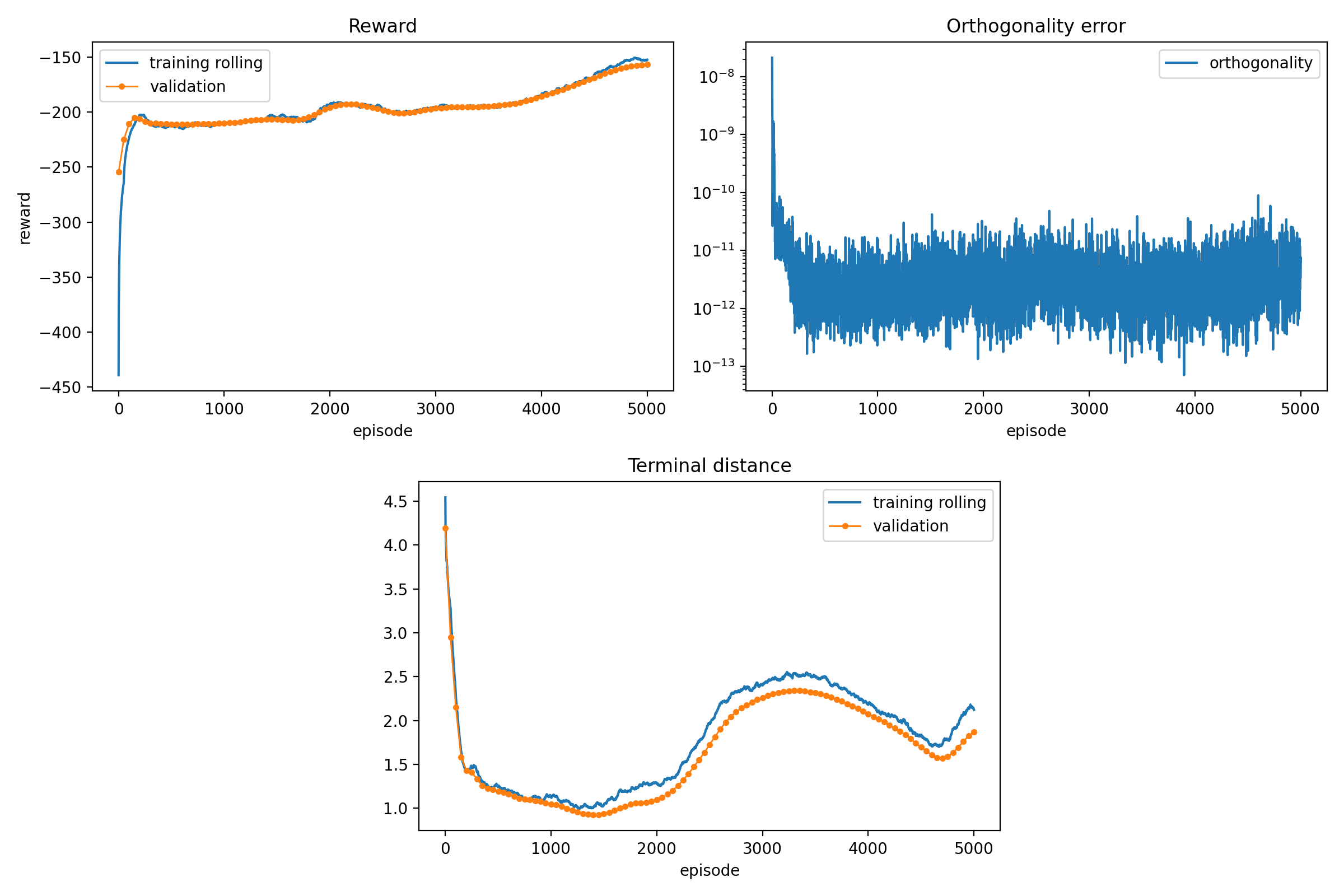}
    \caption{MFC}
\end{subfigure}\quad
\begin{subfigure}{0.48\textwidth}
    \centering
    \includegraphics[width=\linewidth]{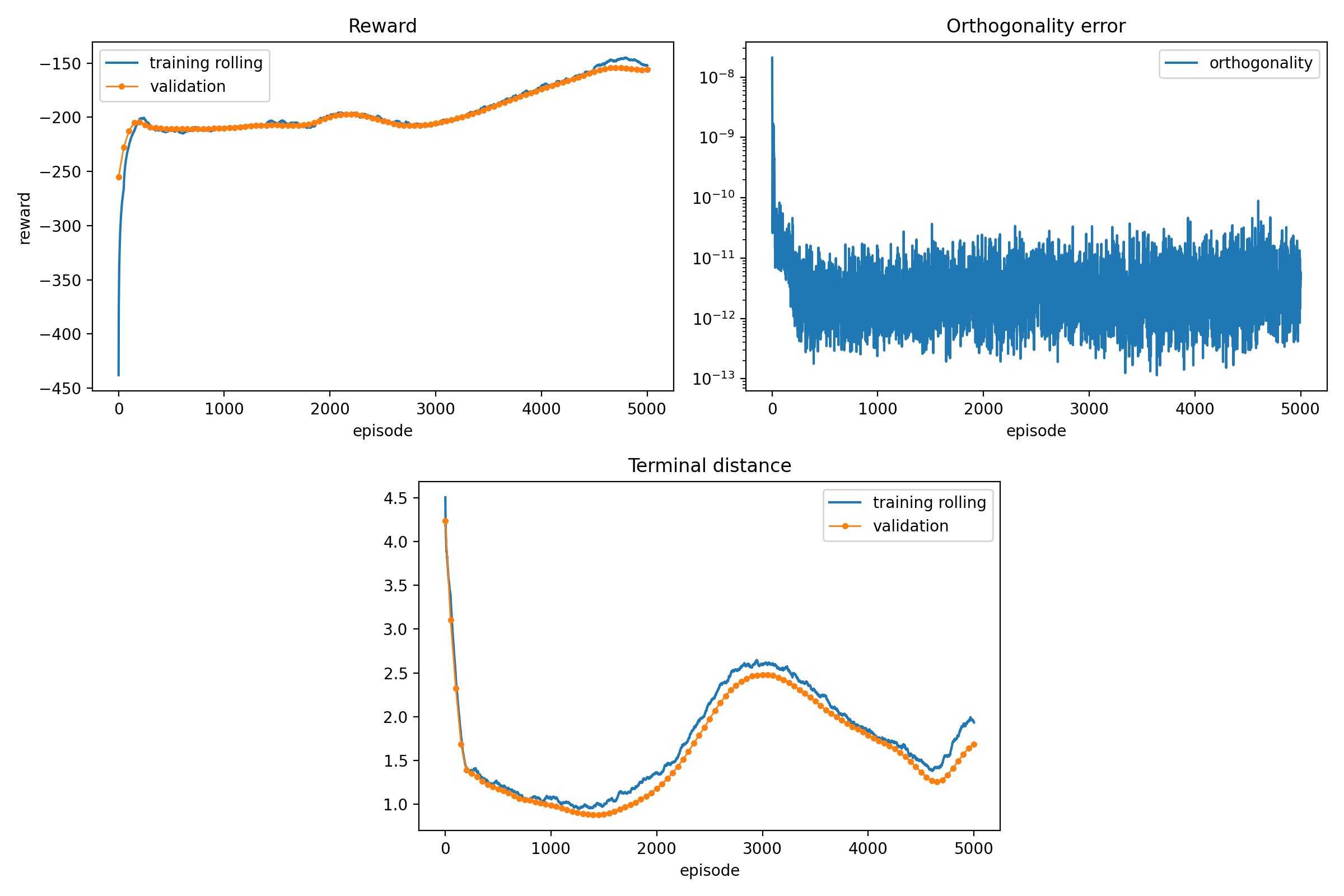}
    \caption{MFG}
\end{subfigure}
\caption{Training diagnostics for the crowd-aversion jump-diffusion example.}
\label{fig:section72-training}
\end{figure}

\begin{figure}[!t]
\centering
\begin{subfigure}{0.48\textwidth}
    \centering
    \includegraphics[width=\linewidth]{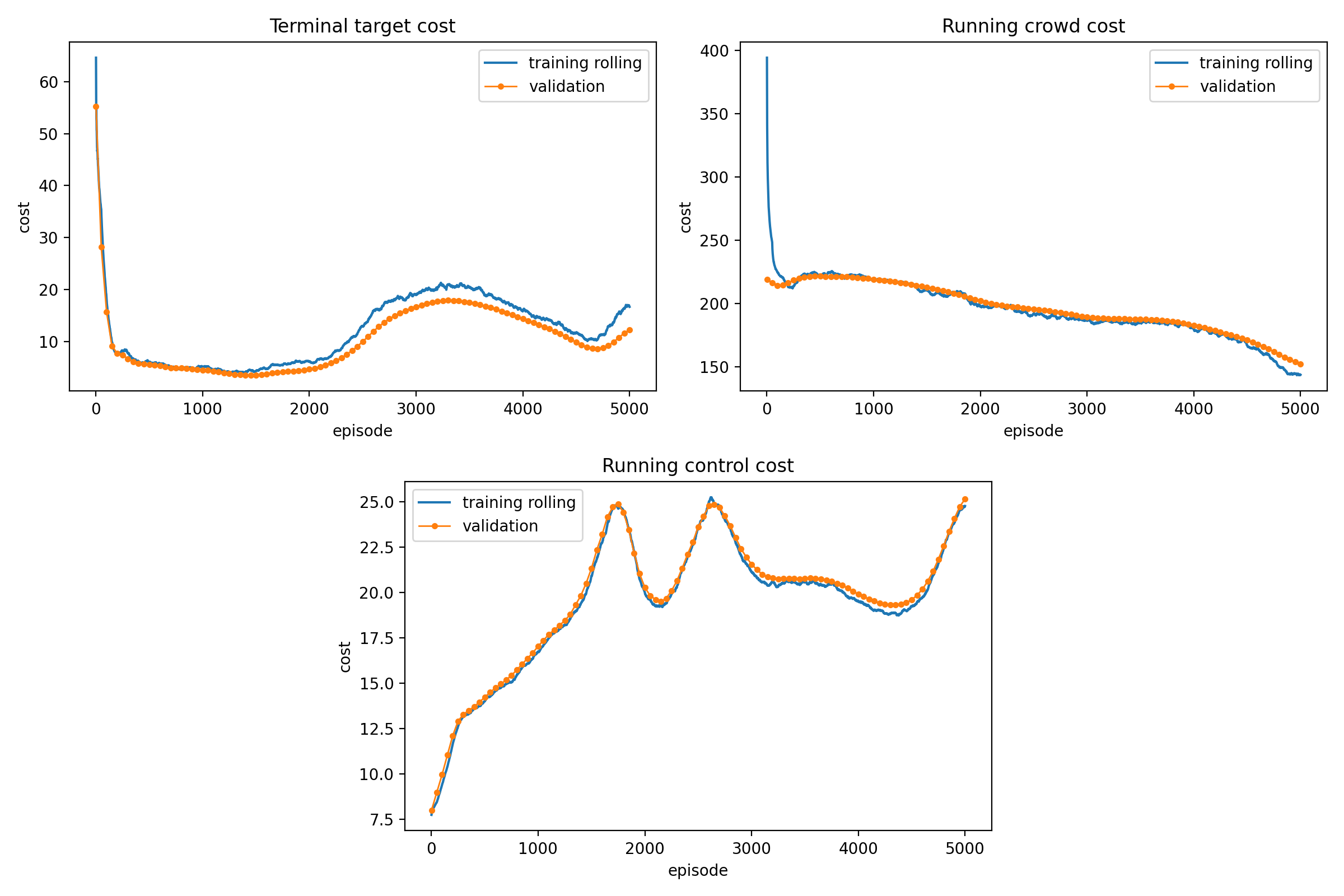}
    \caption{MFC}
\end{subfigure}\quad
\begin{subfigure}{0.48\textwidth}
    \centering
    \includegraphics[width=\linewidth]{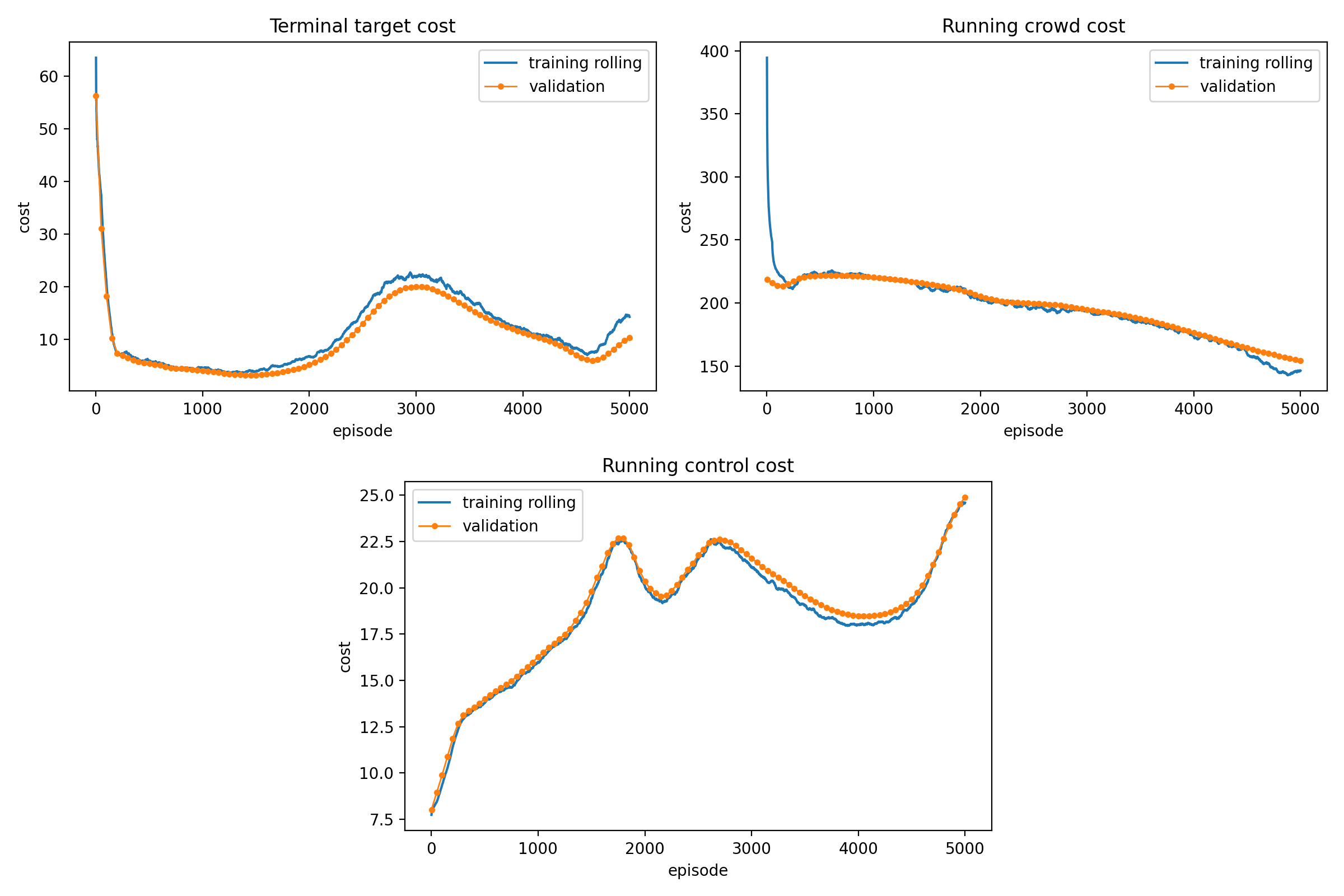}
    \caption{MFG}
\end{subfigure}
\caption{Objective decomposition under the learned randomized policies.}
\label{fig:section72-objective}
\end{figure}

\begin{figure}[!t]
\centering
\begin{subfigure}{0.48\textwidth}
    \centering
    \includegraphics[width=\linewidth]{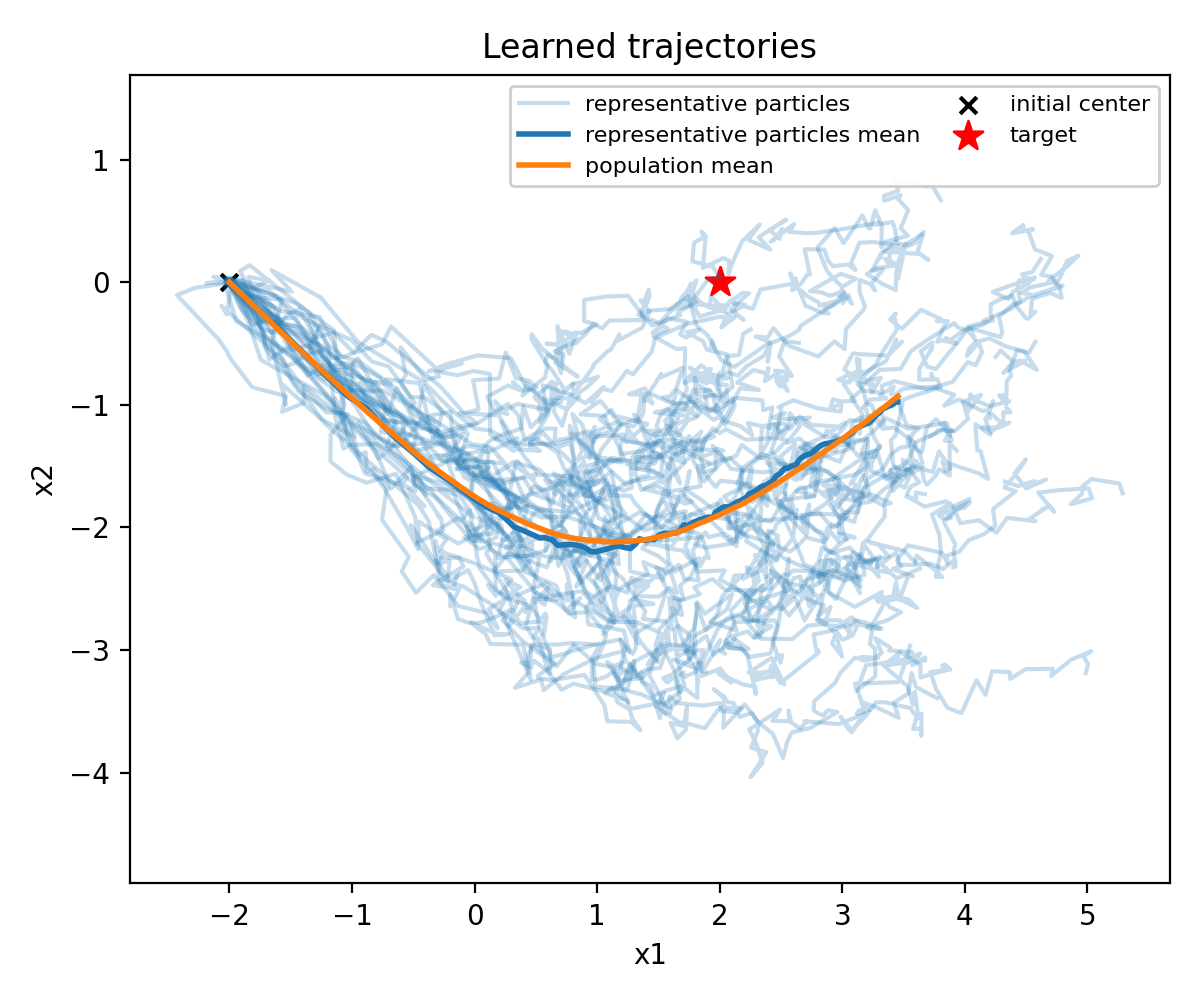}
    \caption{MFC}
\end{subfigure}\quad
\begin{subfigure}{0.48\textwidth}
    \centering
    \includegraphics[width=\linewidth]{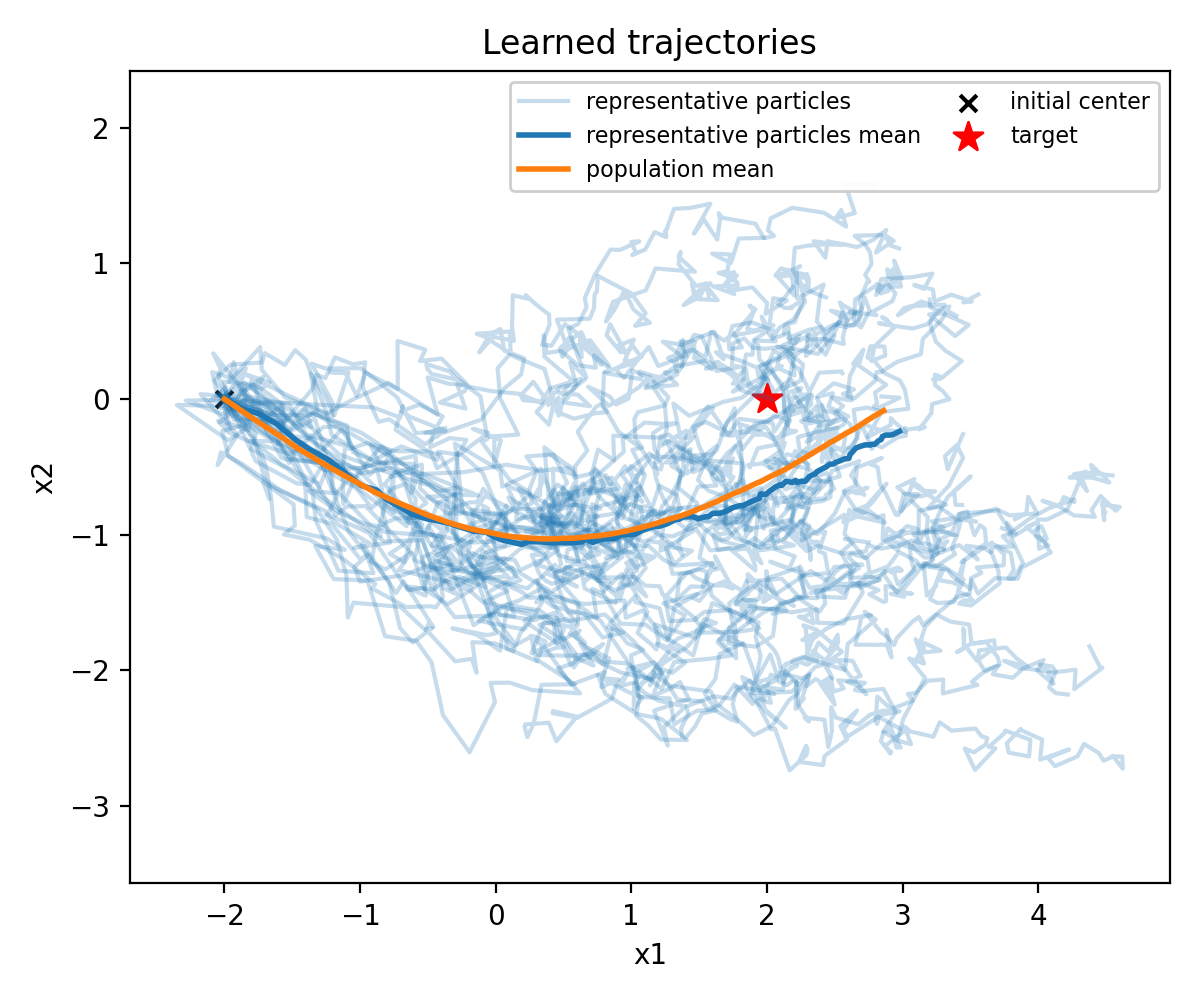}
    \caption{MFG}
\end{subfigure}
\caption{Mean paths and sample paths generated by the learned randomized policies.}
\label{fig:section72-trajectories}
\end{figure}

The main numerical configuration is $T=5$, $\Delta t=0.05$, $\sigma=0.2$, $\beta=0.1$, $x_0=(-2,0)$, $x_{\rm tar}=(2,0)$, $\lambda_J=0.5$, and $\sigma_J=0.1$. We use 64 empirical population particles, 20 test policies, and 5000 training episodes. The cost parameters are $c_a=1$, $w_T=10$, $w_C=100$, $w_R=0.1$, $\ell=0.8$, and $\gamma=2$. Validation is performed every 50 episodes. For both the MFC run and the auxiliary MFG run, the checkpoint with the highest validation regularized reward is used for evaluation. The relative diagnostic gain in Table \ref{tab:section72-policy-comparison} is computed at the selected checkpoint.

The running crowd penalty is intentionally much larger than the terminal target penalty in order to test whether the learned randomized policy responds to crowding while still moving the population toward the target. Figure \ref{fig:section72-training} shows a sharp early improvement in reward and terminal distance, followed by a stable but fluctuating regime, while the orthogonality error remains small after the initial transient. Figure \ref{fig:section72-objective} displays the resulting tradeoff among terminal target cost, running crowd cost, and running control cost.

\begin{table}[H]
\centering
\caption{Performance of the learned randomized policies and the no-control benchmark.\label{tab:section72-policy-comparison}}
{\begin{tabular}{@{}l@{\quad}r@{\quad}r@{\quad}r@{}}
\toprule
Performance measure & MFC & MFG & No control \\ \hline
Unscaled entropy-regularized reward & $-160.06$ & $-162.92$ & $-350.34$ \\
Mean terminal distance & $2.04$ & $1.54$ & $4.01$ \\
Discounted terminal target cost & $16.32$ & $8.79$ & $49.61$ \\
Discounted running crowd cost & $150.83$ & $167.85$ & $300.60$ \\
Discounted running control cost & $25.70$ & $18.81$ & $0.00$ \\
Relative diagnostic gain & $1.83\%$ & $0.00\%$ & --
\end{tabular}}
{The MFC and MFG columns are averaged over 64 independent randomized-policy evaluation rollouts. The no-control benchmark uses $a_t\equiv0$ for both the representative agent and the population particles under the same evaluation seed.}
\end{table}

Table \ref{tab:section72-policy-comparison} compares the learned policies with the deterministic no-control policy. Both learned policies substantially improve the entropy-regularized reward and the terminal accuracy relative to the no-control benchmark: the reward improves from $-350.34$ to $-160.06$ in the MFC run and to $-162.92$ in the MFG run, while the mean terminal distance decreases from $4.01$ to $2.04$ and $1.54$, respectively. The two columns also show the expected difference between the control problem and the game. The MFC policy pays a larger control cost, $25.70$ instead of $18.81$, but achieves a lower running crowd cost, $150.83$ instead of $167.85$, which is consistent with a social planner internalizing the congestion externality. The MFG policy places relatively more weight on individual target accuracy against a fixed population flow, reaching a smaller terminal target cost, $8.79$ instead of $16.32$, while accepting a larger realized crowding cost. The small relative diagnostic gains, $1.83\%$ for MFC and $0.00\%$ for MFG, indicate that the selected checkpoints are locally stable under the tested perturbations. Figure \ref{fig:section72-trajectories} provides the corresponding pathwise picture: under both learned randomized policies, the sample paths and population mean move from the initial cluster toward the target while bending around the crowded region.

}

\section{Proofs}\label{sec:proofs}

This section collects the proofs of all results in previous sections.
{
\subsection{Proof of Proposition \ref{prop:regularity}}\label{app:proof-prop-regularity}

Noting that we are interested in the regularity of value function (with a fixed pair of policies) defined in \eqref{decoupled-J}, it is sufficient for our purpose to work in the {\it weak formulation}. In other words, we consider other stochastic processes, possibly defined on a different probability space, with different filtrations, as long as they have the same distribution as in the main context. In this context, $\P$ and $\E$ in this section are generic probability notations, which may be different from $\P$, $\hat\P$ or $\P^e$ in the main context. We also use $\bar \P$ ($\bar \E$) to denote the probability (expectation, respectively) with respect to independent copies of random variables that are originally under $\P$ ($\E$), and then $\tilde \P$ ($\tilde\E$), $\hat \P$ ($\hat \E$) if we need more independent copies.
%

\begin{proof}[Proof of Proposition \ref{prop:regularity}]
For $\bh\in \Pi$ and any test function $v\in C_b^2(\R^d)$, we have
\begin{equation*}
\begin{aligned}
&\int_{\Ac}\Lc^{t, a, \mu}[v(\cdot)](x) { \bh}(a|t,x,\mu)\dd a\\
=& b_{\bh}(t, x, \mu)\trans \partial_x v(x) + \frac{1}{2}{\rm Tr}\big(\sigma_{\bh}\sigma_{\bh}\trans(t, x, \mu) \partial_{xx} v(x)\big) \\
&+ \int_{\Ac}\big[v(x + \Gamma(t, x, a, \mu)) - v(x) - \Gamma(t, x, a, \mu)\trans \partial_x v(x)\big]\eta(t){\bh}(a|t,x,\mu)\dd a\\
=&b_{\bh}(t, x, \mu)\trans \partial_x v(x) + \frac{1}{2}{\rm Tr}\big(\sigma_{\bh}\sigma_{\bh}\trans(t, x, \mu) \partial_{xx} v(x)\big) \\
&+\int_{\R^d}\big[v(x + y) - v(x) - y\trans \partial_x v(x)\big]\eta(t)\Gamma(t,x,\mu,\cdot)_\#{\bh}(\cdot,t,x,\mu)(\dd y)\\
=& b_{\bh}(t, x, \mu)\trans \partial_x v(x) + \frac{1}{2}{\rm Tr}\big(\sigma_{\bh}\sigma_{\bh}\trans(t, x, \mu) \partial_{xx} v(x)\big) \\
&+\int_{\R^d}\big[v(x + y) - v(x) - y\trans \partial_x v(x)\big]\eta(t)\Gamma_{\bh}(t,x,\mu,\cdot)_\# \mu_0(\dd y)\\
=& b_{\bh}(t, x, \mu)\trans \partial_x v(x) + \frac{1}{2}{\rm Tr}\big(\sigma_{\bh}\sigma_{\bh}\trans(t, x, \mu) \partial_{xx} v(x)\big) \\
&+\int_{\R}\big[v(x + \Gamma_{\bh}(t,x,\mu,z)) - v(x) - \Gamma_{\bh}(t,x,\mu,z)\trans \partial_x v(x)\big]\eta(t)\mu_0(\dd z).
\end{aligned}
\end{equation*}
By the martingale problem characterization of SDE (see Remark \ref{rmk:martingale-problem}), we have constructed another pair of processes $(\hat{X}^{t,x,\mu, \widehat{\bm \pi}}, \hat{X}^{t,\xi, \bpi})$ equal to $({X}^{t,x,\mu, \widehat{\bm \pi}}, {X}^{t,\xi, {\bm \pi}})$ in law, satisfying the following mean-field SDE with jumps:
\begin{align}
\mathrm{d} \hat X_s^{t,\xi,\boldsymbol{\pi}}
  &= b_{\boldsymbol{\pi}}\bigl(s,\hat X_s^{t,\xi,\boldsymbol{\pi}},\mu_s^{\boldsymbol{\pi}}\bigr)\,\mathrm{d}s
   + \sigma_{\boldsymbol{\pi}}\bigl(s,\hat X_s^{t,\xi,\boldsymbol{\pi}},\mu_s^{\boldsymbol{\pi}}\bigr)\,\mathrm{d}W_s \notag \\
  &\quad + \int_{\mathbb{R}}
      \Gamma_{\boldsymbol{\pi}}\bigl(s,\hat X_{s-}^{t,\xi,\boldsymbol{\pi}},\mu_s^{\boldsymbol{\pi}},z\bigr)
      \,\widetilde{\mathcal{N}}(\mathrm{d}s,\mathrm{d}z),
\label{Xhat:x}
\\[0.5em]
\mathrm{d} \hat X_s^{t,x,\mu,\widehat{\boldsymbol{\pi}}}
  &= b_{\widehat{\boldsymbol{\pi}}}\bigl(s,\hat X_s^{t,x,\mu,\widehat{\boldsymbol{\pi}}},\mu_s^{\boldsymbol{\pi}}\bigr)\,\mathrm{d}s
   + \sigma_{\widehat{\boldsymbol{\pi}}}\bigl(s,\hat X_s^{t,x,\mu,\widehat{\boldsymbol{\pi}}},\mu_s^{\boldsymbol{\pi}}\bigr)\,\mathrm{d}W_s \notag \\
  &\quad + \int_{\mathbb{R}}
      \Gamma_{\widehat{\boldsymbol{\pi}}}\bigl(s,\hat X_{s-}^{t,x,\mu,\widehat{\boldsymbol{\pi}}},\mu_s^{\boldsymbol{\pi}},z\bigr)
      \,\widetilde{\mathcal{N}}(\mathrm{d}s,\mathrm{d}z).
\label{Xhat:xi}
\end{align}
where $\hat X_{t-}^{t,\xi, {\bm \pi}}\sim \mu$, $\hat X_{t-}^{t,x,\mu, \widehat{\bm \pi}}=x, $ $\mu_s^{\bm \pi} = \P_{\hat X_s^{t, \xi, {\bm \pi}}}$, and $\tilde \Nc = \Nc(\dd s,dz) - \eta(s) \dd s \mu_0(dz)$ is a compensated random measure with compensator $\mu_0(z)\eta(s) \dd s dz$. In the sequel, we omit all superscripts of $\bpi$ and $\widehat \bpi$. We also drop the hat above $X$, and directly write $(X^{t,x,\mu},X^{t,\xi})$. We shall make use of properties of $(X^{t,x,\mu},X^{t,\xi})$, as well as their derivatives with respect to $\partial_x$, $\partial_\mu$ and $\partial_y\partial_\mu$, which we summarize in Lemmas \ref{app:lemma:Xbound}-\ref{app:lemma:X_mubound}.

Recall that the decoupled value function is given by~\eqref{decoupled-J}.
For simplicity, we denote
$f(t,x,\mu) := e^{-\beta t}\big(r_{\widehat \bpi}(t,x,\mu) - \gamma \Ec_{\widehat{\bpi}}(t,x,\mu)\big)$. To show the regularity of $J_d$, we only showcase the differentiability with respect to $\partial_\mu$ and prove that it has the following explicit form:
\begin{align}
\partial_\mu J_d(t,x,\mu,y) =&e^{\beta t} \E \bigg[e^{-\beta T}\big( \partial_\mu X^{t,x,\mu}_T(y)\big)\trans \partial_x g(X^{t,x,\mu}_T,\mu_T)\nonumber\\
&+ e^{-\beta T}\bar \E\big[ (\partial_x \bar X^{t,y,\mu}_T )\trans \partial_\mu g (X^{t,x,\mu}_T,\mu_T,\bar X^{t,y,\mu}_T) \nonumber\\
&+ \big(\partial_\mu \bar X^{t,\xi}_T(y)\big)\trans \partial_\mu g(X^{t,x,\mu}_T,\mu_T,\bar X^{t,\xi}_T)\big]\nonumber\\
&+ \int_s^T\Big( \big( \partial_\mu X^{t,x,\mu}_r(y) \big)\trans \partial_x f(r, X^{t,x,\mu}_r,\mu_r)\nonumber \\
&+ \bar \E \big[( \partial_x \bar X^{t,y,\mu}_r)\trans \partial_\mu f(r,X^{t,x,\mu}_r,\mu_r)+\big(\partial_\mu \bar X^{t,\xi}_r\big)\trans \partial_\mu f(r,X^{t,x,\mu}_r,\mu_r, \bar X^{t,\xi}_r)\big]\Big)\dd r\bigg] \label{partial_mu_Jd}
\end{align}
For detailed derivations of other derivatives of $J_d$, we refer readers to \cite{FGLPS23} for the case without jumps, and to \cite{Li18} for a more general case with jumps.

To prove \eqref{partial_mu_Jd}, we consider the {\it lifted functional} $\tilde J_d(t,x,\xi):=J_d(t,x,\P_\xi)$ and the {\it lifted process} $\tilde X^{t,x,\xi}:= X^{t,x,\P_\xi}$. We simply write them as $J_d(t,x,\xi)$ and $X^{t,x,\xi}$. By the definition of L-derivative (c.f., Section 5.2 of \cite{CarD}), we need to verify
$|J(t,x,\xi') - J(t,x,\xi) - \tilde \E[\partial_\mu J_d(t,x,\mu,\tilde \xi)(\tilde \xi'-\tilde \xi)]| = o(\|\xi'-\xi\|_{L^2})$ where $\xi' \in L^2_{\Fc_t}$ and $\|\xi'-\xi\|_{L^2}\to 0$. To this end, denote $\{\mu'_s\}_{t\leq s \leq T} := \{\P_{X^{t,\xi'}_s}\}_{t\leq s \leq T}$, and introduce the following notations:
\begin{equation*}
\begin{aligned}
 \Xi = &e^{-rT} g ({X}^{t,x,\xi}_T,\mu_T),\quad \Xi' =   e^{-rT} g ({X}^{t,x,\xi'}_T,\mu'_T)\\
 \partial \Xi (y)= &e^{-rT}\big( \partial_\mu X^{t,x,\xi}_T(y)\big)\trans \partial_x g(X^{t,x,\xi}_T,\mu_T)\\
&+ e^{-rT}\bar \E\big[ (\partial_x \bar X^{t,y,\xi}_T )\trans \partial_\mu g (X^{t,x,\xi}_T,\mu_T,\bar X^{t,y,\xi}_T) + \big(\partial_\mu \bar X^{t,\xi}_T(y)\big)\trans \partial_\mu g(X^{t,x,\xi}_T,\mu_T,\bar X^{t,\xi}_T)\big],\\
\sF_r =& f(r, X^{t,x,\xi}_r,\mu_r), \quad \sF_r' = f(r, X^{t,x,\xi'}_r,\mu'_r),\\
\partial \sF_r(y) =& \big( \partial_\mu X^{t,x,\mu}_r(y) \big)\trans \partial_x f(r, X^{t,x,\mu}_r,\mu_r) \\
&+ \bar \E \big[( \partial_x \bar X^{t,y,\mu}_r)\trans \partial_\mu f(r,X^{t,x,\mu}_r,\mu_r)+\big(\partial_\mu \bar X^{t,\xi}_r\big)\trans \partial_\mu f(r,X^{t,x,\mu}_r,\mu_r, \bar X^{t,\xi}_r)\big].
\end{aligned}
\end{equation*}
We first establish the estimate regarding $\Xi$:
\begin{equation*}
\begin{aligned}
&\Big|\E\Big[ \Xi'-\Xi - \tilde\E[\partial \Xi(\tilde \xi) (\tilde \xi' - \tilde \xi)] \Big] \Big|\\
\leq  & \underbrace{\E \Big[ \big| g(X^{t,x,\xi'}_T,\mu'_T) - g(X^{t,x,\xi}_T,\mu'_T) - \tilde \E\big[ \big(\partial_\mu X^{t,x,\xi}_T(\tilde \xi)\big)\trans \partial_x g (X^{t,x,\xi}_T,\mu_T)\cdot (\tilde \xi'-\tilde \xi)\big] \Big|\Big]}_{\text{(I)}}\\
&+ \underbrace{\E\Big[ \Big|g(X^{t,x,\xi}_T,\mu_T') - g(X^{t,x,\xi}_T,\mu_T) - \bar \E[\partial_\mu g (X^{t,x,\xi}_T,\mu_T,\bar X^{t,\xi}_T) \cdot (\bar X^{t,\xi'}_T - \bar X^{t,\xi}_T)] \Big|\Big]}_{\text{(II)}} \\
&+ \E \Big[ \Big| \tilde \E\big[\bar \E[\partial_\mu g (X^{t,x,\xi}_T,\mu_T,\bar X^{t,\tilde \xi, \xi}_T) \cdot(\bar X^{t,\tilde \xi', \xi'}_T - \bar X^{t,\tilde \xi, \xi'}_T) ]\big] \\
&\underbrace{\ \ \ \ \ \ \  -\tilde \E\big[\bar \E[ (\partial_x \bar X^{t,\tilde \xi,\xi}_T )\trans \partial_\mu g (X^{t,x,\xi}_T,\mu_T,\bar X^{t,\tilde \xi,\xi}_T)\cdot(\tilde \xi' - \tilde \xi)  ]\big]  \Big| \Big]}_{\text{(III)}}\\
&+ \E\bigg[ \bigg| \hat \E \big[ \bar \E[\partial_\mu g (X^{t,x,\xi}_T,\mu_T,\bar X^{t,\hat \xi, \xi}_T) \cdot(\bar X^{t,\hat \xi, \xi'}_T - \bar X^{t,\hat \xi, \xi}_T)]\big]\\
&\underbrace{\ \ \ \ \ \ \  -\hat \E\bigg[ \tilde\E\Big[\bar \E\big[\big(\partial_\mu \bar X^{t,\hat \xi, \xi}_T(\tilde \xi)\big)\trans \partial_\mu g(X^{t,x,\xi}_T,\mu_T,\bar X^{t,\hat \xi, \xi}_T)\cdot(\tilde \xi' - \tilde \xi)  \big]\Big] \bigg]  \bigg|\bigg]. }_{\text{(IV)}}
\end{aligned}
\end{equation*}
For (I), we use fundamental theorem of calculus, Assumption \ref{assump}-(iii) and Lemmas \ref{app:lemma:Xbound} to obtain
\begin{align}
|\text{(I)}|=& \E \bigg[ \bigg| \int_0^1 \partial_x g(\lambda X^{t,x,\xi'}_T+(1-\lambda)X^{t,x,\xi}_T,\mu'_T)\cdot (  X^{t,x,\xi'}_T - X^{t,x,\xi}_T)d \lambda \nonumber\\
&- \tilde \E\big[ \big(\partial_\mu X^{t,x,\xi}_T(\tilde \xi)\big)\trans \partial_x g (X^{t,x,\xi}_T,\mu_T)\cdot (\tilde \xi'-\tilde \xi)\big]\bigg|\bigg]\nonumber\\
\leq& \int_0^1 \E\big[\big|\partial_x  g\big(\lambda X^{t,x,\xi'}_T+(1-\lambda)X^{t,x,\xi}_T,\mu'_T\big) - \partial_x g (X^{t,x,\xi}_T,\mu_T)|\cdot|X^{t,x,\xi}_T - X^{t,x,\xi'}_T\big| \big] d \lambda\nonumber\\
&+\E[|\partial_x g (X^{t,x,\xi}_T,\mu_T)|\cdot |X^{t,x,\xi'}_T - X^{t,x,\xi}_T - \bar \E [\partial_\mu X^{t,x,\xi}_T(\bar \xi)\cdot (\bar \xi' - \bar \xi)] | ]\nonumber\\
\leq & \bigg(\int_0^1  \Big(\E\big[ \big|\partial_x  g\big(\lambda X^{t,x,\xi'}_T+(1-\lambda)X^{t,x,\xi}_T,\mu'_T\big) - \partial_x g \big(X^{t,x,\xi}_T,\mu_T\big)\big|^2\big]\Big)^{1/2} d\lambda\bigg)\|X^{t,x,\xi'}_T - X^{t,x,\xi}_T \|_{L^2}\nonumber\\
&+\|\partial_x g (X^{t,x,\xi}_T,\mu_T)\|_{L^2}\|\| X^{t,x,\xi'}_T - X^{t,x,\xi}_T- \bar \E [\partial_\mu X^{t,x,\xi}_T(\bar \xi)\cdot (\bar \xi' - \bar \xi)]\|_{L^2}\nonumber\\
&\leq C \big(\|X^{t,x,\xi}_T - X^{t,x,\xi'}_T\|_{L^2}+ \Wc_2(\mu_T,\mu_T')\big)\|X^{t,x,\xi}_T - X^{t,x,\xi'}_T\|_{L^2} + o(\|\xi'-\xi\|_{L^2})\nonumber\\
&\leq O(\|\xi'-\xi\|_{L^2}^2) + o(\|\xi'-\xi\|_{L^2})\nonumber\\
& = o(\|\xi'-\xi\|_{L^2}).\label{J_mu_I}
\end{align}
For (II), we note that
\begin{equation*}
\begin{aligned}
&g(X^{t,x,\xi}_T,\mu'_T)-g(X^{t,x,\xi}_T,\mu_T)\\=& \int_0^1\int_{\R^d} \frac{\delta g}{\delta \mu}(X^{t,x,\xi}_T, \lambda \mu'_T+(1-\lambda)\mu_T,y) (\mu'_T - \mu_T)(\dd y)d\lambda\\
=& \int_0^1\bar \E\bigg[\frac{\delta g}{\delta \mu}(X^{t,x,\xi}_T, \lambda \mu'_T+(1-\lambda)\mu_T,\bar X^{t,\xi'}_T)-\frac{\delta g}{\delta \mu}(X^{t,x,\xi}_T, \lambda \mu'_T+(1-\lambda)\mu_T,\bar X^{t,\xi}_T)\bigg]\\
=&\int_0^1\int_0^1 \bar \E\big[ \partial_\mu g \big(X^{t,x,\xi}_T, \lambda \mu'_T+(1-\lambda)\mu_T,\theta  \bar X^{t,\xi'}_T+(1-\theta)\bar X^{t,\xi}_T\big)(\bar X^{t,\xi'}_T - \bar X^{t,\xi}_T)\big]\dd \lambda \dd \theta.
\end{aligned}
\end{equation*}
Therefore,
\begin{align}
|\text{(II)}| \leq& \bar \E\bigg[\bigg(\int_0^1\int_0^1  \E\big[\big|\partial_\mu g \big(X^{t,x,\xi}_T, \lambda \mu'_T+(1-\lambda)\mu_T,\theta  \bar X^{t,\xi'}_T+(1-\theta)\bar X^{t,\xi}_T\big) \nonumber\\
&- \partial_\mu g (X^{t,x,\xi}_T,\mu_T,\bar X^{t,\xi}_T)\big|\big] \dd \lambda \dd \theta \bigg)|\bar X^{t,\xi'}_T - \bar X^{t,\xi}_T|\bigg]\nonumber\\
\leq & C\big(\Wc_2(\mu'_T,\mu_T) + \|\bar X^{t,\xi'}_T - \bar X^{t,\xi}_T\|_{L^2}\big)\|\bar X^{t,\xi'}_T - \bar X^{t,\xi}_T\|_{L^2}\nonumber\\
\leq & C\|\bar X^{t,\xi'}_T - \bar X^{t,\xi}_T\|_{L^2}^2\nonumber\\
\leq & o(\|\xi' - \xi\|_{L^2}). \label{J_mu_II}
\end{align}
For (III), we simplify the notation by $\mathring \E:  = \E^{ \P\times \bar \P }$ and observe that
\begin{align}
\text{(III)} \leq & \tilde \E\big[ \mathring \E [|\partial_\mu g (X^{t,x,\xi}_T,\mu_T,\bar X^{t,\tilde \xi,\xi}_T) |\cdot |\bar X^{t,\tilde \xi',\xi'}_T - \bar X^{t,\tilde \xi, \xi'} - \partial_x \bar X^{t,\tilde \xi, \xi'}_T\cdot(\tilde \xi' - \tilde \xi)|] \big]\nonumber\\
&+ \tilde \E\big[ \mathring \E [|\partial_\mu g (X^{t,x,\xi}_T,\mu_T,\bar X^{t,\tilde \xi,\xi}_T) |\cdot |(\partial_x \bar X^{t,\tilde \xi, \xi}_T - \partial_x \bar X^{t,\tilde \xi, \xi'}_T)\cdot(\tilde \xi' - \tilde \xi)|] \big]\nonumber\\
\leq & C\tilde \E[ (1+|\tilde \xi|)\|\bar X^{t,\tilde \xi',\xi'}_T - \bar X^{t,\tilde \xi, \xi'} _T- \partial_x \bar X^{t,\tilde \xi, \xi'}_T\cdot(\tilde \xi' - \tilde \xi)\|_{L^2}]\nonumber\\
&+\tilde \E[ (1+|\tilde \xi|)|\tilde \xi'- \tilde\xi|\| \partial_x \bar X^{t,\tilde \xi, \xi}_T - \partial_x \bar X^{t,\tilde \xi, \xi'}_T\|_{L^2}]\nonumber\\
\leq & C\tilde \E\bigg[(1+|\tilde \xi|)|\tilde\xi'-\tilde \xi| \int_0^1 \|\partial_x \bar X^{t,\lambda \tilde \xi' + (1-\lambda) \tilde \xi,\xi'}_T - \partial_x\bar X^{t,\tilde \xi, \xi'}_T\|_{L^2}\dd \lambda  \bigg]\nonumber\\
&+C\tilde \E[ (1+|\tilde \xi|)|\tilde \xi'- \tilde\xi|\| \partial_x \bar X^{t,\tilde \xi, \xi}_T - \partial_x \bar X^{t,\tilde \xi, \xi'}_T\|_{L^2}]\nonumber\\
\leq & \tilde C\E[ (1+|\tilde \xi|)|\tilde \xi' - \tilde \xi| (|\tilde \xi' - \tilde \xi|\wedge K)]+C\tilde \E[(1+|\tilde \xi|)|\tilde \xi' - \tilde \xi|]\|\xi' - \xi\|_{L^2}.\label{J_mu_III_i}
\end{align}
Here in the last inequality, we have used
$\int_0^1 \|\partial_x \bar X^{t,\lambda \tilde \xi' + (1-\lambda) \tilde \xi,\xi'}_T - \partial_x\bar X^{t,\tilde \xi, \xi'}_T\|_{L^2}\dd \lambda\leq |\tilde \xi' - \tilde \xi| \wedge K$, which is a direct consequence of Lemma \ref{app:lemma:X_xbound}. Noting that
\begin{align*}
\tilde \E[(1+|\tilde \xi|)|\tilde \xi' - \tilde \xi|(|\tilde \xi' - \tilde \xi|\wedge K)]\leq &  \E[(1+|\tilde \xi|)|\tilde \xi' - \tilde \xi|(|\tilde \xi' - \tilde \xi|\wedge K)I_{\{1+|\tilde \xi|\leq \|\xi' - \xi\|_{L^2}^{-1/2}\}}]\\
&+ \E[(1+|\tilde \xi|)|\tilde \xi' - \tilde \xi|(|\tilde \xi' - \tilde \xi|\wedge K)I_{\{1+|\tilde \xi|> \|\xi' - \xi\|_{L^2}^{-1/2}\}}]\\
\leq & \|\xi' - \xi\|_{L^2}^{3/2} + C\big( \E[(1+|\tilde \xi|^2) I_{\{1+|\tilde \xi|> \|\xi' - \xi\|_{L^2}^{-1/2}\}}]\big)^{1/2}\|\xi'- \xi\|_{L^2}\\
\leq &o(\|\xi'-\xi\|_{L^2}).
\end{align*}
Hence, \eqref{J_mu_III_i} gives
\begin{equation}\label{J_mu_III}
|\text{(III)}|\leq o(\|\xi'-\xi\|_{L^2}) + C\|\xi'-\xi\|_{L^2}^2 = o(\|\xi'-\xi\|_{L^2}).
\end{equation}
Finally, for (IV), we use the simplified notation $ \check \E = \E^{\P\times \bar \P \times \hat \P}$ and obtain
\begin{align}
|\text{(IV)}|\leq & \check \E\big[ \partial_\mu g(X^{t,x,\xi}_T,\mu_T,\bar X^{t,\hat \xi, \xi}_T)\cdot \big( \bar X^{t,\hat \xi, \xi'}_T - \bar X^{t,\hat \xi, \xi}_T - \tilde \E[\partial_\mu \bar X^{t,\hat \xi, \xi}_T(\tilde \xi)\cdot(\tilde \xi'-\tilde \xi)]\big)\big]\nonumber\\
\leq & C \Big(\hat \E \big[ \big\|     \bar X^{t,\hat \xi, \xi'}_T - \bar X^{t,\hat \xi, \xi}_T - \tilde \E[\partial_\mu \bar X^{t,\hat \xi, \xi}_T(\tilde \xi)\cdot(\tilde \xi'-\tilde \xi)]\big)\big\|_{L^2}^2\big]\Big)^{1/2}\nonumber\\
\leq & o(\|\xi'-\xi\|_{L^2}).\label{J_mu_IV}
\end{align}
The last inequality comes from the fact that $\partial_\mu X^{t,x,\xi}$ has a bound that is uniform in $x$ (Lemma \ref{app:lemma:X_mubound}). Combining \eqref{J_mu_I}, \eqref{J_mu_II}, \eqref{J_mu_III} and \eqref{J_mu_IV}, we have
\begin{align}\label{o1estimate-Xi}
\Big|\E\Big[ \Xi'-\Xi - \tilde\E[\partial \Xi(\tilde \xi) (\tilde \xi' - \tilde \xi)] \Big] \Big|\leq o(\|\xi'-\xi\|_{L^2}.
\end{align}
By repeating the same arguments as above, we obtain
\begin{align}\label{o1estimate-Fr}
\bigg| \int_t^T \E \big[\sF_r' - \sF_r - \tilde \E[\partial \sF_r(\tilde \xi)\cdot(\tilde \xi'-\tilde\xi)] \big] \dd r\bigg|\leq o(\|\xi'-\xi\|_{L^2}).
\end{align}
\eqref{o1estimate-Xi} and \eqref{o1estimate-Fr} proves \eqref{partial_mu_Jd} due to the definition of L-derivative.
\end{proof}

In the rest of this appendix we summarize results regarding the properties of $X^{t,x,\mu}$, $\partial_x X^{t,x,\mu}$ and $\partial_\mu X^{t,x,\mu}$. For the proofs, refer to \cite{HL16} and \cite{Li18}.

\begin{lemma}\label{app:lemma:Xbound}
Suppose $(X^{t,x,\mu}, X^{t,\xi})$ is a solution to \eqref{Xhat:x}-\eqref{Xhat:xi}. Then there exists a constant $C$ only depending on the bounds in Assumption \ref{assump}-(i) such that
\begin{align*}
&\E\Big[\sup_{t\leq s\leq T}|X^{t,\xi'}-X^{t,\xi}|^2\Big]\leq C\E[|\xi'-\xi|^2],\\
&\E\Big[\sup_{t\leq s\leq T}|X^{t,x',\xi'}-X^{t,x,\xi}|^2\Big]\leq C(|x'-x|^2+\E[|\xi'-\xi|^2]),\\
&\E\Big[\sup_{t\leq s\leq T}|X^{t,\xi}|^2\Big]\leq C(1+\E[|\xi|^2]),\\
&\E\Big[\sup_{t\leq s\leq T}|X^{t,x,\xi}|^2\Big]\leq C(1+|x|^2),\\
&\sup_{t\leq s\leq T}\Wc_2(\P_{X^{t,\xi'}_s},\P_{X^{t,\xi}_s})^2\leq C\E[|\xi'-\xi|^2],
\end{align*}
for any $t\in [0,T)$, $x',x\in \R^d$ and $\xi,\xi' \in L^2_{\Fc_t}$.
\end{lemma}

\begin{lemma}\label{app:lemma:X_xbound}
Suppose $(X^{t,x,\mu}, X^{t,\xi})$ is a solution to \eqref{Xhat:x}-\eqref{Xhat:xi}. Then, for any $s\in [t,T]$, the mapping $x\mapsto X^{t,x,\xi}_s\in L^2_{\Fc_s}$ is Fr\`echet differentiable and has derivative $\partial_x X^{t,x,\xi}_s$. Moreover, there exists a constant $C$ only depending on the bounds in Assumption \ref{assump}-(i) such that
\begin{align*}
&\E\Big[\sup_{t\leq s\leq T}|\partial_x X^{t,x,\xi}_s|^2\Big]\leq C,\\
& \E\Big[\sup_{t\leq s\leq T}|\partial_x X^{t,x,\xi}_s - \partial_x X^{t,x',\xi'}_s|^2\Big]\leq C(|x-x'|^2 + \E[|\xi'-\xi|^2]),
\end{align*}
for any $t\in [0,T]$, $x',x\in \R^d$ and $\xi,\xi'\in L^2_{\Fc_t}$.
\end{lemma}

\begin{lemma}\label{app:lemma:X_mubound}
Suppose $(X^{t,x,\mu}, X^{t,\xi})$ is a solution to \eqref{Xhat:x}-\eqref{Xhat:xi}. Then, for any $s\in [t,T]$, the mapping $\xi \mapsto X^{t,x,\xi}_s$, from $L^2_{\Fc_t}$ to $L^2_{\Fc_s}$, is Fr\`echet differentiable. The derivative is given by $DX^{t,x,\xi}_s(\eta) = \tilde \E[\partial_\mu X^{t,x,\xi}_s(\tilde \xi)\cdot \tilde \eta]$, where $\eta\in L^2_{\Fc_t}$, and $(\tilde \xi,\tilde \eta)$ are independent copies of $(\xi,\eta)$. Moreover, there exists a constant $C$ only depending on the bounds in Assumption \ref{assump}-(i) such that
\begin{align*}
&\E\Big[\sup_{t\leq s\leq T}|\partial_\mu X^{t,x,\xi}_s(y)|^2\Big]\leq C,\\
& \E\Big[\sup_{t\leq s\leq T}|\partial_\mu X^{t,x,\xi}_s (y)- \partial_x X^{t,x',\xi'}_s(y')|^2\Big]\leq C(|x-x'|^2 + |y-y'|^2+\E[|\xi'-\xi|^2]),
\end{align*}
for any $t\in [0,T]$, $x',x\in \R^d$ and $\xi,\xi'\in L^2_{\Fc_t}$.
\end{lemma}

}

\subsection{Proof of Theorem \ref{thm:martingale-decoupledJq}}\label{app:proof-martingale-decoupledJq}

\begin{proof}We first show that if $\widehat J_d (\cdot)= J_d(\cdot; \widehat{\bm \pi}, {\bm \pi})$ and $\widehat q_d(\cdot) = q_d(\cdot; \widehat{\bm \pi}, {\bm \pi})$, then \eqref{thm:martingale-characterization-consistency}-\eqref{equ:martingale_characterization} holds. It is straightforward to see from \eqref{DP-equation-J}, Definition \ref{def:decoupled-q-function}, and \eqref{DPE-of-J-in-terms-q} that \eqref{thm:martingale-characterization-consistency} holds.

By applying It\^o's formula to $\widehat J_d(t', X_{t'}^{R, {\bm h}}, \mu_{t'}^{{\bm h}})$ between $t$ and $s$, $0 \leq t < s \leq T$, we obtain
\begin{align}
& e^{-\beta s} \widehat J_d(s, X_{s}^{R, {\bm h}}, \mu_{s}^{{\bm h}}) - e^{-\beta t} \widehat J_d(t, x, \mu)
+ \int_{t}^s e^{-\beta t'}\biggl( r(t', X_{t'}^{R, {\bm h}}, \mu_{t'}^{{\bm h}}, {a}_{t'}^{{\bm h}}) \nonumber\\
& - \widehat q_d (t', X_{t'}^{R, {\bm h}}, \mu_{t'}^{{\bm h}}, {a}^{{\bm h}}_{t'}, {\bm h})\biggl)\dd t' \nonumber\\
= & \int_t^s e^{-\beta t'} \biggl(\check q_d(t', X_{t'}^{R, {\bm h}}, \mu_{t'}^{{\bm h}},  a^{{\bh}}_{t'}, {\bm h}) -\widehat q_d (t', X_{t'}^{R, {\bm h}}, \mu_{t'}^{{\bm h}},  a^{{\bh}}_{t'}, {\bm h})\biggl)\dd t'  \label{equ:Ito-martingale-process} \\
&  + \int_t^s e^{-\beta t'} \partial_x \widehat J_d (t', X_{t'}^{R, {\bm h}}, \mu_{t'}^{{\bm h}})\trans \sigma(t', X_{t'}^{R, {\bm h}}, \mu_{t'}^{{\bm h}}, a^{{\bm h}}_{t'}) \dd W_{t'} \nonumber\\
& + \int_t^se^{-\beta t'} \Big\{\widehat J_d(t', X_{t'-}^{R, {\bm h}} + \Gamma(t', X_{t'}^{R, {\bm h}}, \mu_{t'}^{{\bm h}}, a^{{\bm h}}), \mu_{t'}^{{\bm h}})  -\widehat J_d(s, X_{t'-}^{R, {\bm h}}, \mu_{t'}^{{\bm h}})\Big\} \dd \widetilde N_{t'},  \nonumber
\end{align}
where $\check q_d$ is defined by
\begin{equation}\label{checkq}
\begin{aligned}
\check q_d(t, x, \mu, a, {\bm h}) &:=  \frac{\partial \widehat J_d}{\partial t}(t, x, \mu) - \beta \widehat J_d(t, x, \mu) + \mathcal{L}^{t, a, \mu}[\widehat J_d(t, \cdot, \mu)](x) + r(t, x, \mu, a)\\
& + \int_{\R^d \times \Ac}\Big[\mathcal{L}^{t, a, \mu}[\frac{\delta \widehat J_d}{\delta \mu}(t, x, \mu)(\cdot)](y)\Big] {\bm h}(a|t, y, \mu) \dd a \mu(\dd y). \nonumber
\end{aligned}
\end{equation}
Taking $\widehat J_d (\cdot)= J_d(\cdot; \widehat{\bm \pi}, {\bm \pi})$ and $\widehat q_d(\cdot) = q_d(\cdot; \widehat{\bm \pi}, {\bm \pi})$ into the above equality \eqref{equ:Ito-martingale-process}, we obtain
\begin{equation*}
\begin{aligned}
& e^{-\beta s} \widehat J_d(s, X_{s}^{R, {\bm h}}, \mu_{s}^{{\bm h}}) - e^{-\beta t} \widehat J_d(t, x, \mu)
+ \int_{t}^s e^{-\beta t'}\biggl( r(t', X_{t'}^{R, {\bm h}}, \mu_{t'}^{{\bm h}}, {a}_{t'}^{{\bm h}}) \\
& - \widehat q_d (t', X_{t'}^{R, {\bm h}}, \mu_{t'}^{{\bm h}}, {a}^{{\bm h}}_{t'}, {\bm h})\biggl)\dd t'\\
= &   \int_t^s e^{-\beta t'} \partial_x \widehat J_d (t', X_{t'}^{R, {\bm h}}, \mu_{t'}^{{\bm h}})\trans \sigma(t', X_{t'}^{R, {\bm h}}, \mu_{t'}^{{\bm h}}, a^{{\bm h}}) \dd W_{t'}\\
& + \int_t^se^{-\beta t'} \Big\{\widehat J_d(t', X_{t'-}^{R, {\bm h}} + \Gamma(t', X_{t'-}^{R, {\bm h}}, \mu_{t'}^{{\bm h}}, a^{{\bm h}})_{t'-}, \mu_{s}^{{\bm h}})  -\widehat J_d(s, X_{t'-}^{R, {\bm h}}, \mu_{t'}^{{\bm h}})\Big\} \dd \widetilde N_{t'},
\end{aligned}
\end{equation*}
which is a $\F^e$-adapted martingale under Assumptions \ref{assump} and \ref{Gassump}.

Conversely, we prove that if \eqref{thm:martingale-characterization-consistency}-\eqref{equ:martingale_characterization} holds, then $\widehat J_d (\cdot)= J(\cdot; \widehat{\bm \pi}, {\bm \pi})$ and $\widehat q_d(\cdot) = q_d(\cdot; \widehat{\bm \pi}, {\bm \pi})$.

In view of \eqref{equ:Ito-martingale-process}, the equation \eqref{equ:martingale_characterization} implies that
\begin{equation*}
\int_t^s e^{-\beta t'} \biggl(\check q_d(t', X_{t'}^{R, {\bm h}}, \mu_{t'}^{{\bm h}},  a^{{\bm h}}, {\bm h}) -\widehat q_d (t', X_{t'}^{R, {\bm h}}, \mu_{t'}^{{\bm h}},  a^{{\bm h}}, {\bm h})\biggl)\dd t'
\end{equation*}
is a local martingale with finite variation and hence zero quadratic variation. Therefore, $\P^e$-almost surely,
\begin{equation*}
\begin{aligned}
&\int_t^s e^{-\beta t'} \biggl(\check q_d(t', X_{t'}^{R, {\bm h}}, \mu_{t'}^{{\bm h}},  a^{{\bm h}}, {\bm h}) -\widehat q_d (t', X_{t'}^{R, {\bm h}}, \mu_{t'}^{{\bm h}},  a^{{\bm h}}, {\bm h})\biggl)\dd t' \\
=& : \int_t^s e^{-\beta t'} f(t', X_{t'}^{R, {\bm h}}, \mu_{t'}^{{\bm h}},  a^{{\bm h}}, {\bm h})\dd t'
\end{aligned}
\end{equation*}
is zero for all $s \in [t, T]$. That is, there exists $\Omega^c \in \Fc$ with $\P^e(\Omega^c) =1$ such that for any $\omega \in \Omega^c$
\begin{equation} \label{proof:martingale_representation_f}
\int_t^s  e^{-\beta t'}f(t', X_{t'}^{R, {\bm h}}, \mu_{t'}^{{\bm h}},  a^{{\bm h}}_{t'}, {\bm h})\dd t' =0,
\end{equation}

The rest is to show $f(t, x, \mu, a, {\bm h}) \equiv 0$, which will be verified by the argument of contradiction.
Suppose that the claim does not hold. By the continuity of $f$, there exists $(t^*, x^*, \mu^*, a^*, {\bm h}^*) \in [0, T] \times \R^d \times \Pc_2(\R^d) \times \Ac \times \Pi$ and $\epsilon > 0$ such that $f(t^*, x^*, \mu^*, a^*, {\bm h}^*) > \epsilon$.
As $f$ is continuous, there exists $\delta >0$ such that if $\max\{|t-t^*|, |x- x^*|, \Wc_2(\mu, \mu^*), |a- a^*|\} < \delta$, then $f(t, x, \mu, a, {\bm h}^*) > \frac{\epsilon}{2}$.
Now let us consider the process $(X^{R, {\bm h}^*}, X^{P, {\bm h}^*})$ with $X^{R, {\bm h}^*}_{t^*-} = x^*, X^{P, {\bm h}^*}_{t^*}\sim \mu^*$.

Let us define
$\tau := \inf\big\{s > t^*: |X^{R, {\bm h}^*}_{s-} - x^*| > \delta \; \mbox{or} \; \Wc_2(\mu_s^{{\bm h}^*}, \mu^*) > \delta \big\} \wedge (t^* + \delta)$.
By the definition of $\tau$ and property of $\{N_s\}_{t^* \leq s \leq T}$, we have that
\begin{equation}\label{proof:prob-tau-positive}
\P(\tau > t^*) > 0.
\end{equation}
It then follows from Lebesgue's Theorem on \eqref{proof:martingale_representation_f} that
\begin{equation}\label{proof:lebsegue-thm}
f(s, X_{s-}^{R, {\bm h}^*}(\omega), \mu_s^{{\bm h}^*}, {a}_s^{{\bm h}^*}(\omega), {\bm h}^*) = 0, \;\; \mbox{for almost every}\; s \in [t^*, \tau(\omega)] \; \mbox{and any}\; \; \omega \in \Omega^c.
\end{equation}
Consider the set $Z(\omega) = \{s \in [t^*, \tau(\omega)]: {a}_s^{{\bm h}^*}(\omega) \in \mathcal{B}_\delta(a^*)\}$ for every $\omega \in \Omega$. It then holds that  $$f(s, X_{s-}^{R, {\bm h}^*}(\omega), \mu_s^{{\bm h}^*}, {a}_s^{{\bm h}^*}(\omega), {\bm h}^*)>\frac{\epsilon}{2} > 0$$ for every $s \in Z(\omega)$, $\omega \in \Omega$. Combining the result with \eqref{proof:lebsegue-thm}, we derive that
\begin{equation*}
\begin{aligned}
0 &= \int_{\Omega}\int_{t^*}^T {\bf 1}(s \in Z(\omega))\dd s \P(d\omega) = \E^e\Big[\int_{t^*}^T {\bf 1}(s \in Z)\dd s\Big]\\
&= \E^e\Big[\int_{t^*}^T {\bf 1}(s \leq \tau(\omega)) {\bf 1} ({ a}_s^{{\bm h}^*}(\omega) \in  \mathcal{B}_\delta(a^*))\dd s\Big]
= \int_{t^*}^T \E^e\big[{\bf 1}(s \leq \tau) {\bf 1} ({a}_s^{{\bm h}^*} \in  \mathcal{B}_\delta(a^*))\big]\dd s\\
& = \int_{t^*}^T \E^e\big[\E^e\big[{\bf 1}(s \leq \tau) {\bf 1} ({a}_s^{{\bm h}^*} \in  \mathcal{B}_\delta(a^*))\big|\Fc_s\big]\big]\dd s \\
&= \int_{t^*}^T \E^e\Big[{\bf 1}(s \leq \tau)\int_{\mathcal{B}_\delta(a^*)} {\bm h}^*(a|s, X^{R, {\bm h}^*}_{s-}, \mu_s^{{\bm h}^*})\dd a\Big]\dd s\\
& \geq \min_{|s-t^*| <\delta, |x-x^*|<\delta, \Wc_2(\mu, \mu^*) < \delta} \int_{\mathcal{B}_\delta(a^*)} {\bm h}^*(a|s, x, \mu)\dd a \int_{t^*}^T \E^e\big[{\bf 1}(s \leq \tau)\big]\dd s\\
&\geq \min_{|s-t^*| <\delta, |x-x^*|<\delta, \Wc_2(\mu, \mu^*) < \delta} \int_{\mathcal{B}_\delta(a^*)}{\bm h}^*(a|s, x, \mu)\dd a \cdot \E^e[(\tau \wedge T) - t^*]\geq 0,
\end{aligned}
\end{equation*}
which implies that $\min_{|s-t^*| <\delta, |x-x^*|<\delta, \Wc_2(\mu, \mu^*) < \delta} \int_{\mathcal{B}_\delta(a^*)}{\bm h}^*(a|s, x, \mu)\dd a =0$ or $\E[(\tau \wedge T) - t^*] =0$. However,  $\min_{|s-t^*| <\delta, |x-x^*|<\delta, \Wc_2(\mu, \mu^*) < \delta} \int_{\mathcal{B}_\delta(a^*)}{\bm h}^*(a|s, x, \mu)\dd a =0$ contradicts with ${\rm supp}({\bm \pi}) = \Ac$ for any admissible policy ${\bm \pi} \in \Pi$, and {$\E^e[(\tau \wedge T) - t^*] =0$ contradicts with $\P(\tau > t^*) >0$} in \eqref{proof:prob-tau-positive}.
We thus conclude that $f(t, x, \mu, a, {\bm h}) =0$ for every $(t, x, \mu, a, {\bm h}) \in [0, T] \times \R^d \times \Pc_2(\R^d) \times \Ac \times \Pi$, i.e. $\check q_d(t, x, \mu, a, {\bm h}) = \widehat q(t, x, \mu, a, {\bm h})$ for every  $(t, x, \mu, a, {\bm h}) \in [0, T] \times \R^d \times \Pc_2(\R^d) \times \Ac \times \Pi$.

Setting ${\bm h} = {\bm \pi}$ and combining with the consistency condition \eqref{thm:martingale-characterization-consistency}, we obtain
\begin{equation*}
\int_{\Ac}\Big(\check q_d(t, x, \mu, a, {\bm \pi}) - \gamma \log \widehat{\bm\pi}(a|t, x, \mu)\Big)\widehat{\bm \pi}(a|t, x, \mu)\dd a = 0.
\end{equation*}
In view of \eqref{checkq}, this implies that $\widehat J_d$ satisfies \eqref{DP-equation-J}.
This, together with the terminal condition $\widehat J(t, x, \mu) = g(x, \mu)$, yields that $\widehat J_d(t, x, \mu) = J(t, x, \mu; \widehat{\bm\pi}, {\bm \pi})$ by virtue of the Feynman-Kac formula. Hence $\widehat J_d(t, x,\mu) =J_d(t, x, \mu; \widehat{\bm\pi}, {\bm \pi}) $. It then follows from $\check q_d = \widehat q_d$ that $\widehat  q_d(\cdot) = q_d(\cdot; \widehat{\bm \pi}, {\bm \pi})$.
\end{proof}
\subsection{Proof of Proposition \ref{cor:MFG-martingale}}\label{subsec:martingale-MFG}
\begin{proof}
The consistency condition \eqref{equ:mfg-consistency} implies that $\widehat q^*_d(t, x, \mu, a, \widehat{\bm \pi}^*) = \gamma \log \widehat {\bm \pi}^*(a|t, x, \mu)$.
Consequently, $\widehat q^*_d(t, x, \mu, a, \widehat{\bm \pi}^*)$ satisfies the constraint in \eqref{thm:martingale-characterization-consistency}. By Theorem \ref{thm:martingale-decoupledJq}, $\widehat J^*_d$ and $\widehat q^*_d$ are the decoupled value function and decoupled Iq function,respectively , associated with $(\widehat{\bm \pi}^*, \widehat{\bm \pi}^*)$.

Furthermore, \eqref{equ:mfg-consistency} characterizes $\widehat{\bm \pi}^*$ as $\widehat {\bm \pi}^* = \Ic^{MFG}(\widehat {\bm \pi}^*, \widehat {\bm \pi}^*)$, where
the mapping  $\Ic^{MFG}$ in \eqref{MFG-PI-map}  can also be expressed in terms of decoupled Iq-function
\begin{equation*}
\Ic^{MFG}(\widehat{\bm\pi}, {\bm \pi})(a|t, x, \mu) = \frac{\exp\Big\{\frac{1}{\gamma} q_d(t, x, \mu, a, {\bm h}; \widehat{\bm \pi}, {\bm \pi})\Big\}}{\int_{\Ac}\exp\Big\{\frac{1}{\gamma} q_d(t, x, \mu, a, {\bm h}; \widehat{\bm \pi}, {\bm \pi})\Big\}\dd a}.
\end{equation*}
{We emphasize that $\widehat{\bm \pi}^*$ is not necessarily unique. For all ${\bm \pi} \in \Pi$, by It\^o's formula
\begin{align}
& e^{-\beta s} \widehat J_d^*(s, X_{s}^{R, {\bm \pi}}, \mu_{s}^{\widehat {\bm \pi}^*}) \nonumber\\
& = e^{-\beta t} \widehat J_d^*(t, x, \mu) + \int_{t}^s e^{-\beta t'}\biggl(\widehat q_d^* (t', X_{t'}^{R, {\bm \pi}}, \mu_{t'}^{{\widehat {\bm \pi}^*}}, {a}^{{\bm \pi}}_{t'}, \widehat {\bm \pi}^*)- r(t', X_{t'}^{R, {\bm \pi}}, \mu_{t'}^{{\widehat {\bm \pi}^*}}, {a}_{t'}^{{\bm \pi}})\biggl)\dd t' \nonumber\\
&  + \int_t^s e^{-\beta t'} \partial_x \widehat J_d (t', X_{t'}^{R, {\bm \pi}}, \mu_{t'}^{\widehat{\bm \pi}})\trans \sigma(t', X_{t'}^{R, {\bm \pi}}, \mu_{t'}^{\widehat{\bm \pi}}, a^{{\bm \pi}}_{t'}) \dd W_{t'} \nonumber\\
& + \int_t^se^{-\beta t'} \Big\{\widehat J_d(t', X_{t'-}^{R, {\bm \pi}} + \Gamma(t', X_{t'}^{R, {\bm \pi}}, \mu_{t'}^{\widehat {\bm \pi}^*}, a^{{\bm \pi}}), \mu_{t'}^{{\bm \pi}})  -\widehat J_d^*(s, X_{t'-}^{R, {\bm \pi}}, \mu_{t'}^{\widehat{\bm \pi}})\Big\} \dd \widetilde N_{t'}.  \nonumber
\end{align}
By setting $s = T$ and taking the expectation, we obtain
\begin{equation*}
\begin{aligned}
e^{-\beta t}\widehat J_d^*(t, x, \mu) = \mathbb{E} \bigg[ & \int_t^T e^{-\beta t'}  \Big(r(t', X_{t'}^{R, {\bm \pi}}, \mu_{t'}^{{\widehat {\bm \pi}^*}}, {a}_{t'}^{{\bm \pi}})- \widehat q_d^* (t', X_{t'}^{R, {\bm \pi}}, \mu_{t'}^{{\widehat {\bm \pi}^*}}, {a}^{{\bm \pi}}_{t'}, \widehat {\bm \pi}^*)\Big) \dd t' \\
& \qquad  + e^{-\beta T } g(X_T^{R, \bm \pi}, \mu_{T}^{{\widehat {\bm \pi}^*}}) \bigg].
\end{aligned}
\end{equation*}
Recalling the variational property of the Gibbs measure, we have that, for all $\pi \in \Pi$,
\begin{equation*}
\begin{aligned}
&\int_{\Ac} \left(\widehat q_d(t, x, \mu, a, \widehat {\bm \pi}^*) - \gamma \log {\bm \pi}(a|t, x, \mu)\right){\bm \pi}(a|t, x, \mu)\dd a\\
\leq  &\int_{\Ac} \left(\widehat q_d(t, x, \mu, a, \widehat {\bm \pi}^*) - \gamma \log \widehat {\bm \pi}^*(a|t, x, \mu)\right)\widehat {\bm \pi}^*(a|t, x, \mu)\dd a =0.
\end{aligned}
\end{equation*}
It then follows that, for all ${\bm \pi} \in \Pi$,
\begin{equation*}
\begin{aligned}
\widehat J_d^*(t, x, \mu) \geq & \E\bigg[\int_t^T e^{-\beta(t' -t)} \Big(r(t', X_{t'}^{R, {\bm \pi}}, \mu_{t'}^{{\widehat {\bm \pi}^*}}, {a}_{t'}^{{\bm \pi}})- \gamma \mathcal{E}_{\bm \pi}(t', X_{t'}^{R, {\bm \pi}}, \mu_{t'}^{{\widehat {\bm \pi}^*}}) \Big)\dd t'  + e^{-\beta(T - t)} g(X_T^{R, \bm \pi}, \mu_T^{\widehat {\bm \pi}^*})\bigg]\\
=& J_d(t, x, \mu; {\bm \pi}, \widehat{\bm \pi}).
\end{aligned}
\end{equation*}
}
This allows us to conclude that $\widehat {\bm \pi}^*$ is a MFE policy of MFG.
%

\end{proof}

\subsection{Proof of Proposition \ref{cor:MFC-martingale}}\label{subsec:martingale-MFC}
\begin{proof} $\mathcal{I}^{MFC}$ defined in \eqref{equ:policy_improvemet_map} can be written in terms of $q_{e}^{MFC}$ that
\begin{equation*}
\mathcal{I}^{MFC}({\bm \pi})(a|t, x, \mu) = \frac{ \exp\left\{\frac{1}{\gamma} q_{e}^{MFC}(t, x, \mu, a; {\bm \pi})\right\} }{ \int_{\Ac}\exp\left\{ \frac{1}{\gamma} q_{e}^{MFC}(t, x, \mu, a; {\bm \pi})\right\}\dd a },
\end{equation*}
and a MFO policy ${\bm \pi}^{MFC, *}$ is unique and satisfies the fixed point equation ${\bm \pi}^{MFC, *}= \mathcal{I}^{MFC}({\bm \pi}^{MFC, *})$.

First, if $\widehat J^*_d = J^{MFC, *}_d$ and $\widehat q^*_d = q^{MFC, *}_d$, by the same arguments in the proof of Theorem \ref{thm:martingale-decoupledJq}, the consistency condition \eqref{equ:mfc-consistency-q}-\eqref{equ:mfc-consistency-pi} and the martingale condition hold.

Conversely, if $\widehat J^*_d$ and $\widehat q^*_d$ satisfy the martingale condition \eqref{equ:martingale_characterization}, then by Theorem \ref{thm:martingale-decoupledJq}, together with $\widehat J^*_d(T, x, \mu) = g(x, \mu)$, and \eqref{equ:mfc-consistency-q}, we obtain that $\widehat J^*_d(t, x, \mu) = J_d(t, x, \mu; \widehat {\bm \pi}^*)$ and $\widehat q_d^*(t, x, \mu, a, {\bm h}) = q_d(t, x, \mu, a, {\bm h}; \widehat{\bm \pi}^*, \widehat{\bm \pi}^*)$. Then by the integral representation in \eqref{equ:mfc-consistency-pi} between $\widehat q_d^*$ and $\widehat q_e^*$,  it holds that $\widehat q_e^*(t, x, \mu, a) = q_e(t, x, \mu, a; \widehat{\bm \pi}^*)$ and thus \eqref{equ:mfc-consistency-pi} becomes
\begin{equation}\label{app:cor-mfc-pi}
\widehat {\bm \pi}^*(a|t, x, \mu) =  \frac{\exp \left\{\frac{1}{\gamma} q_e(t, x, \mu, a; \widehat{\bm \pi}^*)\right\} }{\int_{\Ac} \exp\left\{\frac{1}{\gamma} q_e(t, x, \mu, a; \widehat{\bm \pi}^*)\right\}\dd a} = \Ic^{MFC}(\widehat{\bm \pi}^*)(a|t, x, \mu).
\end{equation}
In view that $\widehat{\bm \pi}^*$ is the fixed point of $\Ic^{MFC}$, we derive from Theorem \ref{prop:PI-MFC} that $\widehat{\bm \pi}^* = {\bm \pi}^{MFC, *}$, $\widehat J^{*}_d = J^{MFC, *}_d$ and $\widehat q^*_d = q^{MFC, *}_d$.
\end{proof}

\subsection{Proof of Proposition \ref{thm:approximation-M-Deltat}}\label{app:proof-thm-approximation-M-Deltat}
\begin{proof}Assume $\beta = 0$ without loss of generality. Consider an auxiliary probability measure space $(\Omega^\Psi, \P^\Psi, \Fc^\Psi)$ independent of $(\Omega^e, \P^e, \Fc^e)$ such that the random chosen of test policies is independent of $(\Omega^e, \P^e, \Fc^e)$.  On $(\Omega^\Psi, \P^\Psi, \Fc^\Psi)$, we generate $M$ i.i.d. random variables $\tilde\psi_1, \ldots, \tilde\psi_M$ from the uniform distribution on $\Psi \subset \R^{L_\Psi}$. Denote by $\E^\Psi[X]$ the expectation of random variable $X$ and denote the empirical measure $\bar\rho^M = \frac{1}{M} \sum_{m=1}^M \delta_{\tilde\psi_m}$. Recall that $\bar M$ is the discrete-time version of the process $M$.\\
Denote
$\Delta M_{t_{i}}^{\tilde\psi_m, \theta, \psi} = M_{t_{i + 1}}^{\tilde\psi_m, \theta, \psi} - M_{t_{i}}^{\tilde\psi_m, \theta, \psi}, \;\;\Delta \bar M_{t_{i}}^{\tilde\psi_m, \theta, \psi}  = \bar M_{t_{i + 1}}^{\tilde\psi_m, \theta, \psi} - \bar M_{t_{i}}^{\tilde\psi_m, \theta, \psi}.$
According to Lemma 1 in \cite{JZ22b}, we need to estimate the difference
\begin{align}
& \frac{1}{M}\sum_{m=1}^M \sum_{i=0}^{K-1}\E^e \Big[\xi_{t_i}^{\tilde\psi_m}\Delta \bar M_{t_{i+1}}^{\tilde\psi_m, \theta, \psi}\Big] - \int_{\Psi} \E^e\Big[\int_0^T \xi_t^{\tilde\psi} \dd M_t^{\tilde\psi, \theta, \psi}\Big] \Uc(\dd \tilde\psi) \nonumber\\
= &  \int_{\Psi}\Big\{\E^e \Big[\sum_{i=0}^{K-1}\xi_{t_i}^{\tilde\psi}\Delta \bar M_{t_{i+1}}^{\tilde\psi, \theta, \psi}
-\int_0^T \xi_t^{\tilde\psi} d M_t^{\tilde\psi, \theta, \psi} \Big]\Big\} \bar\rho_M(\dd \tilde\psi)\nonumber\\
& + \int_0^T \E^e\Big[\xi_t^{\tilde\psi} \dd M_t^{\tilde\psi, \theta, \psi} \Big] \bar\rho_M(\dd \tilde\psi)
-  \int_{\Psi} \E^e\Big[\int_0^T \xi_t^{\tilde\psi} \dd M_t^{\tilde\psi, \theta, \psi}\Big] \Uc(\dd \tilde\psi) \nonumber\\
= &  \int_{\Psi} \E^e\Big[\sum_{i=0}^{K-1}\Big\{\int_{t_i}^{t_{i+1}}\big(\xi_{t_i}^{\tilde\psi}  - \xi_t^{\tilde\psi}\big) \dd M_t^{\tilde\psi, \theta, \psi}
+ \xi_{t_i}^{\tilde\psi}\big( \Delta \bar M_{t_{i+1}}^{\tilde\psi, \theta, \psi}  - \Delta  M_{t_{i}}^{\tilde\psi, \theta, \psi} \big)\Big\}\Big] \bar\rho_M(\dd \tilde\psi) \label{prop:proof-difference-1}\\
& + \int_0^T \E^e\Big[\xi_t^{\tilde\psi} \dd M_t^{\tilde\psi, \theta, \psi} \Big] \bar\rho_M(\dd \tilde\psi)
-  \int_{\Psi} \E^e\Big[\int_0^T \xi_t^{\tilde\psi} \dd M_t^{\tilde\psi, \theta, \psi}\Big] \Uc(\dd \tilde\psi).
 \label{prop:proof-difference-2}
\end{align}
We first estimate \eqref{prop:proof-difference-1} and obtain that
\begin{align}
&  \left|\int_{\Psi} \E^e\left[\sum_{i=0}^{K-1}\left\{\int_{t_i}^{t_{i+1}}\big(\xi_{t_i}^{\tilde\psi}  - \xi_t^{\tilde\psi}\big) \dd M_t^{\tilde\psi, \theta, \psi}
+ \xi_{t_i}^{\tilde\psi}\big( \Delta \bar M_{t_{i+1}}^{\tilde\psi, \theta, \psi}  - \Delta  M_{t_{i}}^{\tilde\psi, \theta, \psi} \big)\right\}\right] \bar\rho_M(\dd \tilde\psi)\right|\nonumber\\
\overset{(a)}{\leq} & \int_{\Psi}\left| \sum_{i=0}^{K- 1}\E^e\left[ \int_{t_i}^{t_{i+1}} \big(\xi_{t_i}^{\tilde\psi}  - \xi_t^{\tilde\psi}\big) L_t^{\tilde\psi, \theta, \psi}\dd t\right]\right|\bar\rho^M(\dd \tilde\psi)
+ \int_{\Psi} \left|\sum_{i=0}^{K -1} \E^e \left[ \xi_{t_i}^{\tilde\psi}\int_{t_i}^{t_{i+1}} \Delta_{i, t'}^{\tilde\psi}\dd t'\right]\right| \nonumber\\
\overset{(b)}{\leq} & \int_{\Psi}\sum_{i=0}^{K- 1}\E^e\left[\int_{t_i}^{t_{i+1}} \big|\xi_{t_i}^{\tilde\psi}  - \xi_t^{\tilde\psi}\big|^2 \dd t\right]^{1/2} \E^e \left[\int_{t_i}^{t_{i+1}}\big|L_t^{\tilde\psi, \theta, \psi}\big|^2 \dd t\right]^{1/2}\bar\rho^M(\dd \tilde\psi) \nonumber\\
& + \int_{\Psi}  \E^e\left[\left(\sum_{i=0}^{K -1}|\xi_{t_i}^{\tilde\psi}|^2\Big)^{1/2} \Big(\sum_{i=0}^{K-1}\big(\int_{t_i}^{t_{i+1}} |\Delta_{i, t'}^{\tilde\psi}|^2\big)\Delta t\right)^{1/2}\right]\bar\rho_M(\dd \tilde\psi)\nonumber\\
 \overset{(c)}{\leq} &\left(\int_{t_i}^{t_{i+1}} C_{\theta, \psi}\big(t_i - t\big)^\alpha \dd t\right)^{1/2}
\int_{\Psi}\sum_{i=0}^{K- 1}\E^e \left[\int_{t_i}^{t_{i+1}}\big|L_t^{\tilde\psi, \theta, \psi}\big|^2 \dd t\right]^{1/2}\bar\rho^M(\dd \tilde\psi) \nonumber\\
& + (\Delta t)^{1/2}  \int_{\Psi} \E^e\left[ \sum_{i=0}^{K - 1}|\xi_{t_i}^{\tilde\psi}|^2\right]^{1/2} \E^e\left[\sum_{i=0}^{K-1}\left(\int_{t_i}^{t_{i+1}} |\Delta_{i, t'}^{\tilde\psi}|^2 \dd t'\right)\right]^{1/2} \bar\rho^M(\dd \tilde\psi)\nonumber\\
\overset{(d)}{\leq} & \sup_{\tilde\psi \in \Psi}C_{\tilde\psi, \theta, \psi} \big(\Delta t\big)^{\frac{\alpha + 1}{2}}  (\frac{T}{\Delta t})^{1/2} \int_{\Psi}\big\|L_t^{\tilde\psi, \theta, \psi}\big\|_2\bar\rho_M(\dd \tilde\psi) + (\Delta t)^{1/2}\int_{\Psi} \|\bar\xi^{\tilde\psi}\|_2 \cdot \|\Delta_i^{\tilde\psi}\|_2 \bar\rho_M(\dd \tilde\psi),
 \label{equ:Delta-approximation}
\end{align}
where in (a), we have used the expression of $M^{\tilde\psi, \theta, \psi}$ and $\bar M^{\tilde\psi, \theta, \psi}$
\begin{equation*}
\begin{aligned}
d M_t^{\tilde\psi, \theta, \psi} =& \left(\check q_d^\theta(t, X_t^{R, {\bm h}^{\tilde\psi}}, \mu_t^{{\bm h}^{\tilde\psi}}, a_t^{{\bm h}^{\tilde\psi}}, {\bm h}^{\tilde\psi}) - q_d^\psi (t, X_t^{R, {\bm h}^{\tilde\psi}}, \mu_t^{{\bm h}^{\tilde\psi}}, a_t^{{\bm h}^{\tilde\psi}}, {\bm h}^{\tilde\psi})\right)\dd t  + \{\cdots\} \dd W_t + \{\cdots\} \dd \widetilde N_t \\
=:& L_t^{\tilde\psi, \theta, \psi}\dd t + \{\cdots\} \dd W_t + \{\cdots\} \dd \widetilde N_t,
\end{aligned}
\end{equation*}
and
\begin{equation*}
\begin{aligned}
\Delta \bar M_{t_{i+1}}^{\tilde\psi, \theta, \psi}  - \Delta  M_{t_{i}}^{\tilde\psi, \theta, \psi} & := \int_{t_i}^{t_{i+1}} \big(r(t_i, X_{t_i}^{R, {\bm h}^{\tilde\psi}}, \mu_{t_i}^{{\bm h}^{\tilde\psi}}) - q_d^\psi(t_i, X_{t_i}^{R, {\bm h}^{\tilde\psi}}, \mu_{t_i}^{{\bm h}^{\tilde\psi}} , {a}^{{\bm h}^{\tilde\psi}}_{t_i}, {\bm h}^{\tilde\psi}) \\
& - r(t', X_{t'}^{R, {\bm h}^{\tilde\psi}}, \mu_{t'}^{{\bm h}^{\tilde\psi}})  + q_d^\psi(t', X_{t'}^{R, {\bm h}^{\tilde\psi}}, \mu_{t'}^{{\bm h}^{\tilde\psi}} , {a}^{{\bm h}^{\tilde\psi}}_{t'}, {\bm h}^{\tilde\psi}) \big)\dd t' : = \int_{t_i}^{t_{i+1}} \Delta^{\tilde\psi}_{i, t'} \dd t',
\end{aligned}
\end{equation*}
(b) is obtained from H\"older's inequality, (c) follows from Assumption \ref{assu:MOC-converegence}, and in (d), we used Assumption \ref{assu:MOC-converegence} and $\E^e[\sup_{0 \leq t \leq T} |X_t^{x, \mu, {\bm h}^{\tilde\psi}}|^2] \leq C(1 + |x|^2 +  (M_2(\mu))^2)$.

We next estimate \eqref{prop:proof-difference-2}. Denote $\E^e\Big[\int_0^T \xi_t^{\tilde\psi} \dd M_t^{\tilde\psi, \theta, \psi}\Big]$ as $f_{\theta, \psi}(\tilde\psi)$ for some deterministic function $f_{\theta, \psi}: \Phi \to \R$. Under Assumption \ref{assu:MOC-converegence}, $f_{\theta, \psi}$ is measurable in $\tilde\psi \in \Psi$ and $f_{\theta, \psi}$ is integrable.
It follows that
\begin{equation}\label{equ:M-approximation}
\frac{1}{M}\sum_{m=1}^M \E^e \left[\int_0^T \xi_t^{\tilde\psi_m} \dd M_t^{\tilde\psi_m, \theta, \psi}\right] = \frac{1}{M} \sum_{m = 1}^M f_{\theta, \psi}(\tilde\psi_m) \to \int_{\Psi} f_{\theta, \psi}(\tilde\psi) \Uc(\dd \tilde\psi), \; \P^{\Psi}-a.s.,
\end{equation}
where in the last equality,  note that $f(\tilde\psi_1), \ldots, f(\tilde\psi_M)$ are i.i.d. random variables,  their arithmetic mean converges almost surely to $\int_{\Psi} f_{\theta, \psi}(\tilde\psi) \Uc(\dd \tilde\psi)$ by strong law of large number.



Combing  \eqref{equ:Delta-approximation}, \eqref{equ:M-approximation} and Lemma 1 in \cite{JZ22b}, and first letting $M \to \infty$ and then letting $\Delta t \to 0$, we arrive at the desired result.
\end{proof}

\section{Concluding Remarks}\label{sec:conc}

{
In this paper, we have established a unified continuous-time q-learning framework for two different but closely related problems: MFG and MFC. Our framework is unified in three aspects.  First, we introduce a unified decoupled formulation (Section \ref{sec:form}) that separates the representative agent's dynamics \eqref{equ:agent-SDE} from the population flow \eqref{equ:population-SDE}, providing a common analytical framework for both MFG and MFC in Section \ref{sec:mfg-mfc}. To the best of our knowledge, this is the first work to introduce a decoupled formulation in the context of continuous-time RL for mean-field problems. Second, within this formulation, we define a unified decoupled Iq-function (Definition \ref{def:decoupled-q-function}) and establish its martingale characterization (Theorem \ref{thm:martingale-decoupledJq}), which serves as a unified policy evaluation rule applicable to both settings. Third, based on this shared theoretical foundation, we develop a unified parametric q-learning algorithm (Algorithm \ref{algo:decoupled-q}) that can be applied to both MFG and MFC.

Most existing continuous-time RL studies focus on either MFG or MFC separately, and do not address both settings within a single framework, see \cite{FGLPS23,weiyu2025, Liangetal2024}. There are only a few exceptional ones such as \cite{AFL2022} that examines the discrete-time RL by treating both problems with alternately updating the standard Q-function and the population distribution under a two-timescale scheme—for MFG, the Q-function is updated at a faster rate than the population distribution, while for MFC the roles are reversed. In sharp contrast, our unified decoupled Iq-function inherently encodes the population dependence through its lifted structure on the enlarged state-action space, thereby eliminating the need for such alternating updates and providing a genuinely single algorithm for both problems.

We conclude with several remarks on technical aspects and future studies. First, regarding the information structure, our theoretical analysis is conducted under a fully observed mean-field formulation, where decoupled value function is defined on the extended state space. In the RL implementation, however, we do not assume access to a simulator of the exact population distribution; instead, we estimate it empirically via the representative agent's own sampled states using update rules \eqref{equ:distribution-update}-\eqref{equ:distribution-update1}. Thus, the algorithm is a sample-based approximation of the fully observed problem. One interesting future extension is to consider the continuous-time RL approach for the partially observed mean-field setting.

Second, concerning the parametric test policy approximation, Proposition \ref{prop:para-martingale-condition} describes the ideal realizable case where the optimal decoupled value function, decoupled Iq-function, and policy lie within the chosen parameterized classes, in which case the learned parameters recover the true solution with no approximation error. In the present paper, we have established the convergence of our algorithm as the number of test policies $M \to \infty$ and the time discretization $\Delta t \to 0$ (Theorem \ref{thm:approximation-M-Deltat}). A natural next step is to analyze how the choice of the parameterized function class affects the convergence error in the general non-realizable case; this would require additional assumptions on the richness of the function classes and the regularity of the optimal solutions, and is left for future study.

Third, regarding the impact of entropy regularization on MFE uniqueness:  in classical MFG theory, the uniqueness of the MFE is typically guaranteed by monotonicity conditions (e.g. in the sense of Lasry-Lions or displacement monotonicity) \cite{Caretal2019, CarD, GMMZ22}. Recent literature in discrete-time mean-field RL \cite{AKS2023, cuikoeppl2021} suggests that the inclusion of entropy regularization may help ensure the uniqueness of the MFE such that the monotonicity condition might be alleviated. The impact of the temperature rate $\gamma$ on the uniqueness of MFE in the continuous-time regularized MFG deserves some future investigations.

Finally, our work provides an efficient learning algorithm for two different but closely related problems with mean-field interactions. It therefore may provide a suitable and effective RL framework to cope with some learning tasks for mixed-type of mean-field systems such as the $\lambda$-interpolated MFG or the p-partial MFG formulated in \cite{CDDL}.
}

\appendix

\section{Proofs of results in the mean-variance example}\label{appendix:MV}

We first consider the MFC problem. To this end, we recall that the goal of the social planner is to maximize the following objective
\begin{equation*}
J(t, \mu; {\bm \pi}) = \E^e[X_T^{\bm \pi}] -\lambda {\rm Var}(X_T^{\bm \pi}) - \gamma  \E^e\biggl[\int_t^T \int_{\R}\log {\bm \pi}(a|s, X_s^{\bm \pi}, \P^e_{X_s^{\bm \pi}}) {\bm \pi}(a|s, X_s^{\bm \pi}, \P^e_{X_s^{\bm \pi}}) \dd a \biggl].
\end{equation*}

Define ${\bm \pi}^{MFC,*}$ to be the MFO policy, and $J^{MFC, *}(t,\mu):= J(t,\mu;{\bm \pi}^{MFC, *})$. We write $J^{MFC, *}$ as $J^*$ in this section for notational simplicity. Denote ${\rm Var}(\mu): = \int_{\R} (x - \bar\mu)^2 \mu(\dd x), \; \bar\mu := \int_{\R} x\mu(\dd x)$. Using Proposition 5.1 in the main paper and the dynamic programming equation of $J^{*}(t,\mu)$, we obtain the following exploratory HJB equation, satisfied by $J^*$:
\begin{equation}
 \frac{\partial J^*}{\partial t}(t, \mu) + \gamma \int_{\R^d} \log \int_{\R} \exp\bigg(\frac{1}{\gamma}\Lc^{t, a, \mu}\bigg(\frac{\delta J^*}{\delta \mu}(t, \mu)(\cdot)\bigg)(v)\bigg) \dd a \mu(\dd v) - \beta J^*(t, \mu) = 0, \label{equ:exploratory_HJB_2}
\end{equation}
We also have the terminal condition $J^*(T, \mu) = \bar \mu - \lambda {\rm Var}(\mu)$. To solve $J^*$, conjecture that $J^*(t, \mu)$ satisfies the quadratic form
\begin{equation}
J^*(t, \mu) = A(t){\rm Var}(\mu) +  C(t)\bar \mu  +D(t).\label{conject-1}
\end{equation}
It then holds that
\begin{equation*}
\begin{aligned}
\frac{\partial J^*}{\partial t}(t, \mu) &= \dot A(t){\rm Var}(\mu) + \dot C(t)\bar \mu  + \dot D(t),\\
\frac{\delta J^*}{\delta \mu}(t, \mu)(v) &= A(t)(v^2 - 2v\bar\mu) + C(t)v,\\
\partial_\mu J^*(t, \mu)(v) &= 2A(t) (v -\bar\mu) + C(t),\quad\quad \partial_v\partial_\mu J^*(t, \mu)(v) = 2A(t).
\end{aligned}
\end{equation*}
Plugging these into $\Lc^{t, a, \mu}[\frac{\delta J^*}{\delta \mu}(t, \mu)(\cdot)](v)$, we obtain that
\begin{equation*}
\begin{aligned}
& \int_{\R}\exp\bigg(\frac{1}{\gamma}\Lc^{t, a, \mu}\bigg(\frac{\delta J^*}{\delta \mu}(t, \mu)(\cdot)\bigg)(v)\bigg)\dd a\\
= & \int_{\R} \exp\bigg(\frac{1}{\gamma}\bigg(ba\partial_\mu J^*(t, \mu)(v) +\frac{\sigma^2a^2}{2} \partial_v\partial_\mu J^*(t, \mu)(v) + \\
&\;\;\eta\bigg(\frac{\delta J^*}{\delta \mu}(t, \mu)(v + \Gamma a) - \frac{\delta J^*}{\delta \mu}(v) - \Gamma a \partial_\mu J^*(t, \mu)(v)\bigg)\bigg)\bigg) \dd a\\
=& \int_{\R} \exp\bigg(\frac{1}{\gamma}\Big(\big(\sigma^2 + \eta \Gamma^2\big) A(t)a^2 + b\big(2A(t)(v -\bar\mu) + C(t)\big)a\Big)\bigg)\dd a\\
=& \sqrt{-\frac{\pi \gamma}{(\sigma^2 + \eta \Gamma^2)A(t)}} \exp\bigg(-\frac{b^2\big(2A(t)(v - \bar\mu) + C(t)\big)^2}{4\gamma (\sigma^2 + \eta \Gamma^2) A(t)}\bigg).
\end{aligned}
\end{equation*}
Using \eqref{equ:exploratory_HJB_2} we get that
\begin{equation*}
\begin{aligned}
& \bigg(\dot A(t) - \frac{b^2}{\sigma^2  +\eta \Gamma^2}A(t)\bigg){\rm Var}(\mu) + \dot C(t) \bar\mu \\
&+ \bigg(\dot D(t) - \frac{b^2}{\sigma^2 + \eta \Gamma^2} \frac{C(t)^2}{4A(t)} + \frac{\gamma}{2} \log \frac{\pi \gamma}{-(\sigma^2 + \eta \Gamma^2) A(t)}\bigg) =0.
\end{aligned}
\end{equation*}
By the terminal conditions $A(T) = -\lambda$, $C(T) = 1$ and $D(T) =0$, we can obtain the explicit solution of the ODEs that
\begin{equation*}
\begin{aligned}
&A(t) =-\lambda \exp\big(\frac{b^2}{\sigma^2 + \eta \Gamma^2}(t-T)\big), \; C(t) =1,\\
& D(t) = \frac{\gamma b^2}{4(\sigma^2 + \eta \Gamma^2)} (t- T)^2 - (t -T) \frac{\gamma}{2} \log \frac{\pi\gamma}{(\sigma^2+\eta \Gamma^2)\lambda} + \frac{1}{4\lambda} e^{-\frac{b^2}{\sigma^2 + \eta \Gamma^2}(t-T)} - \frac{1}{4\lambda}.
\end{aligned}
\end{equation*}
We thus obtain the optimal policy in equ. (6.2) for MFC.

We still need the decoupled Iq function $q^*_d$ to apply our reinforcement learning algorithm. Assume that $J_d$ is in the following quadratic form
\begin{equation}\label{ex:MV-Jd-form}
J_d^*(t, x, \mu) = A_d(t) (x - \bar\mu)^2 + C_d(t) x + D_d(t),
\end{equation}
then its derivatives are given by
\begin{equation*}
\begin{aligned}
\frac{\partial J_d^*}{\partial t} (t, x, \mu)&= \dot A_d(t)(x -\bar\mu)^2 + \dot C_d(t) x + \dot D_d(t), \\
\partial_x J_d^*(t, x, \mu) &= 2A_d(t) (x - \bar\mu) + C_d, \; \partial_{xx} J_d^*(t, x, \mu) = 2A_d(t),\\
\frac{\delta J_d^*}{\delta \mu} (t, x, \mu)(v)&= 2A_d(t)v (\bar\mu - x), \; \partial_\mu J_d^*(t, x, \mu)(v) =2A_d(t)(\bar\mu - x), \; \partial_v\partial_\mu J_d^*(t, x, \mu)(v) =0.
\end{aligned}
\end{equation*}
Then substituting ${\bm \pi}^{MFC, *}$ and \eqref{ex:MV-Jd-form} into equ. (2.7) (with $\hat \bpi=\bpi=\bpi^{MFC, *}$), we have
an explicit solution that
\begin{equation*}
\begin{aligned}
J_d^*(t, x, \mu) &= -\lambda \exp(\frac{b^2}{\sigma^2 + \eta \Gamma^2}(t-T))(x - \bar\mu)^2 + x +\frac{\gamma b^2}{4(\sigma^2 + \eta \Gamma^2)} (t- T)^2 \\
& \;\;\; - (t -T) \frac{\gamma}{2} \log \frac{\pi\gamma}{(\sigma^2 + \eta \Gamma^2)\lambda}+ \frac{1}{4\lambda} \exp\big(-\frac{b^2}{\sigma^2 + \eta \Gamma^2}(t-T)\big) - \frac{1}{4\lambda},
\end{aligned}
\end{equation*}
and from Definition 3.1, it holds that
\begin{equation*}
\begin{aligned}
q_d^*(t, x, \mu, a, {\bm h})
= & -\lambda(\sigma^2 + \eta \Gamma^2) \exp(\frac{b^2}{\sigma^2 + \eta \Gamma^2}(t-T))\Big(a + \frac{b}{\sigma^2 + \eta \Gamma^2}(x - \bar\mu) \\
&-\frac{b}{2\lambda(\sigma^2 + \eta \Gamma^2)} \exp(-\frac{b^2}{\sigma^2 + \eta \Gamma^2}(t-T))\Big)^2 - \frac{b^2}{\sigma^2 + \eta \Gamma^2} (x - \bar\mu) \\
&- \frac{\gamma}{2} \log \frac{\pi \gamma}{(\sigma^2 + \eta \Gamma^2) \lambda} + \frac{\gamma b^2}{2(\sigma^2 + \eta \Gamma^2)}(t -T) \\
&+2\lambda b \exp\big(\frac{b^2}{(\sigma^2 + \eta \Gamma^2)}(t-T)\big) (\bar\mu- x)\E_{\mu, {\bm h}}\big[a^{\bm h}\big].
\end{aligned}
\end{equation*}

Next, we compute $J^{MFG, *}_d$ and $q^{MFG, *}_d$ for the MFG problem and look for the MFE policy. We still assume $J^{MFG, *}_d$ has the form \eqref{ex:MV-Jd-form}. Then,
\begin{equation*}
\begin{aligned}
\Lc^{t,a,\mu}[J_d^{MFG, *}(t,\cdot,h)](x) & =ba \partial_x J^{MFG, *}_d(t,x,\mu) + \frac{\sigma^2 a^2}{2}\partial_{xx}J^{MFG, *}_d(t,x,\mu)  \\
&+ \eta \big(J^{MFG, *}_d(t,x + \Gamma a,\mu)  - J^{MFG, *}_d(t,x,\mu) - \Gamma a \partial_x J^{MFG, *}_d(t,x,\mu)\big)\\
&= ba\big( 2 A_d(t)(x-\bar \mu) + C_d(t)\big) + b^2a^2 A_d(t) + \eta \Gamma^2 a^2 A_d(t).
\end{aligned}
\end{equation*}
By Definition 3.1, we can write
\begin{align}
q^*_d(t,x,\mu,a,\bh) =& \dot{A}_d(t) (x-\bar \mu)^2 + \dot{C}_d(t)x + \dot{D}_d(t) \nonumber\\
&+ ba\big( 2 A_d(t)(x-\bar \mu) + C_d(t)\big) + (b^2+\eta\Gamma^2)A_d(t)a^2  \nonumber\\
 &+ \E^e_{\xi\sim \mu}\bigg[\int_\R q_2(t,x,\mu,a',\xi)\bh(a'|t,\xi,\mu)\dd a'\bigg].\label{exm:mv:q_d}
\end{align}
where
\begin{equation}
q_2(t,x,\mu,a',y) = \Lc^{t,a',\mu}\bigg[\frac{\delta J^*_d}{\delta \mu}(t,x,\mu)(\cdot)\bigg](v)=2ba'A_d(t) (\bar\mu - x).\label{exm:mv:q_2}
\end{equation}
We conclude that the MFE policy ${\bm \pi}^{MFG, *}$ is
\begin{equation}
{\bm \pi}^{MFG, *}(a|t,x,\mu) = \Nc\bigg( -\frac{b}{b^2+\eta\Gamma^2}\bigg( x-\bar\mu + \frac{C_d(t)}{2A_d(t)}\bigg), -\frac{1}{2(b^2+\eta\Gamma^2)A_d(t)}\bigg).\label{append:MV:pi*}
\end{equation}
When $\bh = \bpi^{MFG, *}$, we have
\begin{equation*}
\begin{aligned}
& \E^e_{\xi\sim \mu}\bigg[\int_\R q_2(t,x,\mu,a',\xi)\bh(a'|t,\xi,\mu)\bigg]\\
 =& 2bA_d(t) (\bar\mu- x) \cdot \E^e_{\xi\sim \mu}\bigg[-\frac{b}{b^2+\eta\Gamma^2}\bigg( \xi-\bar\mu + \frac{C_d(t)}{2A_d(t)}\bigg)\bigg] \\
=& \frac{b^2}{b^2+\eta \Gamma^2}C_d(t)(x-\bar\mu).
\end{aligned}
\end{equation*}
We now plug \eqref{exm:mv:q_d} (with $\bh=\bpi^{MFG, *}$) and \eqref{append:MV:pi*} into the master equation (4.3) to obtain that $A_d$,  $C_d$ and $D_d$ satisfy the following equations in the framework of MFG:
\begin{equation*}
\begin{aligned}
& \dot{A}_d(t) = \frac{b^2}{b^2+\eta\Gamma^2} A_d(t), & A_d(T) = -\lambda\\
& \dot{C}_d(t) = 0, & C_d(T) = 1,\\
& \dot{D}_d(t) = \frac{b^2}{4(b^2+\eta\Gamma^2)A_d(t)}-\frac{\gamma}{2}\log \bigg( -\frac{\gamma\pi}{(b^2+\eta\Gamma^2)A_d(t)}\bigg), &D_d(T) = 0.
\end{aligned}
\end{equation*}
It is straightforward to see that these equations lead to the same $J_d^*$ and $q_d^*$ as in the above MFC problem.

{Finally, we verify the uniqueness of the MFE via a fixed-point argument. For a fixed population mean wealth flow $\bar{\mu} = \{\bar{\mu}_s\}_{s \in [0, T]}$, the representative agent's optimal feedback policy $\bm{\pi}^*$ is uniquely determined by solving the associated exploratory HJB equation:
\begin{equation*}
\bm{\pi}^*(a|t, x, \bar{\mu}) = \mathcal{N}\left(-\frac{b}{\sigma^2 + \eta \Gamma^2} \left(x - \bar{\mu}_t - \frac{1}{2\lambda}e^{\frac{b^2}{\sigma^2 + \eta \Gamma^2}(t-T)}\right), \frac{\gamma}{2\lambda (\sigma^2 + \eta \Gamma^2)}e^{\frac{b^2}{\sigma^2 + \eta \Gamma^2}(t-T)}\right).
\end{equation*}
Let $\Phi$ be the mapping that assigns to each $\bar{\mu} \in C([0, T]; \mathbb{R})$ the induced mean wealth flow $m = \{m_s\}_{s \in [0, T]}$ of the representative agent, i.e., $m = \Phi(\bar{\mu})$. Under the optimal policy $\bm{\pi}^*$, $m_s$ satisfies the linear ODE: $\dot{m}_s = -k(m_s - \bar{\mu}_s) + \phi(s)$ with $m_0 = x$, where $k = \frac{b^2}{\sigma^2 + \eta \Gamma^2} > 0$ and $\phi(s)$ is independent of $\bar{\mu}$. For any two flows $\bar{\mu}^{(1)}$ and $\bar{\mu}^{(2)}$, the difference of their images satisfies $m^{(1)}_s - m^{(2)}_s = \int_0^s k e^{-k(s-r)} (\bar{\mu}^{(1)}_r - \bar{\mu}^{(2)}_r) dr$. Defining the norm as $\|h\|_\infty := \sup_{s \in [0, T]} |h_s|$ for any $h \in C([0, T]; \mathbb{R})$, we have
\begin{equation*}
|m^{(1)}_s - m^{(2)}_s| \leq \left(\int_0^s k e^{-k(s-r)} dr \right) \|\bar{\mu}^{(1)} - \bar{\mu}^{(2)}\|_\infty = (1 - e^{-ks}) \|\bar{\mu}^{(1)} - \bar{\mu}^{(2)}\|_\infty.
\end{equation*}
Taking the supremum over $s \in [0, T]$ yields $\|\Phi(\bar{\mu}^{(1)}) - \Phi(\bar{\mu}^{(2)})\|_\infty \leq (1 - e^{-kT}) \|\bar{\mu}^{(1)} - \bar{\mu}^{(2)}\|_\infty$. It shows tthat $\Phi$ is a contraction mapping on $C([0, T]; \mathbb{R})$ and thus possesses a unique fixed point. This ensures the uniqueness of the MFE policy.
}

{
\section{Numerical details of crowd-aversion transport example}

This section gives the parameterizations and numerical implementation
details for the crowd-aversion experiment in Section 7.2 of the main
paper. The MFC run is the primary experiment, and the MFG run uses the same
simulator and numerical hyperparameters as an auxiliary comparison. The state
and action spaces are both $\R^2$. The terminal reward is
\begin{equation*}
g(x,\mu)=-\frac{w_T}{2}|x-x_{\rm tar}|^2
-\frac{w_R}{2}|x|^2,
\end{equation*}
and the running reward is
\begin{equation*}
r(t,x,\mu,a)=-\frac{c_a}{2}|a|^2-w_CK_\ell(x,\mu),
\end{equation*}
where
\[
K_\ell(x,\mu)=\int_{\R^2}
\exp\left(-\frac{|y-x|^2}{2\ell^2}\right)\mu(dy).
\]
Thus the crowd-aversion term is a running cost rather than a terminal cost.

\subsection{Empirical-measure features and value parameterization}

For a particle approximation
$\widehat\mu^N=N^{-1}\sum_{i=1}^N\delta_{X^i}$, let $\bar\mu$ and $v(\mu)$
denote the empirical mean and componentwise empirical variance. We also use
the kernel embeddings
\begin{equation}\label{app:crowd-measure-embedding}
U_j(\mu)=\int_{\R^2}
\exp\left(-\frac{|y-c_j|^2}{2h^2}\right)\mu(dy),
\qquad j=1,\ldots,M_\mu.
\end{equation}
The value-network input consists of
\[
\left(x,\ |x-x_{\rm tar}|^2,\ |x|^2,\ t,\ T-t,\
\bar\mu,\ v(\mu),\ K_\ell(x,\mu),\
U_1(\mu),\ldots,U_{M_\mu}(\mu)\right).
\]
If $F^\theta$ denotes the resulting MLP output, the decoupled value
approximation is
\begin{equation*}
J_d^\theta(t,x,\mu)
=\left(1-\frac{t}{T}\right)F^\theta(t,x,\mu)
+\frac{t}{T}g(x,\mu).
\end{equation*}
Thus, $J_d^\theta(T,x,\mu)=g(x,\mu)$ holds exactly for every $\theta$.

\subsection{Decoupled and essential q-functions}

The q-network uses the local context
\begin{equation*}
z(t,x,\mu)=
\left(1,t,x,x-x_{\rm tar},\bar\mu,v(\mu),K_\ell(x,\mu),
D_\ell(x,\mu)\right),
\end{equation*}
where the local density-gradient feature is
\begin{equation*}
D_\ell(x,\mu)=\nabla_xK_\ell(x,\mu)
=\int_{\R^2}\frac{y-x}{\ell^2}
\exp\left(-\frac{|y-x|^2}{2\ell^2}\right)\mu(dy).
\end{equation*}
An MLP maps $z$ to
\[
\nu^\psi(t,x,\mu)\in\R^2,\qquad
u^\psi(t,x,\mu),s^\psi(t,x,\mu)\in\R^R,\qquad
C^\psi(t,x,\mu)\in\R^{R\times2},
\]
and to three raw precision coefficients
$(\widetilde d_1^\psi,o^\psi,\widetilde d_2^\psi)$. Define
\[
d_i^\psi={\rm softplus}(\widetilde d_i^\psi)-\log 2,
\qquad i=1,2,
\]
and
\begin{equation*}
P^\psi(t,x,\mu)=
\begin{pmatrix}
(d_1^\psi)^2+\rho_{\min} & d_1^\psi o^\psi\\
d_1^\psi o^\psi & (o^\psi)^2+(d_2^\psi)^2+\rho_{\min}
\end{pmatrix}.
\end{equation*}
This construction makes $P^\psi$ positive definite for
$\rho_{\min}>0$. Let
\[
\bar u^\psi(t,\mu)
=\int_{\R^2}u^\psi(t,y,\mu)\mu(\dd y),
\]
and for a generic test policy ${\bm h}$, define
\begin{align*}
&m_{\bm h}(t,y,\mu)
=\int_{\R^2}\alpha\,{\bm h}(d\alpha|t,y,\mu),\\
&\bar s_{\bm h}^\psi(t,\mu)
=\int_{\R^2}\left\{
s^\psi(t,y,\mu)+C^\psi(t,y,\mu)m_{\bm h}(t,y,\mu)
\right\}\mu(\dd y).
\end{align*}
On the normalized reward scale, let
$\widetilde\gamma=\gamma/S_r$. The MFC essential q-function is parameterized as
\begin{equation}\label{app:crowd-essential-q}
\begin{aligned}
q_e^\psi(t,x,\mu,a)
=c^\psi(t,x,\mu)+\widetilde\gamma\bigg\{
&(\nu^\psi(t,x,\mu))^\top a
-\frac12a^\top P^\psi(t,x,\mu)a\\
&+(\bar u^\psi(t,\mu))^\top
\left(s^\psi(t,x,\mu)+C^\psi(t,x,\mu)a\right)
\bigg\}.
\end{aligned}
\end{equation}
Terms independent of $a$ do not affect the Gibbs policy. Consequently, the
MFC target policy is
\begin{equation}\label{app:crowd-gibbs-policy}
{\bm \pi}^\psi(\cdot|t,x,\mu)
=\Nc\left(
(P^\psi)^{-1}
\left\{\nu^\psi+(C^\psi)^\top\bar u^\psi\right\},
\ (P^\psi)^{-1}
\right),
\end{equation}
where all local quantities on the right-hand side are evaluated at
$(t,x,\mu)$. The corresponding decoupled Iq-function is
\begin{equation*}
\begin{aligned}
q_d^\psi(t,x,\mu,a,{\bm h})
=c^\psi(t,x,\mu)+\widetilde\gamma\bigg\{
&(\nu^\psi(t,x,\mu))^\top a
-\frac12a^\top P^\psi(t,x,\mu)a\\
&+(u^\psi(t,x,\mu))^\top\bar s_{\bm h}^\psi(t,\mu)
\bigg\}.
\end{aligned}
\end{equation*}
The scalar normalization is set analytically so that the MFC consistency
condition
\begin{align*}
\int_{\R^2}\left(
q_d^\psi(t,x,\mu,a,{\bm \pi}^\psi)
-\widetilde\gamma\log{\bm \pi}^\psi(a|t,x,\mu)
\right){\bm \pi}^\psi(\dd a|t,x,\mu)=0
\end{align*}
holds. For the auxiliary MFG run, the same local coefficients and precision
matrix are used, but the Gibbs mean is
$(P^\psi)^{-1}\nu^\psi$ and the MFG fixed-point consistency condition is
used in place of the MFC essential-q correction.

\subsection{Test policies and population-distribution update}

For each training episode, we draw $M$ independent shifts
$\delta_m\sim\Nc(0,\sigma_h^2I_2)$ and use
\begin{equation*}
{\bm h}^{\delta_m}(\cdot|t,x,\mu)
=\Nc\left(
m^\psi(t,x,\mu)+\delta_m,\
(P^\psi(t,x,\mu))^{-1}
\right),
\qquad m=1,\ldots,M,
\end{equation*}
where $m^\psi$ denotes the current target-policy mean. The shift changes only
the Gaussian mean and is held fixed over the entire trajectory associated
with that test policy.

For each time point and test policy, the implementation stores an empirical
measure with $N$ atoms. At episode $j$, each stored atom is independently
replaced by the current representative observation with probability
$\rho^j=j^{-0.6}$. Conditional on that observation,
\begin{equation*}
\E\left[
\widehat\mu_{t_k}^{N,j,m}\mid X_{t_k}^{j,m},
\widehat\mu_{t_k}^{N,j-1,m}
\right]
=(1-\rho^j)\widehat\mu_{t_k}^{N,j-1,m}
+\rho^j\delta_{X_{t_k}^{j,m}},
\end{equation*}
which is a particle implementation of the cross-episode distribution update
in the main algorithm.

\subsection{Discrete martingale update and numerical configuration}

For a sampled transition at time $t_k$, define
\begin{equation*}
\begin{aligned}
G_{t_k}^{m}
={}&J_d^\theta(t_{k+1},X_{t_{k+1}}^m,\mu_{t_{k+1}}^m)
-J_d^\theta(t_k,X_{t_k}^m,\mu_{t_k}^m)\\
&+\left(
r_{t_k}^m-\beta J_d^\theta(t_k,X_{t_k}^m,\mu_{t_k}^m)
-q_d^\psi(t_k,X_{t_k}^m,\mu_{t_k}^m,a_{t_k}^m,
{\bm h}^{\delta_m})
\right)\Delta t.
\end{aligned}
\end{equation*}
The implementation uses the parameter gradients of $J_d^\theta$ and
$q_d^\psi$ as test functions. Equivalently, it differentiates the
first-order surrogate
\begin{equation*}
-\frac1{MK}\sum_{m=1}^M\sum_{k=0}^{K-1}
e^{-\beta t_k}\,
{\rm stopgrad}(G_{t_k}^m)
\left\{
J_d^\theta(t_k,X_{t_k}^m,\mu_{t_k}^m)
+q_d^\psi(t_k,X_{t_k}^m,\mu_{t_k}^m,a_{t_k}^m,
{\bm h}^{\delta_m})
\right\}.
\end{equation*}
This produces the same parameter-update direction as the discrete averaged
martingale orthogonality update, up to the optimizer convention and averaging
constants.

The experiment uses
\[
\begin{gathered}
T=5,\quad \Delta t=0.05,\quad \sigma=0.2,\quad \beta=0.1,\quad
x_0=(-2,0),\quad x_{\rm tar}=(2,0),\\
\lambda_J=0.5,\quad \sigma_J=0.1,\quad
c_a=1,\quad w_T=10,\quad w_C=100,\quad w_R=0.1,\quad
\ell=0.8,\quad \gamma=2.
\end{gathered}
\]
The learning quantities are divided by $S_r=1000$, giving
$\widetilde\gamma=0.002$. We take $N=64$ empirical atoms,
$M=20$ test policies, rank $R=16$, and 5000 training episodes. The value and
Q networks have two hidden layers of width 64 with hyperbolic-tangent
activations. The 25 embedding centers in
\eqref{app:crowd-measure-embedding} form a $5\times5$ Cartesian grid on
$[-8,4]^2$, with $h=1$. We use $\rho_{\min}=0.01$,
$\sigma_h=0.05$, initialize the population heads with standard deviation
$0.02$, and initialize the target Gaussian policy with componentwise standard
deviation $\sqrt{2}$. The Gibbs precision is learned. Adam learning rates are
$5\times10^{-4}$ for both the value network and the Q-network, both
multiplied at episode $j$ by $j^{-0.6}$. Validation uses 32 rollouts every 50
episodes with mean-shift perturbation scale $0.1$. For both MFC and MFG problems, the checkpoint with the
highest validation regularized reward is evaluated using 64 independent
randomized-policy rollouts.

}

\ \\
\textbf{Acknowledgement}:{We are grateful to the Associate Editor and two anonymous referees for their helpful comments. X. Wei is supported by National Natural Science Foundation of China grant under no.12201343.  X. Yu is supported by the Hong Kong RGC General Research Fund (GRF) under grant  no. 15211524.}

\end{document}